\documentclass[11pt]{article}

\usepackage[colorlinks=true, pdfstartview=FitV, linkcolor=blue,
citecolor=blue,
urlcolor=blue,
breaklinks=true]{hyperref}
\usepackage{amsmath,amsfonts,amssymb,amsthm,paralist,xspace,a4wide,amscd,array}%,amsxtram,array}
\usepackage[all]{xy}
\usepackage[utf8]{inputenc}
\usepackage{import}
\usepackage[dvipsnames]{xcolor}
\usepackage{enumitem}
\usepackage{mathrsfs}
\usepackage{mathtools}
\usepackage{ bbold }
\usepackage[square]{natbib}
\setcitestyle{numbers}
\DeclareMathAlphabet{\mathpzc}{OT1}{pzc}{m}{it}

\input xy
\xyoption{all}

\usepackage{etoolbox}

\usepackage{enumitem}
\usepackage{stmaryrd}

\newcommand\C{\mathbb{C}}

\newcommand\Z{\mathbb{Z}}
\newcommand\N{\mathbb{N}}

\newcommand\PP{\mathbb{P}}

\newcommand\g{\mathfrak{g}}
\newcommand\h{\mathfrak{h}}
\newcommand\m{\mathfrak{m}}
\newcommand\p{\mathfrak{p}}
\newcommand\lie[1]{\mathfrak{#1}}
\newcommand\frkW{\mathfrak{W}}

\newcommand\admW{{\mathfrak{W}^\textup{adm}}}
\newcommand\frkA{\mathfrak{A}}

\newcommand\M{\mathcal{M}}

\newcommand\ba{\textbf{a}}
\newcommand\wt{\textbf{wt}}
\newcommand\wtvec{\textup{wt}}

\newcommand\bb{\textbf{b}}

\newcommand\bx{\textbf{x}}

\newcommand\A{\mathcal{A}}
\newcommand\U{\mathcal{U}}

\newcommand\cS{\mathcal{S}}
\newcommand\cC{\mathcal{C}}
\newcommand\cE{\mathcal{E}}

\newcommand\cL{\mathcal{L}}

\newcommand\W{\mathcal{W}}
\newcommand\cO{\mathcal{O}}

\newcommand\X{\mathfrak{X}}

\newcommand{\fl}{\textup{fl}}
\newcommand{\mdfg}{\textup{-mod}}
\newcommand{\md}{\textup{-Mod}}

\newcommand{\fmd}{\textup{-fmod}}

\newcommand{\semisim}{\textup{ss}}

\newcommand\irr{\textup{irr}}
\newcommand\wlongest{w_{\textup{long}}}

\newcommand{\del}{\partial}

\newcommand{\cEXT}{\mathcal{EXT}}

\DeclareMathOperator{\Diff}{Diff}

\DeclareMathOperator{\Sing}{Sing}

\DeclareMathOperator{\tr}{tr}
\DeclareMathOperator{\bigtr}{Tr}
\DeclareMathOperator{\ev}{ev}
\DeclareMathOperator{\ad}{ad}

\DeclareMathOperator{\Ann}{Ann}
\DeclareMathOperator{\im}{Im} %Image of a map
\DeclareMathOperator{\Ker}{Ker}
\DeclareMathOperator{\Hom}{Hom}

\DeclareMathOperator{\codim}{codim}

\DeclareMathOperator{\Vspan}{span}
\DeclareMathOperator{\rank}{rank}

\DeclareMathOperator{\Frac}{Frac}
\DeclareMathOperator{\End}{End}

\DeclareMathOperator{\Tor}{Tor}

\DeclareMathOperator{\Aut}{Aut}
\DeclareMathOperator{\Id}{Id}

\DeclareMathOperator{\Spec}{Spec}
\DeclareMathOperator{\Specm}{Specm}
\DeclareMathOperator{\Supp}{Supp}
\DeclareMathOperator{\Suppca}{supp}

\DeclareMathOperator{\ess}{ess}

\newcommand{\setn}[1]{\llbracket #1 \rrbracket}

\DeclareMathOperator{\Res}{Res}

\makeatletter
\def\@endtheorem{\endtrivlist}% NEW
\makeatother

\newtheorem{theo}{Theorem}[section]
\newtheorem*{theo*}{Theorem}
\newtheorem{prop}[theo]{Proposition}
\newtheorem{lem}[theo]{Lemma}
\newtheorem{cor}[theo]{Corollary}

\newtheorem{conj}{Conjecture}[]
\newtheorem*{conj*}{Conjecture}

\theoremstyle{definition}
\newtheorem{defin}[theo]{Definition}
\newtheorem{rem}[theo]{Remark}
\newtheorem*{rem*}{Remark}
\newtheorem*{cor*}{Corollary}
\newtheorem{eg}[theo]{Example}

\newcommand{\keywords}[1]
{
  \small	
  \textbf{\textit{Keywords:}} #1
}

\newcommand{\arxiv}[1]{\href{http://arxiv.org/abs/#1}{\tt arXiv:\nolinkurl{#1}}}

\newcommand{\function}[3]{#1\colon #2 \rightarrow #3}
\newcommand{\sfunction}[2]{#1 \rightarrow #2}
\newcommand{\map}[3]{#1\colon #2 \mapsto #3}

\newcommand{\restr}[2]{\ensuremath{\left.#1\right|_{#2}}}

\newcommand{\nocontentsline}[3]{}
\newcommand{\tocless}[2]{\bgroup\let\addcontentsline=\nocontentsline#1{#2}\egroup}

\newcommand{\rascunho}[1]{}
\newtoggle{lembrete}
\newtoggle{comments}
\newtoggle{details}
\newtoggle{detailsnote}
\newtoggle{prelimnote}
\newtoggle{predetails}

\toggletrue{comments}   % To include comments (default is false)

\toggletrue{details}   % To include details (default is false)

\toggletrue{predetails}   % For including text in case there are no details (default is false)

\iftoggle{lembrete}{%
  \newcommand{\lembrete}[1]{
    \ \\
    {\color{orange}
      \textbf{Lembrar/Pensar depois:} #1
    }
    \\
    }
}{%
  \newcommand{\lembrete}[1]{}
}

\iftoggle{comments}{%
  \newcommand{\comments}[1]{
    \ \\
    {\color{red}
      \textbf{Comments:} #1
    }
    \\
    }
}{%
  \newcommand{\comments}[1]{}
}

\iftoggle{detailsnote}{%
  \newcommand{\detailsnote}[1]{
      \ \\
      {\color{Green}
        \textbf{Details/Notes:} #1
      }
      \\
  }
}{%
  \newcommand{\detailsnote}[1]{}
}

\iftoggle{details}{%
  \newcommand{\details}[1]{
      \ \\
      {\color{OliveGreen}
        \textbf{Details:} #1
      }
      \\
  }
}{%
  \newcommand{\details}[1]{}
}

\iftoggle{predetails}{%
    \iftoggle{details}{
        \newcommand{\predetails}[1]{
            {\color{Gray} #1 }
        }
    }
    {
        \newcommand{\predetails}[1]{
            {#1}
        }
    }
  
}{%
  \newcommand{\predetails}[1]{}
}

\iftoggle{prelimnote}{%
  
}{%
  
}

\newcommand{\Addresses}{{% additional braces for segregating \footnotesize
  \bigskip
  \footnotesize

  (E.M.~Mendonça) 
  \textsc{Institute of Mathematics and Statistics, University of São Paulo, Brazil}\par\nopagebreak
  and \textsc{Institut Camille Jordan, Université Claude Bernanrd Lyon 1, France}\par\nopagebreak
  \textit{E-mail address:} \texttt{edummend@ime.usp.br}
}}

\begin{document}
%{\color{red}{\centering Draft of M.Sc. Thesis\par}}
\title{$\U(\h)$-finite modules and weight modules I: weighting functors, almost-coherent families and category $\frkA^{\irr}$}
\author{Eduardo Monteiro Mendonça}
\date{}

\maketitle \thispagestyle{empty}
\renewcommand{\labelenumi}{(\alph{enumi})}
\numberwithin{equation}{section}
%\numberwithin{equation}{subsection}

\begin{abstract}

This paper builds upon Nilsson's classification of rank one $\U(\h)$-free modules \cite{N16} by extending the analysis to modules without rank restrictions, focusing on the category $\mathfrak{A}$ of $\U(\h)$-finite $\g$-modules. A deeper investigation of the weighting functor $\W$ and its left derived functors, $\W_*$, led to the proof that simple $\U(\h)$-finite modules of infinite dimension are $\U(\h)$-torsion free.  Furthermore, it is shown that these modules are $\U(\h)$-free if they possess non-integral or singular central characters. It is concluded that the existence of $\U(\h)$-torsion-free $\g$-modules is restricted to Lie algebras of types A and C.

The concept of an \emph{almost-coherent family}, which generalizes Mathieu's definition of coherent families \cite{Mat00}, is introduced. It is proved that $\W(M)$, for a $\U(\h)$-torsion-free module $M$, falls within this class of weight modules. Furthermore, a notion of \emph{almost-equivalence} is defined to establish a connection between irreducible semi-simple almost-coherent families and Mathieu's original classification \cite{Mat00}.

Progress is also made in classifying simple modules within the category $\frkA^{\irr}$, which consists of $\U(\h)$-finite modules $M$ with the property that $\W(M)$ is an irreducible almost-coherent family. A complete classification is achieved for type C, with partial classification for $\g$ of type A. Finally, a conjecture is presented asserting that all simple $\lie{sl}(n+1)$-modules in $\frkA^{\irr}$ are isomorphic to simple subquotients of exponential tensor modules defined in \cite{GN22}, and supporting results are proved. In particular, the complete classification of simple $\lie{sl}(3)$-modules in $\frkA^{\irr}$ is obtained.

\end{abstract}


\setcounter{tocdepth}{2}
\tableofcontents
%\tableofcontents

%%%%%%%%%%%%%%%%%%%%%%%%%%%%%%%%%%%%%%%%%%%%%%%%%%%%%%%%%%%%%%%%%%%
% Notes
%%%%%%%%%%%%%%%%%%%%%%%%%%%%%%%%%%%%%%%%%%%%%%%%%%%%%%%%%%%%%%%%%%%
%\section{Notes}
%\import{sections}{notes-24102023}
%\newpage

%%%%%%%%%%%%%%%%%%%%%%%%%%%%%%%%%%%%%%%%%%%%%%%%%%%%%%%%%%%%%%%%%%%
% Introduction
%%%%%%%%%%%%%%%%%%%%%%%%%%%%%%%%%%%%%%%%%%%%%%%%%%%%%%%%%%%%%%%%%%%
%\section*{Introduction}
%%%%%%%%%%%%%%%%%%%%%%%%%%%%%%%%%%%%%%%%%%%%%%%%%%%%%%%%%%%%%%%%%%%
% Introduction
%%%%%%%%%%%%%%%%%%%%%%%%%%%%%%%%%%%%%%%%%%%%%%%%%%%%%%%%%%%%%%%%%%%
\section*{Introduction}

\subsubsection*{The big picture}
In representation theory, a crucial first step in understanding a given category of modules is the complete classification of its simple objects. For the category of all representations of finite-dimensional complex Lie algebras $\g$, this problem is highly challenging and only fully solved for $\g= \mathfrak{sl}(2)$ (see \cite{Blo81}). Consequently, a natural approach is to restrict the problem to specific subclasses of modules, solve it within these restricted categories, and then attempt to generalize the results to larger, more general categories of modules. 

Following this philosophy, the classification of simple representations has been completed for several key categories: finite-dimensional modules (see \cite{Car13}), highest-weight modules (see \cite{Dix96}), the BGG category $\mathcal{O}$ and its parabolic version (see \cite{Hum08}). All these categories have a common feature: they consist of \emph{weight modules}, which are representations where the action of a Cartan subalgebra $\h$ is diagonalizable. The complete classification of simple weight modules with finite-dimensional weight spaces started with Futorny \cite{Fut87} and Fernando \cite{Fer90} and was completed by Mathieu in his breakthrough paper \cite{Mat00}.

Moreover, the motivation for studying representations of Lie algebras has often been influenced by problems in theoretical physics, as well as connections to other branches of modern mathematics, such as quiver theory, D-modules, vertex algebras, and quantum field theory. This influence has been particularly strong in the study of weight modules, making them one of the most well understood classes of representations today.

Other classes of weight modules have also been studied, although non-weight modules they still share some conceptual similarities with weight modules: they are modules where a subalgebra of $\g$ (or of the enveloping algebra $\U(\g)$) acts diagonalizable or locally finitely. It includes the well-known Whittaker modules \cite{Kos78} and Gelfand-Tsetlin modules \cite{DFO94}. This approach has also been extended to representations of noncommutative algebras related to Lie theory, such as the Witt algebras and Weyl algebras.

Recently, a new class of non-weight modules, fundamentally opposite in nature to weight modules, has gained attention. For these modules, it is required that the Cartan subalgebra (or its analog in noncommutative algebras) acts freely. We highlight the work of Nilsson \cite{N15, N16}, who classified $\g$-modules that are free of finite rank over $\U(\h)$. Additionally, Grantcharov and Nguyen \cite{GN22} constructed classes of $\U(\h)$-free $\lie{sl}(n+1)$-modules via modules over Weyl algebras and finite-dimensional $\lie{gl}(n)$-modules. For some infinite-dimensional Lie algebras, modules that are free over Cartan-like subalgebras have been explored (see, e.g. \cite{TZ15, HCL20, FLM24}).

Surprisingly, despite being of a completely opposite nature, weight modules and $\U(\h)$-free modules have interesting connections. These connections are mostly obtained due to the \emph{weighting functor} $\W$ that, as the name suggests, assigns to a $\U(\h)$-free module $M$ a weight module $\W(M)$. This functor was suggested by Mathieu to Nilsson, who introduced it for the first time in \cite{N16} and used it as the main tool to his classification.

As a major goal, this paper aims to explore the interplay between weight modules and $\U(\h)$-free modules, and the role of the weighting functor in linking these two seemingly opposite module categories.

\subsubsection*{Overview of the results}

The notion of a \emph{coherent family} played a crucial role in the classification of simple weight modules presented by Mathieu. Coherent families are ``big'' weight modules whose support is the entire $\h^*$. Although they are certainly not simple modules, coherent families still possess a notion of irreducibility. The classification of irreducible semi-simple coherent families led to the conclusion of the classification of simple weight modules. It also played an important role in Nilsson's work \cite{N16}.

Nilsson proved that the weighting functor maps a $\U(\h)$-free module to a coherent family. Moreover, when $M$ is a $\U(\h)$-free module of rank one, $\W(M)$ is, in fact, an irreducible coherent family. Using Mathieu's classification, Nilsson was able to apply information about irreducible coherent families to complete the classification of $\U(\h)$-free modules of rank one.

In this paper, we extend Nilsson's approach without the restriction on the rank. In particular, to take advantage of working within an abelian category, we expend the study to the category $\frkA$ of $\U(\h)$-finite modules, i.e., $\g$-modules that are finitely generated as $\U(\h)$-modules. In Section~\ref{sec:U(h)-finite-modules}, the first properties of the category $\frkA$ are given and important classes of objects in $\frkA$ are introduced.

The first step in the study was a deeper understanding of the weighting functor $\W$ and, in fact, of their left derived functors $\W_*$, which we also call \emph{weighting functors}. This is discussed in Section~\ref{sec:weighting-functors}. It turns out that $\W_{k>0}$ shares equal properties with $\W_0 = \W$, such as preserving central character and commuting with translation functors, but it also sends a module $M \in \frkA$ to a finite-dimensional module (Proposition~\ref{prop:finite-dim-Wk}) with support exactly where $M$ fails to be locally free (Theorem~\ref{theo:almost-free}). This allowed us to conclude that every simple $\U(\h)$-finite module with either non-integral, or integral and singular central character is, in fact, $\U(\h)$-free (Theorem~\ref{theo:almost-free}), and that, independently of its central character, a simple infinite-dimensional module of $\frkA$ is always $\U(\h)$-torsion-free (Corollary~\ref{cor:simple-tor-free}). Moreover, we proved that, as it happens for coherent families, $\U(\h)$-torsion-free modules only exist for Lie algebras of type A and C (Theorem~\ref{theo:U(h)-free-only-type-A-and-C}).

These results allowed us to investigate, in Section~\ref{sec:weigting-functor-and-a.c.f}, the $\g$-module structure of $\W(M)$ whenever $M \in \frkA$ is $\U(\h)$-torsion-free, leading to the conclusion that it is sufficiently close to being a coherent family. It fails to be a coherent family only in a finite subset of $\h^*$. This leads to the definition of an \emph{almost-coherent family} of degree $d$: a weight module $\M$ such that, for a cofinite $U \subseteq \h^*$, 
\begin{itemize} 
    \item $\dim \M_\lambda = d$ for all $\lambda \in U$; 
    \item for every $a$ in $\A = \U(\g)_0$, the commutator of $\U(\h)$ in $\U(\g)$, the map $\lambda \mapsto \tr \restr{a}{\M_\lambda}$ is polynomial in $U$.
\end{itemize}
Then, Corollary~\ref{Cor:W(M)-U-coherent-family} concludes that $\W(M)$ is an almost-coherent family of degree equal to the rank of $M$ (i.e., $\dim_K M \otimes_{\U(\h)} K$, where $K$ is the field of fractions of $\U(\h)$).
By defining the concept of \emph{almost-equivalence} -- which essentially identifies two semi-simple weight modules that differ only by a finite dimensional module -- one can view Mathieu's classification as a classification of \emph{irreducible almost-coherent families}, up to almost-equivalence.

Then, as a natural follow-up to Nilsson's work, the next step is to classify the simple modules from $\frkA^{\irr}$, the category which consists of $\U(\h)$-finite modules $M$ such that $\W(M)$ is an irreducible almost coherent family. This problem is addressed in Section~\ref{sec:almost-comp-class}. A suitable translation transforms the problem to the classification of simple modules in $\frkA^{\irr}$ with rank one and then use Nilsson's classification.
For $\g$ of type C, the classification was complete:
\begin{theo*}[\ref{theo:class-frkA-typeC}]
    Every simple $\lie{sp}(2n)$-module of $\frkA^{\irr}$ is isomorphic to $T(M)$, where $M$ is a $\U(\h)$-free module of rank one and $T$ is a suitable translation functor.
\end{theo*}
For $\g$ of type A, the complete classification was achieved only for central characters that are not integral and regular:
\begin{theo*}[\ref{theo:class-simple-multipli-one-non-integral-reg}]
    Every simple $\lie{sl}(n+1)$-module of $\frkA^{\irr}$ with non-integral or integral singular central character is isomorphic to an exponential tensor module introduced by Grantcharov and Nguyen in \cite{GN22}.
\end{theo*}

It naturally leads to the following conjecture, stated in the last section:
\begin{conj*}[\ref{conj:exponential-mod}+\ref{conj:full-classification-of-frkA1}]
    Every simple $\lie{sl}(n+1)$-module of $\frkA^{\irr}$ is isomorphic to a simple subquotient of an exponential tensor module.
\end{conj*}
As supporting evidence for the conjecture, it is shown that the candidates for such subquotients indeed lie in $\frkA^{\irr}$ (Proposition~\ref{prop:W(cE_{k,lambda})=EXT(L(w_klambda))}). Furthermore, the conjecture holds when the almost-coherent families in the image of $\W$ can be translated to an almost-coherent family of degree one (Theorem~\ref{theo:class-frkA^(1)-integ-reg-k=1,n}). As a consequence, the classification of simple modules of $\frkA^\irr$ is completed for $\g = \lie{sl}(3)$ (Corollary~\ref{cor:full-class-frkA-irr-sl(3)}).

\subsubsection*{Acknowledgements}
I am very grateful to Olivier Mathieu for proposing the problem and for his valuable support through many discussions that made this work possible.
This study was financed, in part, by the São Paulo Research Foundation (FAPESP) 2020/14313-4 and 2022/05915-6.

%\newpage

\setcounter{section}{0}
%%%%%%%%%%%%%%%%%%%%%%%%%%%%%%%%%%%%%%%%%%%%%%%%%%%%%%%%%%%%%%%%%%%
% Section 0 - Gerneral conventions and definitions
%%%%%%%%%%%%%%%%%%%%%%%%%%%%%%%%%%%%%%%%%%%%%%%%%%%%%%%%%%%%%%%%%%%

\section{General conventions and preliminary results}\label{sec:gen-conv-and-def}

We write by $\Z$ the set of integers and, for every $a \in \Z$ and $\star \in \{\le,\ge,>,<\}$, the set $\Z_{\star a}$ is simply $\{b \in \Z\mid b \star a\}$. And $\N =\Z_{\ge 0}$ denotes the set of non-negative integers. 
The ground field is always assumed to be the complex numbers $\C$ and $V^* = \Hom_\C(V,\C)$ denotes the dual of any vector space $V$.

For an associative algebra $A$, we denote the category of $A$-module by $A\md$. Given a morphism of algebras $\function{\varphi}{A}{B}$, we associate a functor $\function{(-)^\varphi}{B\md}{A\md}$ that assign a $B$-module $M$ to the $A$-module $M^\varphi$ with action defined by $a\cdot m = \varphi(a)m$, for all $a\in A$ and $m \in M$. And, as usual, $\Hom_A(M,N)$ denotes the set of morphism of $A$-modules $M\rightarrow N$. 

Given a category $\cC$ and two full subcategories $\cC_1,\cC_2$ of $\cC$, we define $\cC_1\cap \cC_2$ to be the full subcategory of $\cC$ consisting of objects common to $\cC_1$ and $\cC_2$. 

\subsection{Lie algebras and weight modules}
Throughout this paper, $\g$ will be a given simple finite dimensional Lie algebra (which, after Section~\ref{sec:weigting-functor-and-a.c.f}, is assumed to be of type A or C) and $\h$ a Cartan subalgebra of $\g$. We denote by $\U(\mathfrak{a})$ the universal enveloping algebra of a lie algebra $\mathfrak{a}$ and we define $\mathfrak{a}\mdfg$ to be the category of $\mathfrak{a}$-module that is finitely generated by $\U(\mathfrak{a})$. We also denote by $\mathfrak{a}\fmd$ the category of finite dimensional $\mathfrak{a}$-modules. We fix $n = \dim \h = \rank\g$.

We denote by $\Delta\subseteq\h^*$ the set of roots relative to $\h$ and by $Q = \Z\Delta$ the root lattice. For every $\alpha \in \Delta$ we let $h_\alpha \in\h$ be the corresponding coroot. Then, as a subgroup of $GL(\h^*)$, the Weyl group $W$ of $\g$ is generated by the reflections $s_\alpha \colon \lambda \mapsto \lambda - \lambda(h_\alpha)\alpha$. We fix a basis $\pi \subseteq \Delta$ of simple roots and denote $\Delta^+ = \Z_{\ge 0 }\pi\cap \Delta$ and $\Delta^- = \Z_{\le 0 }\pi\cap \Delta$ the sets of positive and negative roots of $\g$, respectively. We set $P = \{\lambda \in \h^* \mid \lambda(h_\alpha) \in \Z,\, \forall \alpha \in \Delta\}$, $D = \{\lambda \in \h^* \mid \lambda(h_\alpha) \notin \Z_{\le -1},\, \forall \alpha \in \Delta^+\}$ and $P^+ = P\cap D$. An element of $\h^*$ is called \emph{weight} and elements of $P$ and $D$ are called, respectively, \emph{integral weights} and \emph{dominant weights}. Moreover, we set $\rho = \frac{1}{2}\sum_{\alpha \in \Delta^+}$, so that a weight $\lambda$ will be called \emph{regular} if $(\lambda-\rho)(h_\alpha) \ne 0$, for all $\alpha \in \Delta^+$, and called \emph{singular} if it were not regular.

Let $M$ be an $\h$-module. The subspaces $M_{\lambda} = \{m \in M \mid hm = \lambda(h)m,\,\forall h \in \h\}$, for all $\lambda \in \h^*$ are called \emph{weight spaces} of $M$, and its elements \emph{weight vectors}. We set $M[t] = \bigoplus_{\lambda\in t} M_{\lambda}$, for any $t \subseteq \h^*$. Then, for any lie algebra $\mathfrak{a}$ containing $\h$, a \emph{weight $\mathfrak{a}$-module} is an $\mathfrak{a}$-module $M$ such that $M = M[\h^*]$ and all weight spaces are finite dimensional. It is a known fact that every finite dimensional $\mathfrak{a}$-module is in fact a weight $\mathfrak{a}$-module.
When $\mathfrak{a}=\g$, the module $M$ is simply called \emph{weight modules} and the \emph{support} of $M$ is defined to be the set $\Supp M = \{\lambda \in \h^*\mid M_\lambda \ne 0\}$. We denote by $\frkW$ the full subcategory of $\g\mdfg$ consisting of weight modules. We fix the notation $T^* = \h^*/Q$ and called its elements \emph{$Q$-cosets}. Then every weight module $M \in \frkW$ has its support in finitely many $Q$-cosets.

Given a Lie subalgebra $\mathfrak{a} \subseteq \g$ containing $\h$, we see $\mathfrak{a}$ as an $\h$-module by adjoint action. Then $\Delta(\mathfrak{a})$ denotes the set of non-zero weights of the support of the $\h$-module $\mathfrak{a}$. Then, for any parabolic subalgebra $\mathfrak{p}\subseteq \g$ that contains $\h$, we let $\mathfrak{u}^+$ be its nilradical and define $\mathfrak{l} = \h \oplus \bigoplus_{\alpha \in \Delta(\p)\setminus \Delta(\mathfrak{u}^+)}\g_\alpha$. The decomposition $\p = \mathfrak{l}\oplus \mathfrak{u}^+$ is called \emph{Levi decomposition} of $\p$ and $\mathfrak{l}$, defined in this way, will be called \emph{the Levi factor} of $\p$. For any finite dimensional $\p$-module $V$ we define an induced $\g$-module $M_\p(V) \coloneqq \U(\g)\otimes_{\U(\p)}V$ and call it a \emph{parabolic Verma module}. It is a known fact that $M_\p(V)$ is a weight module and when $V$ is simple, it has a unique simple quotient, which will be denoted by $L_\p(V)$. If $\mathfrak{b}'\subseteq \p$ is a Borel subalgebra containing $\h$, and $\mathfrak{n}^+$ is the nilradical of $\mathfrak{b}'$, we say that a weight $\p$-module $M$ is a \emph{$\mathfrak{b}'$-highest weight} module, of highest weight $\lambda \in \h^*$, if there exists a weight vector $v_\lambda \in M_\lambda$ that generates $M$ and $\mathfrak{n}^+v_\lambda = 0$. In fact, every finite dimensional simple $\mathfrak{p}$-module is a \emph{$\mathfrak{b}'$-highest weight module} of weight $\lambda \in \h*$, for some Borel subalgebra $\mathfrak{b}'\subseteq \p$ containing $\h$ and $\mathfrak{u}^+$, and so that we denote it by $L_{\mathfrak{l}}(\lambda)$. When $\p=\mathfrak{b}$ is Borel such that $\Delta(\mathfrak{b}) = \Delta^+$ (and, in particular, its Levi factor is $\h$), we denote the weight module $L_{\mathfrak{b}}(L_\h(\lambda))$ simply by $L(\lambda)$, for all $\lambda \in\h^*$. It is well known that every finite dimensional $\g$-module is isomorphic to $L(\lambda)$ for some integral and dominant weight $\lambda \in P^+$.

\subsection{Central characters and translation functor}

Write $Z(\g)$ for the center of $\U(\g)$. A \emph{central character} of $\g$ is a homomorphism of algebras $\function{\chi}{Z(\g)}{\C}$. We denote by $\X$ the set of all central character of $\g$. For any $\g$-module $M$ and a central character $\chi$, we set the submodule $M^\chi = \{m \in M \mid (z-\chi(z))^N m = 0,\,\forall z \in Z(\g),\, N \gg 0\}$. Then we say that $M$ has \emph{central character} (respectively \emph{generalized central character}) $\chi$ if $(z-\chi(z))M = 0$, for all $z \in Z(\g)$ (respectively if $M = M^\chi$). The category $\g\mdfg^\chi$ is defined to be the full subcategory of $\g\mdfg$ consisting of modules of generalized central character $\chi$.

The known Schur's Lemma for representation of finite dimensional associative algebras admits a version for modules over $\U(\mathfrak{a})$ for a finite dimensional Lie algebra $\mathfrak{a}$ (see \cite{Dix96}). In particular, every simple $\g$-module $M$ admits a central character, denoted by $\chi_M$. Thus, for every $\lambda \in \h^*$, we fix the notation $\chi_\lambda = \chi_{L(\lambda)}$. By Duflo \cite{Duf77}, every central character $\chi$ is equal to $\chi_\lambda$, for some weight $\lambda$. Then we say that a central character $\chi$ is \emph{integral} (respectively \emph{non-integral, regular, singular}) if $\chi = \chi_\lambda$, for some  \emph{integral} (respectively \emph{non-integral, regular, singular}) weight $\lambda$.
The \emph{dot action} (also known as affine action) of $W$ on $\h^*$ is defined by $w \cdot \lambda \coloneqq w(\lambda + \rho) - \rho$, for all $w \in W$ and $\lambda \in \h^*$. By Harish-Chandra theorem (see \cite{Harish51}), two central character $\chi_\lambda$ and $\chi_\mu$ are equal if, and only if, $\lambda \in W\cdot \mu$, i.e., $\lambda$ lies in the dot action orbit of $W$ on $\mu$.

Given a pair of weights $(\lambda,\mu)\in (\h^*)^2$ such that $\lambda - \mu \in P$, there exists a unique integral dominant weight $\widehat{\lambda - \mu} \in P^+$ in the $W$-orbit of $\lambda -\mu$, so that $L(\widehat{\lambda - \mu}) \in \g\fmd$.
Then we may define the endofunctor of $\g\mdfg$:
\[
     T_\mu^\lambda (- ) = ( (-)^{\chi_\mu} \otimes L(\widehat{\lambda - \mu}))^{\chi_\lambda}.
\]
Any functor from the family $\{T_\mu^\lambda\mid \lambda-\mu \in P\}$ is called a \emph{translation functor}.

For any weight $\lambda \in \h^*$, let $W_\lambda = \{w \in W\mid w\cdot \lambda = \lambda\}$ be the stabilizer of $\lambda$ in $W$. Then an integral weight $\lambda$ is regular if, and only if $W_\lambda$ is the identity. 
Then define $\Upsilon'$ to be the set of pairs of weights $(\lambda,\mu)$ satisfying: $\lambda + \rho, \mu + \rho \in D$, $\lambda-\mu \in P$ and $W_\lambda = W_\mu$. % Under this condition, there exists a unique integral dominant weight $\widehat{\lambda - \mu} \in P^+$ in the $W$-orbit of $\lambda -\mu$, so that $L(\widehat{\lambda - \mu}) \in \g\fmd$.

\begin{prop}[{\citealp[Theo.~4.2]{BG80}}]\label{prop:Trans-functor-equiv}
    Let $(\lambda, \mu) \in \Upsilon'$ and let $\mathcal{C}$ be an abelian subcategory of $\g\mdfg$ on which $Z(\g)$ acts locally finitely. If for any finite dimensional module $L$, the functor $-\otimes V$ preserves $\cC$, then
    \[
        T_\mu^\lambda \colon \mathcal{C}\cap\g\mdfg^{\chi_\mu} \rightarrow \mathcal{C}\cap\g\mdfg^{\chi_\lambda}
    \]
    is an equivalence of category with inverse $T_\lambda^\mu$.
\end{prop}

For a weight $\lambda$, define the set $\Delta_{[\lambda]} = \{\alpha \in\Delta \mid\lambda(h_\alpha) \in \Z\}$ and the subgroup $W_{[\lambda]}$ of $W$ generated by reflection $s_\alpha$, with $\alpha$ running in $\Delta_{[\lambda]}$. Then, we define $\lambda^\natural \in \Z\Delta_{[\lambda]}$ uniquely by requiring that $(\lambda^\natural + \rho)(h_\alpha) = (\lambda + \rho)(h_\alpha)$, for all $\alpha \in \Delta_{[\lambda]}$. The following result is due to Jantzen \cite{Jan83}. See also \cite{Hum08} for the proof.

\begin{prop}\label{prop:trans-functor-L(lambda)}
    For every $(\lambda,\mu) \in \Upsilon'$ such that $(\lambda^\natural +\rho)(h_\alpha)$ and $(\mu^\natural +\rho)(h_\alpha)$ have the same sign for all $\alpha \in \Delta_{[\lambda]}=\Delta_{[\mu]}$ (i.e., lies $\lambda^\natural$ and $\mu^\natural$ lies in the same facet), we have the isomorphism
    \[
        T_\mu^\lambda(L(w\cdot \mu)) \simeq L(w\cdot \lambda),
    \]
    for all $w \in W_{[\lambda]}=W_{[\mu]}$.
\end{prop}

For convenience, we define $\Upsilon$ to be the set of pairs of weights $(\lambda,\mu)$ such that $(w\cdot \lambda,w\cdot \mu)\in \Upsilon'$ satisfies the statement of Proposition~\ref{prop:trans-functor-L(lambda)}, for some $w\in W$. Then, for $(\lambda,\mu)\in \Upsilon$ the functor $T_\mu^\lambda$ is defined to be $T_{w\cdot\mu}^{w\cdot \lambda}$, for some $w\in W$ such that $(w\cdot \lambda,w\cdot \mu) \in \Upsilon'$. Since $w\cdot \lambda -w\cdot \mu = w(\lambda - \mu)$, such definition does not depend on the chosen $w$. %Therefore, whenever $(\lambda,\mu) \in \Upstion$, we may apply Propositions~\ref{prop:Trans-functor-equiv} and \ref{prop:trans-functor-L(lambda)}.

Whenever we say that a module $M$ \emph{can be translated to} another module $N$, we mean that there exists a suitable translation functor, which is also an equivalence of categories, that sends $M$ to a module isomorphic to $N$.

\subsection{Automorphism}
An automorphism of Lie algebras $\function{\varphi}{\g}{\g}$ is a linear isomorphism satisfying $\varphi([x,y]) = [\varphi(x),\varphi(y)]$, for all $x,y \in \g$. We denote by $\Aut(\g)$ the set of automorphisms of $\g$. Every automorphism $\varphi \in \Aut(\g)$ extends uniquely to an automorphism of associative algebra $\function{\varphi}{\U(\g)}{\U(\g)}$. We denote by $\Aut(\g,\h)$ the subgroup of $\Aut(\g)$ consisting of all automorphisms that stabilize $\h$.

Given an ad-nilpotent $x \in \g$ (i.e, an element $x \in \g$ such that $\map{\ad(x)}{\g\ni y}{[x,y]\in \g}$ is nilpotent), the linear map 
$\function{\exp(\ad(x))}{\g}{\g}$, $y \mapsto \sum_{i=0}^\infty \frac{1}{i!}\ad(x)^i(y)$, is well-defined and is an automorphism of lie algebras. We call such automorphism \emph{elementary}, and we denote by $\Aut_e(\g) \subseteq \Aut(\g)$ the subgroup consisting of compositions of elementary automorphisms. Define $\Aut_e(\g,\h) = \Aut_e(\g)\cap \Aut(\g,\h)$.

For every $\alpha \in \Delta$, fix $e_\alpha\in\g_\alpha$ and $f_\alpha \in \g_{-\alpha}$ such that $[e_\alpha,f_\alpha] = h_\alpha$. Then we define the automorphism of Lie algebras $\function{\tau}{\g}{\g}$ by $\tau(h_\alpha) = -h_\alpha$, $\tau(e_\alpha) = -f_\alpha$ and $\tau(f_\alpha) = -e_\alpha$, for all $\alpha \in \Delta$. We call $\tau$ \emph{negative transpose}, since by looking to $\g$ in a Lie subalgebra $\lie{gl}(m)$, we may find a suitable basis such that $\tau$ becomes exactly the restriction of the negative transpose of matrices. Clearly $\tau \in \Aut(\g,\h)$.

Seeing $W$ as a Coxeter group, we may consider the length of its elements, so that the longest element $\wlongest \in W$ is well-defined.
\begin{lem}\label{lem:centra-char-with-automorphisms}
    Let $\varphi \in \Aut_e(\g,\h)$ and $\lambda \in \h^*$. Then $(-)^\varphi$ preserves $\g\mdfg^{\chi_\lambda}$ and $(-)^\tau(\g\mdfg^{\chi_\lambda}) = \g\mdfg^{\chi_{-\wlongest\lambda}}$.
\end{lem}
\begin{proof}
    For any $\phi \in \Aut(\g)$, $\chi \in \X$ and $M \in \g\mdfg^\chi$ we have $M^\phi \in \g\mdfg^{\chi\circ \phi}$. Indeed, for all $z \in Z(\g)$, $(z - \chi(\phi(z)))$ acts in $M^\phi$ by $(\phi(z) - \chi(\phi(z)))$ which is a nilpotent action, hence the claim. So, to prove the Lemma, it is enough to show that $\chi_\lambda\circ \varphi = \chi_\lambda$ and $\chi_\lambda\circ \tau = \chi_{-\wlongest\lambda}$, for all $\lambda \in \h^*$.
    
    Consider the triangular decomposition $\g = \mathfrak{n}^- \oplus \h \oplus \mathfrak{n}^+$, where $\mathfrak{n}^\pm = \sum_{\alpha \in \Delta^\pm} \g_\alpha$, one can define the projection $\sfunction{\U(\g) = \U(\h) \oplus (\mathfrak{n}^-\U(\g) + \U(\g)\mathfrak{n}^+)}{\U(\h)}$. The restriction of this projection to $Z(\g)$ defines a commutative algebra homomorphism $\function{\xi}{Z(\g)}{\U(\h)}$ called \emph{Harish-Chandra homomorphism}. Then, for every $\lambda \in \h^*$, we have $\chi_\lambda(z) = \lambda(\xi(z))$, for all $z \in Z(\g)$. 

    Let $z \in Z(\g)$. Note that $\chi_\lambda\circ \varphi(z) = \chi_\lambda(z)$ if and only if $\lambda(\xi(\varphi(z))) = \lambda(\xi(z))$. Seeing elements of $\U(\h)$ as polynomials in $\h^*$, and using the fact that $P^+$ is dense in $\h^*$, we have that $\xi\circ\varphi = \xi$ if and only if $\lambda(\xi(\varphi(z)) = \lambda(\xi(z))$ for all $\lambda \in P^+$. Thus the equality $\chi_\lambda\circ \varphi = \chi_\lambda$ is satisfied for all $\lambda \in \h^*$ if and only if it is satisfied for all $\lambda \in P^+$. Analogously, $\chi_\lambda\circ \tau = \chi_{-\wlongest\lambda}$ for all $\lambda \in \h^*$ if and only if it is satisfied for all $\lambda \in P^+$.

    Let $\lambda \in P^+$ and $\varphi \in \Aut_e(\g,\h)$. Note that $L(\lambda)^\varphi$ is simple and finite dimensional. Furthermore $L(\lambda)^\varphi_\mu = L(\lambda)_{\mu \circ \varphi}$. If $\varphi \in \Aut_e(\g,\h)$ then, by \cite[Ch. VIII, $\S$5, Prop.~2.4]{BourbakiLieAlg}, the map $\Delta \ni \alpha \mapsto \alpha \circ \varphi \in \Delta$ is equal to the action of some element $w \in W$. Therefore $\Supp L(\lambda)^\varphi = w\Supp L(\lambda)$. Since the support of finite dimensional modules are preserved by action of the Weyl group (see \cite[Theo.~1.6]{Hum08}) and simple finite dimensional $\g$-modules are uniquely determined by its support, it follows that $L(\lambda)^\varphi = L(\lambda)$. Then $\chi_\lambda\circ \varphi = \chi_{L(\lambda)^\varphi} = \chi_{L(\lambda)} = \chi_\lambda$. Finally, notice that $\mu \circ \tau = -\mu$, for all $\mu \in \h^*$, which implies $\Supp L(\lambda)^\tau = -\Supp L(\lambda)$. In particular, the heighest weight of $L(\lambda)^\tau$ is $-\lambda'$, where $\lambda'$ is the lowest weight of $L(\lambda)$. But $\lambda' = \wlongest\lambda$. Thus $L(\lambda)^\tau = L(-\wlongest\lambda)$ and thus $\chi_\lambda\circ \tau = \chi_{-\wlongest\lambda}$, which conclude the Lemma.
\end{proof}

\subsection{Admissible modules and coherent families}
A weight module $M\in \frkW$ is called \emph{admissible} if its support is contained in a unique $Q$-coset and its weight spaces have uniformly bounded dimensions, i.e, there exists $d \in \N$ such that $\dim{M_\lambda} \le d$ for all $\lambda \in \h^*$.
The \emph{degree} of an admissible weight module $M$ is the maximal dimension of its weight spaces and denoted by $\deg M$. 
The \emph{essential support} of an admissible weight module $M$ of degree $d$ is $\Supp_{\ess}(M) = \{ \lambda \in \Supp(M) \mid \dim M_\lambda = d\}$.
%We define the category $\admW$ to be the full subcategory of $\frkW$ consisting of finite direct sum of admissible modules.

%We list some interesting results about admissible modules:
\begin{lem}\cite{Mat00}\label{lem:first-prop-admW}
    \begin{enumerate}[label=(\roman*),ref=(\roman*)]
        \item \label{lem:first-prop-admW:finite-lenght} Every admissible module has finite length;
        \item \label{lem:first-prop-admW:type-a-c} There exist infinite dimensional admissible weight modules only for $\g$ of type A and C;
        \item \label{lem:first-prop-admW:supp-ess-dense} Let $M \in \frkW$ be admissible of infinite dimension. Then $\Supp_{\ess} M$ is Zariski dense in $\h^*$;
        %\item \label{lem:first-prop-admW:project-functor-preserv} Projective functors preserves $\admW$.
    \end{enumerate}
\end{lem}
%\begin{proof}
%    The first two items are proven in \cite{Mat00}. The third is a consequence of the first item and \cite[Prop.~3.5]{Mat00}. The last item is straightforward from the fact that
%    $M \otimes V$ is admissible of degree $\deg M \dim V$, for all $M \in \admW$ and $V \in \g\fmd$.
%\end{proof}
\begin{rem}[Admissible weight modules are always finitely generated]\label{rem:adm-implies-g-fg}
    Admissible weight modules are subfamily of object of $\frkW$, where all modules are, by definition, finitely generated over $\U(\g)$. However, this property is a consequence of the uniform bound in the dimension of their weight spaces. In fact, the proof that admissible modules have finite length (in \cite{Mat00}) is independent of the cardinality of the chosen set of generators. And clearly any finite length module is finitely generated. Therefore, any $\g$-module with decomposition in direct sum of weight spaces that have uniformly bounded dimensions is also $\U(\g)$-finitely generated.
\end{rem}

%Simple cuspidal modules are examples of admissible modules. Indeed, if $L$ is a simple cuspidal module, it has the support equal exactly one $Q$-coset and since all weight vectors of $\g$ act bijectively on $L$ it implies that all weight spaces of $L$ must have the same dimension. So, by studying the whole class of admissible modules, O. Mathieu  classified all cuspidal modules in \cite{Mat00}. To mention, they are twisted localizations of admissible highest weight modules. But our main interest for paper lies in Mathieu's approach via coherent families, which we will discuss below.

Define $\A= \U(\g)_0$ to be the commutant of $\U(\h)$ in $\U(\g)$.
\begin{defin}
    A coherent family $\M$ of degree $d = \deg \M$ is a weight module satisfying
    \begin{itemize}
        \item $\dim \M_\lambda = d$, for all $\lambda \in \h^*$;
        \item For all $u \in \A$, the map $\h^* \ni \lambda \mapsto \tr \restr{u}{\M_\lambda} \in \C$ is polynomial.
    \end{itemize}
    $\M$ is called \emph{irreducible} if there is $\lambda \in \h^*$ such that $\M_\lambda$ is a simple $\A$-module.
\end{defin}
 
A $\g$-module is called semi-simple if it is a direct sum of \textbf{finitely or infinitely many} simple $\g$-modules. The main connection between admissible weight modules and coherent family is rooted in the following results:

\begin{prop}\label{prop:EXT(L)}\cite[Prop.~4.8]{Mat00}\cite[Prop.~6.2ii]{Mat00}
    Let $L$ be a simple admissible weight module of infinite dimension, and let $d$ be its degree.
    \begin{enumerate}[label=(\roman*),ref=(\roman*)]
        \item\label{prop:EXT(L):uniq-ss-irr-c-f} There exists a unique semi-simple coherent family $\cEXT(L)$ of degree $d$ that contains $L$;
        \item\label{prop:EXT(L):equality-ss-irr-c-fs} $\cEXT(L)$ is irreducible and for any infinite dimension simple submodule $L' \subseteq \cEXT(L)$ is admissible of degree $d$ and $\cEXT(L) \simeq \cEXT(L')$; and
        \item\label{prop:EXT(L):central-char} The central character of the simple submodules of $\cEXT(L)$ is all the same.
    \end{enumerate}
\end{prop}
As a consequence of Prop.~\ref{prop:EXT(L)}\ref{prop:EXT(L):central-char}, we say that a semi-simple irreducible coherent family $\M$ \emph{has central character} $\chi \in \X$, if it is the central character of its simple submodules, and write $\chi_\M = \chi$.

Conversely to the previous proposition, we have:
\begin{prop}\label{prop:EXT(L)=EXT(L(lambda))}\cite[Prop.~6.2ii]{Mat00}
    Every semi-simple irreducible coherent family is isomorphic to $\cEXT(L)$ for some simple admissible weight module of infinite dimension $L$.
    %Let $\M$ be a semi-simple irreducible coherent family. Then there exist $\lambda \notin P^+$ such that $\M \simeq \cEXT(L(\lambda))$.
\end{prop}

In particular, Lemma~\ref{lem:first-prop-admW}~\ref{lem:first-prop-admW:type-a-c} and Prop.~\ref{prop:EXT(L)=EXT(L(lambda))} implies the existence of irreducible coherent families only for $\g$ of type A or C. 
%A \emph{cuspidal module} is a weight module on each every weight vector of $\g$ acts injectively (and consequently bijectively). In particular, cuspidal modules are special cases of admissible weight modules. Then, by classifying semi-simple irreducible coherent family, O. Mathieu obtained the classification of simple cuspidal modules. This, together with the Fernando's result \cite{Fer90} that shows that every simple weight module is the simple quotient of a parabolic induced module from a simple cuspidal, completed the classification of simple weight modules.
%Next, recall the classification for each of these cases. 

\subsubsection{Irreducible semi-simple coherent families for type A}\label{sec:class-irr-ss-cf-A}
Assume $\g = \lie{sl}(n+1)$. Let $\pi = \{\alpha_1,\dotsc, \alpha_n\}$ be a basis of $\Delta$, indexed in such a way that two consecutive roots are connected, and let $h_1,\dotsc,h_n$ be its corresponding coroots.

Define $\X^\infty(A_n)\subseteq \X$ to be the subset consisting of central characters of \emph{infinite-dimensional admissible weight $\lie{sl}(n+1)$-modules}.
\begin{prop}\label{prop:central-char-adm-sln+1}\cite[Prop.~8.5]{Mat00}
    Let $\chi \in \X^\infty(A_n)$. Then $\chi$ is either integral and regular, integral and singular, or non-integral and regular. Furthermore, there exists a unique $\wt(\chi) \in \h^*$ such that $\chi = \chi_{\wt(\chi)}$ and
    \begin{itemize}
        \item if $\chi$ is integral and regular: $\wt(\chi) \in P^+$;
        \item if $\chi$ is integral and singular: $(\wt(\chi) + \rho)(h_{i_\chi}) = 0$ for a unique $i_\chi \in \{1,\dotsc, n\}$ and $(\wt(\chi) + \rho)(h_i) \in \Z_{> 0}$, for all $i \ne i_\chi$;
        \item and if $\chi$ is non-integral: $\wt(\chi)(h_1) \notin \Z$ and $\wt(\chi)(h_i) \in \Z_{\ge 0}$, for all $2 \le i \le n$.
    \end{itemize}
\end{prop}

Let $s_i = (i, i+1)$, for $i = 1,\dotsc, n$, be the simple transpositions in $W = S_{n+1}$, such that $s_i = s_{\alpha_i}$.
\begin{theo}\label{theo:class-c.f.-sl(n+1)}\cite[Lemma~8.3, Theo.~8.6]{Mat00}
    Let $\chi \in \X^\infty(A_n)$ and $\lambda = \wt(\chi)$ (and $i_\chi$ if exists) as defined above.
    \begin{itemize}
        \item{If $\chi$ is integral regular, the following form a complete list (up to isomorphism) of mutually non-isomorphic irreducible semi-simple coherent families with central character $\chi$:
        \[
            \cEXT(L(s_1 \cdot \lambda)),\quad \cEXT(L(s_2 \cdot \lambda)),\quad \dotsc,\quad \cEXT(L(s_n \cdot \lambda)),
        \]
        Furthermore, let $w \in W$, then $L(w \cdot \lambda)$ is admissible and $\cEXT(L(w \cdot \lambda)) \simeq \cEXT(L(s_i \cdot \lambda))$ if, and only if,
        \[
            w = s_j s_{j-1} \dotsb s_i,\text{ for some } j\ge i; \qquad \text{or} \qquad w=  s_j s_{j+1} \dotsb s_i,\text{ for some } j\le i.
            %\begin{cases}
            %    s_j s_{j-1} \dotsb s_i,\quad &\text{ for some } j\ge i; \text{ or }\\ 
            %    s_j s_{j+1} \dotsb s_i,\quad &\text{ for some } j\le i.
            %\end{cases}
        \]
        }
        \item{If $\chi$ is integral singular, then $\cEXT(L(\lambda))$ is the unique, up to isomorphism, semi-simple irreducible coherent family with central character $\chi$. And for $w \in W$, we have that $L(w \cdot \lambda)$ is admissible if, and only if, 
        \[
            w = s_j s_{j+1} \dotsb s_{i_\chi - 1},\text{ for some } j< i_\chi; \qquad \text{or}\qquad w= s_j s_{j-1} \dotsb s_{i_\chi + 1},\text{ for some } j > i_\chi.
            %\begin{cases}
            %    s_j s_{j+1} \dotsb s_{i_\chi - 1},\quad &\text{ for some } j< i_\chi; \text{ or }\\ 
            %    s_j s_{j-1} \dotsb s_{i_\chi + 1},\quad &\text{ for some } j > i_\chi.
            %\end{cases}
        \]
        }
        \item{If $\chi$ is non-integral, then $\cEXT(L(\lambda))$ is the unique, up to isomorphism, semi-simple irreducible coherent family with central character $\chi$. And for $w \in W$, we have that $L(w \cdot \lambda)$ is admissible if, and only if, 
        \[
            w = s_ks_{k-1}\dotsb s_1, \quad \text{ for some }1 \le k \le n.
        \]}
    \end{itemize}
\end{theo}

We describe next how to compute the degree of irreducible semi-simple coherent families of integral regular central character.
Let $\p \subseteq \g$ be a parabolic subalgebra with root system given by $\Delta(\p) = \Delta(\g) \setminus \{-\beta_1,\dotsc,-\beta_n\}$ where $\beta_i = \alpha_i + \alpha_{i+1} + \dotsb + \alpha_n$, for all $1 \le i \le n$. Let $\mathfrak{u}^+$ be its nilradical and consider the Levi decomposition $\p = \mathfrak{l}\oplus \mathfrak{u}^+$.
Its nilradical $\mathfrak{u}^+$ has root system $\Delta(\mathfrak{u}^+) = \{\beta_1,\dotsc,\beta_n\}$. Every finite dimensional $\p$-module is uniquely determined by a highest weight $\lambda \in P^+_{\mathfrak{l}} = \{\lambda \in \h^*\mid \lambda(h_i) \in \Z_{\ge 0},\, i = 1,\dotsc, n-1\}$.

Define
\begin{equation}\label{eq:def-w_k}
    w_0 = 1\quad\text{and}\quad w_i = s_ns_{n-1}\dotsb s_{n-i+1},\quad 1 \le i \le n.
\end{equation}
\begin{theo}\cite[Lemma~11.2, Theo.~11.4]{Mat00}\label{theo:deg-c.f.}
    For all $\lambda \in P^+$ and $0 \le k \le n$, the weight $w_k\cdot \lambda$ lies in $P^+_{\mathfrak{l}}$ and there exists an exact sequence
        \[
            0 \rightarrow L_\mathfrak{p}(L_{\mathfrak{l}}(w_{k+1} \cdot \lambda)) \rightarrow M_{\mathfrak{p}}(L_{\mathfrak{l}}(w_{k} \cdot \lambda)) \rightarrow L_{\mathfrak{p}}(L_{\mathfrak{l}}(w_{k} \cdot \lambda)) \rightarrow 0
        \]
        where $L_\mathfrak{p}(L_{\mathfrak{l}}(w_{n+1} \cdot \lambda)) = 0$.
    Furthermore, for $1 \le k \le n$, the degree of the coherent family $\cEXT(L(w_k \cdot \lambda))$ is $\sum_{i=0}^{n-k}(-1)^{i} \dim L_{\mathfrak{l}}(w_{k+i} \cdot \lambda)$.
\end{theo}

\subsubsection{Irreducible semi-simple coherent families for type C}
Assume $\g = \lie{sp}(2n)$. Let $\pi = \{\alpha_1,\dotsc, \alpha_n\}$ be a basis of $\Delta$ such that two roots are connected, $\alpha_1,\dotsc, \alpha_{n-1}$ are the short roots and $\alpha_n$ is the longest root. Let $h_1,\dotsc,h_n$ be its corresponding coroots.

Define $\X^\infty(C_n)\subseteq \X$ to be the subset consisting of central characters of \emph{infinite-dimensional admissible weight $\lie{sp}(2n)$-modules}.
%\begin{prop}\label{prop:central-char-adm-sp2n}\cite[Lemmas~9.1 and 9.2]{Mat00}
%    Let $\chi \in \X^\infty(C_n)$. Then there exists a unique $\wt(\chi) \in \h^*$ such that $\chi = \chi_{\wt(\chi)}$ and
%    \begin{itemize}
%        \item $\wt(\chi)(h_i) \in \Z_{\ge 0}$, for all $1 \le i \le n-1$;
%        \item $\wt(\chi)(h_n) \in 1/2 + \Z$ and $\wt(\chi)(h_n) \ge -1/2$; and
%        \item $\wt(\chi)(h_{n-1} + 2h_n) + 2 \in \Z_{\ge0}$.
%    \end{itemize}
%\end{prop}

\begin{theo}\label{theo:class-c.f.-sp2n}\cite[Lemmas~9.1 and 9.2,\,Theo.~9.3]{Mat00}
    Let $\chi \in \X^\infty(C_n)$. Then there exists a unique $\wt(\chi) \in \h^*$ such that $\chi = \chi_{\wt(\chi)}$ and
    \begin{itemize}
        \item $\wt(\chi)(h_i) \in \Z_{\ge 0}$, for all $1 \le i \le n-1$;
        \item $\wt(\chi)(h_n) \in 1/2 + \Z$ and $\wt(\chi)(h_n) \ge -1/2$; and
        \item $\wt(\chi)(h_{n-1} + 2h_n) + 2 \in \Z_{\ge0}$.
    \end{itemize}
    Furthermore, up to isomorphism, the module $\cEXT(L(\wt(\chi)))$ is the unique irreducible semi-simple coherent family with central character $\chi$.
\end{theo}

\subsection{Commutative algebras}
Throughout this paper, we will often see a $\U(\h)$-module with a point of view from theory of commutative algebra and its modules. So basic knowledge of this theory is expected from the reader. We cite \cite{Mats80} as a reference. However, it is worth establishing some notations.

The set of maximal ideals of a commutative ring $R$ is denoted by $\Specm R$. For any maximal ideal $\m \in \Specm R$, the localization of an $R$-module $M$ on the multiplicative set $R\setminus \m$ is denoted by $M_\m$.
For any $\lambda \in \h^*$, $\mathfrak{m}_\lambda \in \Specm \U(\h)$ denotes the maximal ideal equal to the kernel of the morphism of algebra $\function{\overline{\lambda}}{\U(\h)}{\C}$ uniquely extended from $\lambda$. 
By Hilbert's Nullstellensatz, the map 
\[
    \h^* \ni \lambda \mapsto \mathfrak{m}_\lambda \in \Specm \U(\h),
\]
is a homeomorphism of topological spaces (with Zariski topology), and both sets will be constantly identified through this paper. Finally, for any commutative algebra $R$ and an $R$-module $M$, $\Tor_k^R(-,M)$ denotes the $k$-th left derived functor of the functor $- \otimes_R M$.

%From now on, the Cartan sualgebra $\h$ will be fixed. Others Cartan subalgebras will appear in this paper, but we will always denote it differently.

%\newpage

%%%%%%%%%%%%%%%%%%%%%%%%%%%%%%%%%%%%%%%%%%%%%%%%%%%%%%%%%%%%%%%%%%%
% Section - Admissible weight modules
%%%%%%%%%%%%%%%%%%%%%%%%%%%%%%%%%%%%%%%%%%%%%%%%%%%%%%%%%%%%%%%%%%%
%\import{sections/}{admissible-modules}
%\newpage

%%%%%%%%%%%%%%%%%%%%%%%%%%%%%%%%%%%%%%%%%%%%%%%%%%%%%%%%%%%%%%%%%%%
% Section - Category \frkA
%%%%%%%%%%%%%%%%%%%%%%%%%%%%%%%%%%%%%%%%%%%%%%%%%%%%%%%%%%%%%%%%%%%
%%%%%%%%%%%%%%%%%%%%%%%%%%%%%%%%%%%%%%%%%%%%%%%%%%%%%%%%%%%%%%%%%%%
%%%%%%%%%%%%%%%%%%%%%%%%%%%%%%%%%%%%%%%%%%%%%%%%%%%%%%%%%%%%%%%%%%%
% U(h)-finite modules
%%%%%%%%%%%%%%%%%%%%%%%%%%%%%%%%%%%%%%%%%%%%%%%%%%%%%%%%%%%%%%%%%%%
%%%%%%%%%%%%%%%%%%%%%%%%%%%%%%%%%%%%%%%%%%%%%%%%%%%%%%%%%%%%%%%%%%%
\section{$\U(\h)$-finite modules}\label{sec:U(h)-finite-modules}

For any subalgebra $\mathfrak{a} \subseteq \g$, we denote the restriction functor $\sfunction{\U(\g)\md}{\U(\mathfrak{a})\md}$ by $\Res_\mathfrak{a}^\g$.
\begin{defin}
    A $\U(\g)$-finitely generated modules $M$ is said to be a \emph{$\U(\h)$-finite} module if $\Res_\h^\g M$ is a finitely generated $\U(\h)$-module. The category $\frkA$ is defined to be the full subcategory of $\U(\g)\mdfg$ consisting of modules that are $\U(\h)$-finite.
\end{defin}

The category $\frkA$ satisfies the following properties:
\begin{prop}\label{prop:properties-frkA}%\cite[Prop.3]{N15}
    \begin{enumerate}[label=(\roman*),ref=(\roman*)]
        \item\label{prop:properties-frkA:abelian} $\frkA$ is an abelian category;
        \item\label{prop:properties-frkA:serre} $\frkA$ is a Serre subcategory of $\U(\g)\mdfg$, i.e, closed by extensions;
        \item\label{prop:properties-frkA:intersec-weight-mod} $\frkA \cap \frkW = \g\fmd$;
        \item\label{prop:properties-frkA:projective-functors} Translation functors preserves $\frkA$.
    \end{enumerate}
\end{prop}
\begin{proof}
    Clearly, $\frkA$ is closed under taking subquotients and direct sums, since both actions preserves finitely-generation properties. Furthermore, since $\U(\h)$ is noetherian, $\U(\h)\mdfg$ is closed under taking submodules, and thus so is $\frkA$. Therefore, $\frkA$ is abelian.

    The fact that $\frkA$ is a Serre subcategory of $\U(\g)\mdfg$, follows from the fact that $\Res_\h^\g$ is exact and that $\U(\h)\mdfg$ is closed by extensions.

    Modules of finite dimension are finitely generated by any of its basis and finite-dimensional $\g$-modules are weight modules, thus $\g\fmd \subseteq \frkA\cap \frkW$. Now let $M \in \frkA \cap \frkW$, $B = \{v_i \in M \mid i \in I\}$ a $\C$-basis of $M$ consisting of weight vectors and $\{m_1,\dotsc, m_r\} \subseteq M$ be a set of generators over $\U(\h)$. Then, there is a finite subset $J \subseteq I$ such that $m_j \in \bigoplus_{i \in J} \C v_i$, for all $1 \le j \le r$. Then, since $M$ is a weight module, it follows
    \[
        \bigoplus_{i \in I} \C v_i = M = \bigoplus_{j=1}^r \U(\h)m_j \subseteq \bigoplus_{i \in J} \C v_i.
    \]
    Since $B$ is a basis, it follows that $I = J$ and $M \in \g\fmd$. Thus, the third item follows.

    Since $\frkA$ is abelian, to show that translation functors preserve $\h$-finite type $\g$-modules, it is enough to show that $M \otimes V \in \frkA$, for all $M \in \frkA$ and $V \in \g\fmd$. Let $v \in V$ be a weight vector of weight $\lambda \in \h^*$ and $m \in M$. For all $h \in \h$ it follows that $(h - \lambda(h))(m \otimes v) = h m \otimes v$. Thus, if $X \subseteq M$ is a finite set of $\U(\h)$-generators of $M$ and $B$ is a basis of $V$ consisting of weight vectors, then $X\otimes B = \{x \otimes v \mid x \in X, v \in B\}$ is a finite set of $\U(\h)$-generators of $M\otimes V$, hence $M \otimes V \in \frkA$.
\end{proof}

A priory, we do not know yet if $\U(\h)$-finite modules have or not finite length. So it is worth to restrict the class of objects of $\frkA$ for such condition:
\begin{defin}
    The category $\frkA_\fl$ is defined to be the full subcategory of $\U(\g)\mdfg$ consisting of $\U(\h)$-finite modules of finite length.
\end{defin}

Next, we will give some classes of examples of objects in $\frkA$ for $\g = \lie{sl}(n+1)$ and $\g = \lie{sp}(2n)$. In fact, later in the paper, we show that there exist infinite dimensional modules in $\frkA$ only for $\g$ of type A and C.

\subsection{$\U(\h)$-finite modules for type A}

In the following, let $E_{i,j}$ stand as the $(i,j)$-elementary matrix of $\lie{gl}(n+1)$. We see $\lie{sl}(n+1)$ as subalgebra of $\lie{gl}(n+1)$ and fix the following basis 
\[
    \lie{sl}(n+1) = \Vspan\{h_k,E_{i,j}\mid 1 \le i,j \le n+1,\, i \ne j,\, 1 \le k \le n\}
\]
where $h_k= E_{k,k} - E_{k+1,k+1}$. And we set the Cartan subalgebra $\h = \Vspan\{h_k\mid 1 \le k \le n\} \subseteq \lie{sl}(n+1)$.

\subsubsection{Parabolic induced}\label{sec:parabolic-induced}
The first class of examples we are going to give are the \emph{parabolic-induced} ones.
Before defining, we need to construct a class of parabolic subalgebras of $\lie{sl}(n+1)$ and show some interesting properties.

By fixing the canonical basis $\{e_1,\dotsc,e_n,e_{n+1}\}$ of $\C^{n+1}$ we identify $\lie{gl}(n+1)$ with $\lie{gl}(\C^{n+1})$ by assigning $E_{i,j}$ to the linear map $e_k \mapsto \delta_{k,j}e_i$. Then $\lie{sl}(n+1)$ is identified with $\lie{sl}(\C^{n+1})$.
Now, for every subspace $V\subseteq \C^{n+1}$, define the subalgebra
\[
    \p_V = \{f \in \lie{sl}(\C^{n+1}) \mid f(V) \subseteq V\}.
\]
Let $\tau \in \Aut(\lie{sl}(n+1))$ be the negative transpose automorphism associated to the Borel subalgebra $\h + \Vspan_\C\{E_{i,j}\mid 1 \le i < j \le n\}$.
Through the identification $\lie{sl}(n+1)=\lie{sl}(\C^{n+1})$, $\tau$ is in fact the negative transpose of matrices, i.e, $E_{i,j}\mapsto -E_{j,i}$.

\begin{lem}\label{lem:parabolic-subalgs}
    \begin{enumerate}[label=(\roman*),ref=(\roman*)]
        \item\label{lem:parabolic-subalgs:iten1} For every subspace $V \subseteq \C^{n+1}$, $\p_V$ is a parabolic subalgebra of $\lie{sl}(n+1)$;
        \item\label{lem:parabolic-subalgs:iten2} The parabolic subalgebra $\p_V$ has maximal dimension (which is equal to $n^2 + n$) if and only if $\dim V = 1$ or $\codim V = 1$;
        \item\label{lem:parabolic-subalgs:item3} Every parabolic subalgebra of maximal dimension is of the form $\p_V$ for some subspace $V \subseteq \C^{n+1}$ of dimension or co-dimension one.
        \item\label{lem:parabolic-subalgs:iten4} Let $\cS$ be the set of all one-dimensional subspaces of $C^{n+1}$. Then $\{\p_V,\tau(\p_V)\mid V\in \cS\}$ forms a complete list of mutually distinct parabolic subalgebras of $\lie{sl}(n+1)$ of maximal dimension.
    \end{enumerate}
\end{lem}
\begin{proof}
    Fix a basis $\{v_1,\dotsc, v_k\}$ of the subspace $V \subseteq \C^{n+1}$ and complete it to a basis $B = \{v_1,\dotsc,v_k\}\cup\{v_{k+1},\dotsc,v_{n+1}\}$ of $\C^{n+1}$. This defines an algebra automorphism $\function{\phi_B}{\lie{gl}(n+1)}{\lie{gl}(\C^{n+1})}$ such that
    \begin{equation}\label{eq:p_V-matrix}
        \phi_B(\p_V) = 
            \begin{bmatrix} 
                M_{k\times k}(\C) & M_{k \times(n+1-k)}(\C)\\
                0       & M_{(n+1 - k) \times (n+1 - k)}(\C)
            \end{bmatrix} \cap \lie{sl}(\C^{n+1}),
    \end{equation}
    where $M_{r\times s}(\C)$ denotes the algebra of $(r\times s)$-matrices with entry in $\C$. Clearly $\phi_B(\p_V)$ is a parabolic subalgebra, and thus so is $\p_V$, proving \ref{lem:parabolic-subalgs:iten1}. Furthermore, equation \eqref{eq:p_V-matrix} easily implies \ref{lem:parabolic-subalgs:iten2}.

    Now let $\p\subseteq \lie{sl}(n+1)$ be a parabolic subalgebra of maximal dimension. 
    Then it contains a Borel subalgebra $\mathfrak{b}'\subseteq \p$.
    Consider $\pi = \{\alpha_1,\dotsc,\alpha_n\}$ the basis of $\Delta$ such that $h_i$ is the co-root of $\alpha_i$. Then $\mathfrak{b} = \h \oplus \Vspan_\C\{E_{i,j}\mid 1 \le i < j \le n\}$ is the Borel subalgebra of $\lie{sl}(n+1)$ with root system $\Z_{\ge 0 }\pi \cap \Delta$.
    By \cite[Ch. VIII,\S 3, Corollary 3.10]{BourbakiLieAlg}, there is an elementary automorphism $\varphi \in \Aut_e(\g)$ such that $\varphi(\mathfrak{b}') = \mathfrak{b}$. Then by \cite[Ch. VIII,\S 3, Remmark 4]{BourbakiLieAlg}, there exists $\Sigma \subseteq \pi$ such that
    \[
        \varphi(\p) = \mathfrak{b} \oplus \bigoplus_{\beta \in \Z_{\le 0}\Sigma \cap \Delta}\g_\beta.
    \]
    Since $\dim \p = n+ |\Z_{\ge 0}\pi \cap \Delta| + |\Z_{\le 0}\Sigma \cap \Delta|$ and it is maximal, we must have $|\Z_{\le 0}\Sigma \cap \Delta| = n(n-1)/2$. A straightforward computation shows that $\Sigma$ must be equal to $\pi\setminus \{\alpha_1\}$ or $\pi\setminus \{\alpha_n\}$. Then, one can check that
    \[
        \varphi(\p) = 
        \begin{cases}
            \p_{\C e_1}, \quad& \text{if } \Sigma = \pi\setminus \{\alpha_1\}\\
            \p_{\Vspan_\C\{e_1,\dotsc,e_n\}}, \quad& \text{if } \Sigma = \pi\setminus \{\alpha_n\}.
        \end{cases}
    \]
    Moreover, seeing $\C^{n+1}$ as the natural representation of $\lie{sl}(n+1)$, \cite[Ch. VIII,\S 7, Prop.~1.2]{BourbakiLieAlg} implies the existence of a linear isomorphism $\function{F}{\C^{n+1}}{\C^{n+1}}$ such that $\varphi(x)(v) = F^{-1}xFv$, for all $x \in \lie{sl}(n+1),v\in \C^{n+1}$. Therefore,
    \[
        \p = 
        \begin{cases}
            \p_{F(\C e_1)}, \quad& \text{if } \Sigma = \pi\setminus \{\alpha_1\}\\
            \p_{F(\Vspan_\C\{e_1,\dotsc,e_n\})}, \quad& \text{if } \Sigma = \pi\setminus \{\alpha_n\},
        \end{cases}
    \]
    and item \ref{lem:parabolic-subalgs:item3} is proved.

    To prove the last item, we start by claiming that for every subspace $U \subseteq \C^{n+1}$ of co-dimension one, there exists a unique subspace $V\in \cS$ such that $\tau(\p_U) = \p_V$. Indeed, let $\langle,\rangle$ be the complex inner product in $\C^{n+1}$ such that $\langle e_i,e_j\rangle = \delta_{i,j}$, for all $1 \le i, j \le n+1$. Then for every $A \in \lie{gl}(n+1)$, we have $\langle Av,w\rangle = - \langle v,\overline{\tau(A)}w\rangle$, for every $v,w\in \C^{n+1}$ (where $\overline{\cdot}$ denotes the function that conjugates the entry of a matrix). Consider $V = U^\perp$, the subspace of $\C^{n+1}$ orthogonal to $U$, which clearly lies in $\cS$, and let $v\in V$ non-zero. Then for every $A \in \p_U$ and $w \in U$ we have $\langle \tau(A)v,w\rangle = - \langle v,\overline{A}w\rangle = 0$, since $\overline{A} \in \p_U$. This implies that $\tau(A)v \in U^\perp = V$, for all $A \in \p_U$, and thus $\tau(\p_U) = \p_V$. 
    
    It remains to prove the non-redundancy of the list $\{\p_V,\tau(\p_V)\mid V\in \cS\}$. Let $V,U\in \cS$ distinct, $v \in V$, $u \in U$ non-zero, and a basis $\{v,v_2,\dotsc,v_n,u\}$ of $\C^{n+1}$. The linear map 
    \[
        v \mapsto v,\quad v_2 \mapsto -v_2,\quad v_i \mapsto 0 \quad (\text{for }3 \le i \le n),\quad u \mapsto v,
    \]
    lies in $\p_V$ but is not contained in $\p_U$, since does not stabilize $U$. Therefore $\p_U \ne \p_V$. Moreover, for any $V_1,V_2 \in \cS$, $\tau(\p_{V_1}) = \tau(\p_{V_2})$ if and only if $\p_{V_1} = \p_{V_2}$ that is satisfied only when $V_1 = V_2$. Now let $\{v'_1,\dotsc,v'_n\}$ be a basis of $V^\perp$. The linear map
    \[
        v \mapsto v + v'_1,\quad v'_1 \mapsto -v'_1,\quad  v'_i \mapsto 0\quad (\text{for }2 \le i \le n),
    \]
    lies in $\tau(\p_V) = \p_{V^\perp}$ but it is not contained in $\p_V$. Hence $\p_V \ne \tau(\p_V)$. Now write $u = av + a_1v'_1 + \dotsb + a_nv'_n$. Since $U \ne V$, it must have an index $1 \le j \le n$ such that $a_j \ne 0$. Then the linear map
    \[
        v \mapsto v,\quad v'_j \mapsto -v'_j,\quad  v'_i \mapsto 0\quad (\text{for }2 \le i \le n,\, i \ne j),
    \]
    lies in $\tau(\p_V) = \p_{V^\perp}$ but not in $\p_U$ since it does not stabilize $U$. Thus $\p_U \ne \tau(\p_V)$. And item \ref{lem:parabolic-subalgs:iten4} is proved.   
\end{proof}

Let $\PP^n$ denotes the $n$-dimensional projective variety over $\C$ and write by $[b_1\colon\dotsb\colon b_{n+1}] \in \PP^n$ the equivalence class of $(b_1,\dotsc,b_{n+1})\in \C^{n+1}$.
It is easy to see that the assignment
\[
    \PP^n \ni \widetilde{\bb} = [b_1\colon \dotsb \colon b_{n+1}] \mapsto V_{\widetilde{\bb}} = \Vspan_\C\{b_1e_1 + \dotsb + b_{n+1}e_{n+1}\} \in \cS
\]
defines a bijection. 
For every $\bb \in (\C\setminus \{0\})^n$, let $\widetilde{\bb} = [b_1\colon \dotsb \colon b_n \colon 1] \in \PP^n$ and define $\p_{\bb} \coloneqq \p_{V_{\widetilde{\bb}}}$. Then, from \ref{lem:parabolic-subalgs}\ref{lem:parabolic-subalgs:iten4}, we have that $\{\p_\bb,\,\tau(\p_\bb)\mid \bb \in (\C\setminus\{0\})^n\}$ forms a non-redundant list of parabolic subalgebras of $\lie{sl}(n+1)$ of maximal dimension.
\begin{lem}\label{lem:dec-sln-h+p}
    The set 
    \begin{equation*}%\label{eq:lem:dec-sln-h+p}
        \left\{\p_\bb,\,\tau(\p_\bb)\quad\colon \quad \bb \in (\C\setminus\{0\})^n\right\},
    \end{equation*}
    forms a complete and non-redundant list of parabolic subalgebras $\p\subseteq \lie{sl}(n+1)$ satisfying $\lie{sl}(n+1) = \h \oplus \p$.
\end{lem}
\begin{proof}
    Let $\p$ be a parabolic subalgebra such that $\g = \h \oplus \p$. Then $\p$ has co-dimension $n$ and by Lemma~\ref{lem:parabolic-subalgs} there exists $\overline{\bb} = [b_1\colon \dotsb \colon b_{n+1}] \in \PP^n$ such that $\p = \p_{V_{\overline{\bb}}}$ or $\p = \tau(\p_{V_{\overline{\bb}}})$.

    Suppose $\p = \p_{V_{\overline{\bb}}}$ and let $v_{\overline{\bb}} = \sum_{i = 1}^{n+1}b_ie_i$. Then, since $\h \cap \p = 0$, for every $a =(a_1,\dotsc,a_n)\in \C^n$ we must have $\sum_{i=1}^{n}a_ih_i(v_{\overline{\bb}}) \notin \C v_{\overline{\bb}}$. But, setting $a_0 = a_{n+1} = 0$, we have $\sum_{i=1}^{n}a_ih_i(v_{\overline{\bb}}) = \sum_{j = 1}^{n+1}(a_j-a_{j-1})b_je_j$. Thus for every $a$, the following must be satisfied
    \[
        [a_1b_1\colon (a_2-a_1)b_2 \colon \dotsb \colon (a_i - a_{i-1})b_i\colon \dotsb \colon (a_n-a_{n-1})b_n \colon -a_nb_{n+1}] \ne {\overline{\bb}}.
    \]
    Now, suppose $b_i = 0$ for some $i \in \{1,\dotsc,n+1\}$, one can find $a \in \C^n$ such that $a_j - a_{j-1} = 1$ for all $1 \le j \le n+1$, $j \ne i$. And in particular $[a_1b_1\colon \dotsb \colon (a_i - a_{i-1})b_i\colon \dotsb \colon -a_nb_{n+1}] = \overline \bb$, which is a contradiction. Hence $\overline{\bb} = [b_1/b_{n+1}\colon \dotsb \colon b_{n}/b_{n+1}\colon 1]$ and $\p$ is in the list described in the lemma.

    Suppose $\p = \tau(\p_{V_{\overline{\bb}}})$. Then $0 = \tau(\h \cap \tau(\p_{V_{\overline{\bb}}})) = \tau(\h)\cap \p_{V_{\overline{\bb}}}$, since $\tau$ is an involution. But $\tau(h_i) = - h_i$, for all $1 \le i \le n$, and thus $\tau(\h) = \h$. Therefore, we have $\h \cap \p_{V_{\overline{\bb}}} = 0$, which implies that $\tau(\p) = \p_{V_{\overline{\bb}}}$ is in the list described in the lemma's statement. Thus so is $\p$.

    The non-redundancy of the list follows from Lemma~\ref{lem:parabolic-subalgs}\ref{lem:parabolic-subalgs:iten4}.
\end{proof}

\begin{eg}[Parabolic induced modules in $\frkA$]\label{eg:parabolic-induced}
    Let $\mathfrak{q}$ be any parabolic subalgebra such that $\g = \lie{sl}(n+1) = \h \oplus \mathfrak{q}$, which are described in Lemma~\ref{lem:dec-sln-h+p} and let $V$ be a finite dimensional $\mathfrak{q}$-module. Consider the parabolic Verma module $M_{\mathfrak{q}}(V)$.
    By the decomposition of $\g$, we obtain $\U(\g) \simeq \U(\h)\otimes_\C\U(\mathfrak{q})$. Then, as $\U(\h)$-module, the module $M_{\mathfrak{q}}(V)$ is isomorphic to $\U(\h) \otimes_\C \widetilde V$, where $\widetilde V = V$ as vector space with trivial $\h$-action. Therefore $M_{\mathfrak{q}}(V) \simeq \U(\h)^{\oplus \dim V}$, and clearly, $M_{\mathfrak{q}}(V) \in \frkA$.
\end{eg}

Let $\p\subseteq \lie{sl}(n+1) = \g$ be a parabolic subalgebra. Then it contains a Borel subalgebra $\mathfrak{b}'\subseteq \p$, which contains a Cartan subalgebra $\h' \subseteq\mathfrak{b}'$. This determines a basis of simple roots $\pi' =\{\alpha_1', \dotsc, \alpha_n'\}$ of root system $\Delta'$ associated to $\h'$ such that $\Delta'(\mathfrak{b}') = \Delta'^+$. Consider then the Levi decomposition $\p = \mathfrak{l} \oplus \mathfrak{u}^+$
associated with $\pi'$ (as defined in Section~\ref{sec:gen-conv-and-def} - so that $\mathfrak{l}$ contains $\h'$), and define the opposite nilradical $\mathfrak{u}^-$ to be Lie subalgebra $\mathfrak{u}^- = \oplus_{\beta' \in \Delta'(\mathfrak{u}^-)} \g_{-\beta'}$. Then we have a triangular decomposition $\g = \mathfrak{u}^- \oplus\mathfrak{l}\oplus \mathfrak{u}^+$.
With this setting, we consider the \emph{parabolic BGG category $\cO^{\p}$} which is the full subcategory of $\U(\g)\mdfg$ consisting of finitely generated $\g$-modules $M$ satisfying
\begin{itemize}
    \item $\Res^\g_{\mathfrak{l}}M$ is completely decomposable in direct sum of finite dimensional $\mathfrak{l}$-modules; and
    \item The action of $\mathfrak{u}^+$ in $M$ is locally finite.
\end{itemize}
This category has been extensively studied, and we refer the reader to \cite[Chap.~9]{Hum08} for quick overview of this category.

Every simple object of $\cO^\p$ are simple quotients of parabolic Verma modules, i.e, $L_\p(V)$ for some finite dimensional $\mathfrak{p}$-module $V$. In particular, if $\p$ is such that $\g = \h \oplus \p$, by Example~\ref{eg:parabolic-induced} and the fact that $\frkA$ is abelian (Proposition~\ref{prop:properties-frkA}\ref{prop:properties-frkA:abelian}), we have that all simple objects of $\cO^\p$ are contained in $\frkA$. Furthermore, every object of $\cO^\p$ has a finite composition series with simple subquotients in $\cO^\p$. Therefore, since $\frkA$ is closed by extensions (Proposition~\ref{prop:properties-frkA}\ref{prop:properties-frkA:serre}) we obtain:

\begin{cor}\label{cor:parabolic-bgg-O^p}
    For every $\bb \in (\C\setminus\{0\})^n$, the parabolic BGG categories $\cO^{\p_\bb}$ and $\cO^{\tau(\p_\bb)}$ are subcategories of $\frkA$.
\end{cor}

Later, in Corollary~\ref{cor:parabolic-bgg-O^p-if-only-if}, we prove a reciprocal statement of Corollary~\ref{cor:parabolic-bgg-O^p}. Precisely, we prove that if there is an infinite dimensional module $M \in \cO^{\p}\cap \frkA$, for some parabolic subalgebra $\p$ of $\lie{sl}(n+1)$, then $\p$ must be equal to $\p_\bb$ or $\tau(\p_\bb)$, for some $\bb \in (\C\setminus\{0\})^n$.

\subsubsection{Tensor modules}\label{sec:tensor-modules}
First defined in \cite{GN22}, we recall the construction of the \emph{tensor modules}. Then, under some condition in the parameters, we show that we obtain another class of examples of objects of $\frkA$.

We denote by $A_n$ the $n$-th Weyl algebra over $\C$, i.e, the algebra $A_n = \Diff \C[x_1,\dotsc,x_n]$ of differential operators on $\C[x_1,\dotsc,x_n]$ with polynomial coefficients. Denoting the differential operators $\frac{\partial}{\partial x_i}$ by $\partial_i$ we can describe $A_n$ as the $\C$-algebra generated by $x_1,\dotsc,x_n$, $\partial_1,\dotsc,\partial_n$, satisfying the relations
\begin{equation}\label{eq:rel-An}
    \partial_ix_j - x_j\partial_i = \delta_{i,j}\quad \text{(where $\delta_{i,j}$ is the Kronecker delta)},\quad [x_i,x_j] =[\partial_i,\partial_j] = 0,
\end{equation}
for all $1 \le i,j \le n$. 

Define $\setn{n} = \{1,2,\dotsc,n\}$ and let $S \subseteq \setn{n}$. Consider the basis of $\h$ given by $\{\widetilde{h_k} \coloneqq E_{k,k} - \frac{1}{n+1}\sum_{i = 1}^{n+1}E_{i,i}\mid k\in \setn{n}\}$. By \cite[Prop.~2.2]{GN22}, there is a morphism of algebras $\function{\omega_S}{\U(\lie{sl}(n+1))}{A_n\otimes \U(\lie{gl}(n))}$ defined by (for $1 \le i,j,k \le n$):
\begin{align*}
    \widetilde{h_k} 
        &\mapsto
        \begin{cases}
            -x_kd_k \otimes 1 + 1 \otimes E_{k,k} - 1\otimes 1,\quad
                &k \notin S,\\
            x_kd_k \otimes 1 + 1 \otimes E_{k,k},\quad
                &k \in S;
        \end{cases}\\
    E_{i,j} 
        &\mapsto
        \begin{cases}
            1 \otimes E_{i,j} - x_j\partial_i \otimes 1,\quad
                &i,j \notin S,\\
            1 \otimes E_{i,j} + x_i\partial_j \otimes 1,\quad
                &i,j \in S,\\
            1 \otimes E_{i,j} + x_it_j \otimes 1,\quad
                &i \in S\,j \notin S,\\
            1 \otimes E_{i,j} - \partial_i\partial_j \otimes 1,\quad
                &i \notin S,\,j \in S;
        \end{cases}\\
    E_{n+1,j}
        &\mapsto
        \begin{cases}
            -x_j\otimes 1,\quad&i \notin S\\
            -\partial_j\otimes 1,\quad&i \in S;
        \end{cases}\\
    E_{i,n+1}
        &\mapsto
        \begin{cases}
            \sum_{j\notin S}(x_j\partial_j\partial_i \otimes 1 - \partial_j \otimes E_{i,j}) - \sum_{r \in S} (x_rd_r\partial_i \otimes 1 - x_r \otimes E_{i,r})\quad&\\
            \quad-\sum_{j=1}^n \partial_i\otimes E_{j,j} + ((n+1) - |S|)\partial_i \otimes 1,\quad &i \notin S\\
            \sum_{r\in S}(x_ix_rd_r \otimes 1 + x_j \otimes E_{i,r}) - \sum_{j \notin S} (x_ix_j\partial_j \otimes 1 - \partial_j \otimes E_{i,j})\quad&\\
            \quad+\sum_{j=1}^n x_i\otimes E_{j,j} - (n - |S|)x_i \otimes 1,\quad &i \in S.
        \end{cases}
\end{align*}

\begin{defin}[Tensor modules]
    Let $P \in A_n$-module, $V \in \lie{gl}(n)\mdfg$ of finite dimension, and $S \subseteq \setn{n}$. The \emph{tensor module} $T(P, V, S)$ is the $\lie{sl}(n+1)$-module equal to $P \otimes V$ as vector space with action of $\lie{sl}(n+1)$ given by $\omega_S$.
\end{defin}

Let $\mathfrak{t} \subseteq A_n$ be the subspace of $A_n$ with basis $\{x_i\partial_i \mid 1 \le i \le n \}$. Seeing $\mathfrak{t}$ as an abelian lie algebra, one can prove that $\U(\mathfrak{t})$ is isomorphic to the subalgebra of $A_n$ generated by $\mathfrak{t}$ and $1$. From the relations \eqref{eq:rel-An}, one can see that
\[
    [x_i\partial_i,\partial_j] = -\delta_{i,j}\partial_j \quad \text{and}\quad [x_i\partial_i,x_j] = \delta_{i,j}x_j,
\]
for all $1 \le i,j \le n$. So $\mathfrak{t}$ can be seen as an analogous of Cartan subalgebra for $A_n$, and $x_i$, $\partial_i$ are the $\mathfrak{t}$-weight vectors with weights $e_i$ and $-e_i$, respectively, where $\{e_1,\dotsc,e_n\}\subseteq \C^n$ is the canonical basis. This gives us an insight for the construction of the next example:

\begin{eg}[$\mathfrak{t}$-finite tensor modules]\label{eg:tensor-modules}
    Let $P$ be an $A_n$-module, $V$ be a $\lie{gl}(n)$-module of finite dimension, and let $S \subseteq \setn{n}$. We say that the module $T(P, V, S)$ is a $\mathfrak{t}$-finite tensor module if $P$ is finitely generated as $\U(\mathfrak{t})$-module. We claim that a $\mathfrak{t}$-finite tensor module is of $\h$-finite type.

    Let $\{p_1\,\dotsc,p_r\} \subseteq P$ be a finite set of generators of $P$ (as $\U(\mathfrak{t})$-module) and $B \subseteq V$ be a basis consisting of weight vectors. For every $1 \le i,k \le n$ and $v \in B$, notice that
    \[
        (x_kd_k p_i)\otimes v = 
            \begin{cases}
                (-\widetilde{h_k} + \lambda(E_{k,k})-1) (p_i\otimes v),\quad &k \notin S\\
                (\widetilde{h_k} - \lambda(E_{k,k})) (p_i\otimes v),\quad &k \in S,
            \end{cases}
    \]
    where $\lambda \in \h^*$ is the weight of $v$. Then $(\U(\mathfrak{t})p_i)\otimes v \subseteq \U(\h)(p_i\otimes v)$, which implies $P \otimes V = \U(\h)\{p_i \otimes v \mid 1 \le i \le n,\, v \in B\}$. Thus $T(P,V,S) \in \frkA$. 
\end{eg}

\subsection{$\U(\h)$-finite modules for type C}\label{sec:example-typeC}

Next we will give an example constructed by Nilsson in \cite{N16}.

Seeing $\lie{sp}(2n)$ as a subalgebra of $\lie{gl}(2n)$, we fix a Cartan subalgebra $\h \subseteq \lie{sp}(2n)$ with basis 
\[
    \{\widetilde{h_i} = E_{i,i}- E_{n+i,n+i}\mid 1 \le i \le n\}.
\]
Let $\{\epsilon_i\}$ be the basis of $\h^*$ dual to $\{h_i\}$. Then the root system of $\lie{sp}(2n)$ related to $\h$ is 
\[
    \Delta = \{\pm \epsilon_i \pm \epsilon_j \mid 1 \le i,j \le n\} \setminus \{0\}.
\]
We fix the root vector of $\lie{sp}(2n)$ as follows
\[
    e_{\epsilon_i + \epsilon_j} = E_{i,n+j} + E_{j,n+i},\quad e_{-\epsilon_i - \epsilon_j} = -E_{n+i,j} - E_{n+j,i}
    \quad \text{and}\quad e_{\epsilon_i - \epsilon_k} = E_{i,k} - E_{n+k,n+i},
\]
for $1 \le i,j,k \le n$ and $i \ne k$.

Consider the algebra automorphisms $\{\function{\sigma_i}{\U(\h)}{\U(\h)}\mid 1 \le i \le n\}$ where $\sigma_i(\widetilde{h_j}) =\widetilde{h_j} - \delta_{i,j}$. We define a $\lie{sp}(2n)$-module $M_0$ as the regular $\U(\h)$-module $\C[\widetilde{h_1},\dotsc,\widetilde{h_n}] \simeq \U(\h)$, with action of the weight vectors given by
\begin{align*}
    e_{2\epsilon_i} \cdot f &= \left(\widetilde{h_i} - \frac{1}{2}\right)\left(\widetilde{h_i} - \frac{3}{2}\right)\sigma_i^2(f),\\
    e_{-2\epsilon_i} \cdot f &= \sigma_i^{-2}(f),\\
    e_{\epsilon_i + \epsilon_j}\cdot f &= \left(\widetilde{h_i} - \frac{1}{2}\right)\left(\widetilde{h_j} - \frac{1}{2}\right)\sigma_i\sigma_j(f),\\
    e_{-\epsilon_i - \epsilon_j} \cdot f &= \sigma_i^{-1}\sigma_j^{-1}(f),\\
    e_{\epsilon_i - \epsilon_j} \cdot f &= \left(\widetilde{h_i} - \frac{1}{2}\right)\sigma_i\sigma_j^{-1}(f),
\end{align*}
for all $f \in \C[\widetilde{h_1},\dotsc,\widetilde{h_n}]$, $1 \le i, j \le n$, with $i \ne j$.
Clearly $M_0 \in \frkA$. In fact, $M_0$ is $\h$-free of rank $1$.

\begin{prop}[{\citealp[Prop.~13]{N16}}]\label{prop:M_0-simple}
    The module $M_0$ is simple.
\end{prop}
\begin{proof}
    Notice that $(1-e_{-2\epsilon_i})$ decrease the degree of the variable $\widetilde{h_i}$ by one. So, applying actions from $\Vspan\{1-e_{-2\epsilon_i}\mid i \in \setn{n}\}$, one can reduce any element $f \in M_0$ to one. Since $M_0 = \U(\h)$, it follows the proposition.
\end{proof}

%\newpage

%%%%%%%%%%%%%%%%%%%%%%%%%%%%%%%%%%%%%%%%%%%%%%%%%%%%%%%%%%%%%%%%%%%
% Section - Weighting Functor
%%%%%%%%%%%%%%%%%%%%%%%%%%%%%%%%%%%%%%%%%%%%%%%%%%%%%%%%%%%%%%%%%%%
%%%%%%%%%%%%%%%%%%%%%%%%%%%%%%%%%%%%%%%%%%%%%%%%%%%%%%%%%%%%%%%%%%%
%%%%%%%%%%%%%%%%%%%%%%%%%%%%%%%%%%%%%%%%%%%%%%%%%%%%%%%%%%%%%%%%%%%
% Chapter 5 - Weighting functor
%%%%%%%%%%%%%%%%%%%%%%%%%%%%%%%%%%%%%%%%%%%%%%%%%%%%%%%%%%%%%%%%%%%
%%%%%%%%%%%%%%%%%%%%%%%%%%%%%%%%%%%%%%%%%%%%%%%%%%%%%%%%%%%%%%%%%%%

\section{Weighting functors}\label{sec:weighting-functors}

In this chapter, we investigate another connection between $\U(\h)$-finite $\g$-modules and admissible weight modules. Such connection is given by the weighting functor $\W$ that, as the name says, assigns any $\g$-module to a module with weight decomposition. In fact, when applied to objects of $\frkA$, the weighting functor produced an admissible weight module.

We start by recalling the definition of the weighting functor and construct their left-derived functors. We continue by investigating their properties. After, we apply $\W$ to the category $\frkA$ and obtain some statements about their freeness over $\U(\h)$. Specifically, we prove that a simple $\U(\h)$-finite $\g$-module $M$, of infinite dimension, is locally $\U(\h)$-free except in a finite set of points of $\h^*$. Moreover, we also prove that $M$ must be $\U(\h)$-torsion free. As a consequence, we conclude that there exist infinite dimensional simple $\U(\h)$-finite modules only for $\g$ of type A and C.

\subsection{Definitions of $\W_*$ and its properties}

The \emph{weighting functor}, introduced by J. Nilsson in \cite{N16} following a suggestion by O. Mathieu, is defined as follows:
\[
    \function{\W}{\g\mdfg}{\h\mdfg}, \quad M \mapsto \W(M) = \bigoplus_{\lambda \in \h^*} M/\mathfrak{m}_\lambda M.
\]
In fact, we can define a $\g$-action on $\W(M)$ such that $\W$ becomes an endofunctor of $\g$-modules, in the following way:
any weight vector $x_\alpha \in \g_\alpha$ acts on $\W(M)$ by
\[
    x_\alpha \cdot (m + \mathfrak{m}_\lambda M) \coloneqq x_\alpha m + \mathfrak{m}_{\lambda + \alpha}M,
\]
for all $m \in M$.
On the level of morphism, given $\function{f}{M}{N}$ in $\U(\g)\md$, we define $\function{\W(f)}{\W M}{\W N}$ by $\W(f)(m + \mathfrak{m}_\lambda M) = f(m) + \mathfrak{m}_\lambda M$, for all $m \in M$.

A $\g$-module $M$ has \emph{weight decomposition} if it has decomposition in weight spaces, which not necessarily have finite dimension. Let $\widetilde{\frkW}$ be the subcategory of $\U(\g)\mdfg$ consisting of $\g$-modules with a weight decomposition.

\begin{lem}[\cite{N16}, Lemma~9]\label{lem:W}
    The weighting functor $\W$ satisfies:
    \begin{enumerate}[label=(\roman*),ref=(\roman*)]
        \item $\W$ assigns to any $\g$-module $M$ a module with a weight decomposition, i.e., $\W$ is an additive functor from $\g\mdfg$ to $\widetilde{\frkW}$.
        \item The restriction of $\W$ to the category $\widetilde\frkW$ is the identity, and hence, $\W$ is an idempotent functor.
        \item $\W$ preserves the (generalized) central character: if $M$ admits (generalized) central character $\chi_M$, then so does $\W(M)$ and $\chi_\W(M) = \chi_M$.
    \end{enumerate}
\end{lem}

Notice that $\W\circ \Res_\h^\g$ is naturally isomorphic to the functor
    \[
        \bigoplus_{\lambda \in \h^*} - \otimes_{\U(\h)} \U(\h)/\mathfrak{m}_{\lambda},
    \]
which is right exact. 
Thus:

\begin{lem}
    The weight functor $\W$ is right exact.
\end{lem}

Let $\function{\W_*}{\g\mdfg}{\g\mdfg}$ be the left-derived functors of $\W$. Since any $\g$-module has a projective resolution, the functor $\W_k$ is well-defined for all $k \in \N$. Furthermore, since $\W$ is an additive functor, so is $\W_k$, for all $k \in \N$. Recall that $\W_0 = \W$.

The following lemma gives two description of $\W_*$, one in terms of $\Tor$ functors and the other in terms of Lie algebra homologies.
While the first description is more natural when applying commutative algebra theory, the second one is better for Lie-theoretical perspectives. 

For a weight $\lambda \in \h^*$, let $\C_\lambda$ be the one dimensional $\h$-module of weight $\lambda$. Then $\C_\lambda$ is the $\U(\h)$-module isomorphic to $\U(\h)/\m_\lambda$.

\begin{lem}\label{lem:LW-is-Tor}
    The functor $\Res_\h^\g \circ \W_*$ is naturally isomorphic to the functors 
    \[
        \bigoplus_{\lambda \in \h^*} \Tor_*^{\U(\h)}(\Res_{\h}^{\g}(-), \C_\lambda)\simeq \bigoplus_{\lambda \in \h^*} H_*(\h, -\otimes \C_{-\lambda})\otimes \C_\lambda.
    \]  
    Furthermore, for any $M \in \g\mdfg$ and $k\in \N$, $\W_k(M)$ is a module with weight decomposition with respective weight spaces
    \[
        \W_k(M)_\lambda \simeq \Tor_*^{\U(\h)}(M,\C_\lambda) \simeq H_*(\h, M\otimes \C_{-\lambda})\otimes \C_\lambda,
    \]
    for all $\lambda \in \h^*$. And in particular, $\W_*$ is a functor from $\g\mdfg$ to $\widetilde{\frkW}$. 
\end{lem}
\begin{proof}
    By PBW theorem, we have that $\U(\g)$ is free as $\U(\h)$-module. Thus, for any $M \in \g\mdfg$, a projective resolution of $M$ in $\g\mdfg$ induces a projective resolution of $M$ in $R\md$.
    The isomorphism of functor in terms of $\Tor$ follows when we observe that $\function{\Res_\h^\g \circ \W}{\g \md}{\U(\h)\md}$ is naturally isomorphic to $\bigoplus_{\lambda \in \h^*} - \otimes_{\U(\h)} \U(\h)/\frak{m}_\lambda$. To obtain the characterization of $\Res_\h^\g \circ \W$ via homology of Lie algebras, we first consider the isomorphism of functors of $\h$-modules $\otimes_{\U(\h)} \C_\lambda \simeq ((-\otimes\C_{-\lambda}) \otimes_{\U(\h)} \C_0) \otimes \C_{\lambda}$.
    Moreover, one can check that $-\otimes \C_{\lambda}$ is an exact endofunctor of $\h$-modules. Thus, we have natural isomorphisms
    \[
        \Tor_*^{\U(\h)}(-,\C_\lambda) \simeq \Tor_*^{\U(\h)}(-\otimes\C_{-\lambda},\C_0) \otimes \C_\lambda \simeq H_*(\h, -\otimes\C_{-\lambda}) \otimes \C_\lambda,
    \]
    where the second isomorphism follows from Chevalley-Eilenberg resolution, which implies the characterization of homology if Lie algebras via Tor functor. Then the first claim is proved.
        
    The second claim follows from the fact that, for any $\lambda \in \h^*$ and $\h$-module $P$, the module $P\otimes_{\U(\h)} \U(\h)/\mathfrak{m}_\lambda$ is annihilated by the action of $h-\lambda(h)$, $h \in \h$, since $h-\lambda(h)\in \mathfrak{m}_\lambda$. So $h-\lambda(h)$ also annihilates any subquotient of $P\otimes_{\U(\h)} \U(\h)/\mathfrak{m}_\lambda$. Then we have that $h-\lambda(h)$ acts by zero in the homology of $(P_* \otimes_{\U(\h)} \C_\lambda,d_*)$ for any free $\U(\g)$ resolution $(P_*,d_*)$ of $M$. Thus, $\W_*(M)_\lambda \simeq \Tor_*^{\U(\h)}(M,\C_\lambda)$ and the Lemma is proved.
\end{proof}

\begin{lem}\label{lem:tensor-free}
    Let $M \in \g\mdfg$, then $F \otimes M$ is a free $\U(\g)$-module whenever $F$ is $\U(\g)$-free.
\end{lem}
\begin{proof}
    It is enough to show that $\U(\g)\otimes M$ is free. Let $M^0$ denote the vector space $M$ with trivial $\g$-action, we claim that $\U(\g)\otimes M \simeq \U(\g)\otimes M^0$ as $\g$-module.
    
    Let $\function{\Delta}{\U(\g)}{\U(\g) \otimes \U(\g)}, \,\function{S}{\U(\g)}{\U(\g)}$ be respectively the coproduct and the antipode of $\U(\g)$. In what follows, we will use the Sweedler notation for $\Delta$, i.e, we write $\Delta(u) = u_{(1)}\otimes u_{(2)}$, for $u \in \U(\g)$. In these notations, it is satisfied that
    \begin{equation}\label{eq:Delta-xu}
        \Delta(xu) = xu_{(1)} \otimes u_{(2)} + u_{(1)}\otimes xu_{(2)},
    \end{equation}
    for all $x \in \g$ and $u \in \U(\g)$.
    Now, define the following maps
    \[
        \function{\varphi}{\U(\g)\otimes M}{\U(\g)\otimes M^0},\quad u\otimes m \mapsto u_{(1)}\otimes S(u_{(2)})m,
    \]
    \[
        \function{\psi}{\U(\g)\otimes M^0}{\U(\g)\otimes M},\quad u\otimes m \mapsto u_{(1)}\otimes u_{(2)}m,
    \]
    for all $u \in \U(\g)$, $m \in M$.
    A straightforward computation, using equation \eqref{eq:Delta-xu} and properties of the antipode of $\U(\g)$, shows that $\varphi$ and $\psi$ are morphisms of $\g$-modules.
    To see that $\psi$ and $\varphi$ are mutual inverse, we compute:
    \[
        \psi(\varphi(u\otimes m)) = \psi(u_{(1)}\otimes S(u_{(2)})m) = u_{(1)_{(1)}} \otimes u_{(1)_{(2)}}S(u_{(2)})m = u_{(1)} \otimes u_{(2)_{(1)}}S(u_{(2)_{(2)}})m,
    \]
    for all $u \in \U(\g)$ and $m \in M$, where the last equality follows from the fact that $(\Delta \otimes \Id)\Delta = (\Id \otimes \Delta)\Delta$. Furthermore, we have that $u_{(1)}S(u_{(2)}) = \epsilon(u)$, for all $u \in U$, where $\function{\epsilon}{\U(\g)}{\C}$ is the counite of $\U(\g)$. Then $\psi(\varphi(u\otimes m)) = u_{(1)} \otimes \epsilon(u_{(2)})m = u\otimes m$. And hence $\psi \circ \varphi = \Id$. Analogously, using now that $S(u_{(1)})u_{(2)} = \epsilon(u)$ for all $u \in \U(\g)$, we show that $\varphi\circ \psi  = \Id$. And the lemma is proved.
\end{proof}

As expected, such functors have similar properties as $\W$:

\begin{lem}\label{lem:prop-Wk}
    The left-derived $\W_*$ of the weight functor $\W$ satisfies:
    \begin{enumerate}[label=(\roman*),ref=(\roman*)]
        \item\label{lem:prop-Wk-commut-tensor-weight}{For any module $W \in \overline{\frkW}$ the functor $\W_*$ commutes with the functor $- \otimes W$, i.e, 
            \[
                \W_*( - \otimes W) \simeq \W_*(-)\otimes W
            \]
            are naturally isomorphic. In particular, $\W_k(W) = W^{\oplus \binom{n}{k}}$, for $k > 0$;}
        \item\label{lem:prop-Wk-annihilatator}{Let $M\in \g\mdfg$, then $\Ann_{\U(\g)} M \subseteq \Ann_{\U(\g)} \W_*(M)$;}
        \item\label{lem:prop-Wk-central-char}{$\W_*$ preserves (generalized) central character: if $M$ admits (generalized) central character $\chi_M$, then so does $\W_*(M)$ and $\chi_{\W_*(M)} = \chi_M$;}
        \item\label{lem:prop-Wk-projective-functor}{$\W_*$ commutes with translation functors.}
    \end{enumerate}
\end{lem}
\begin{proof}    
    \ref{lem:prop-Wk-commut-tensor-weight}
    We first prove the claim for $\W_0 = \W$. Let $W \in \widetilde{\frkW}$ and $M \in \g\mdfg$. We will show that $\W(M \otimes W) \simeq \W(M)\otimes W$.
    For any $\mu \in \h^*$, define $\function{\widetilde\varphi_\mu}{M\otimes W}{\W(M)\otimes W}$ by $\widetilde\varphi_\mu(m \otimes w_\lambda) = (m + \mathfrak{m}_{\mu-\lambda} M)\otimes w_\lambda$, where $w_\lambda \in W$ is a weight vector of weight $\lambda \in \h^*$. A straightforward computation shows that $\widetilde\varphi_\mu$ is a well-defined morphism of $\U(\h)$-module with kernel containing $\mathfrak{m}_{\mu}(M\otimes W)$. Thus, we have the induced morphism of $\U(\h)$-modules
    \[
        \function{\varphi_\mu}{\dfrac{M\otimes W}{\mathfrak{m}_{\mu}(M \otimes W)}}{\W(M)\otimes W}.
    \]
    Define the morphism of $\U(\h)$-module
    \[
        \function{\varphi}{\W(M \otimes W)}{\W(M)\otimes W},\quad  m\otimes w + \mathfrak{m}_\mu (M\otimes W) \mapsto \varphi_\mu(m\otimes w).
    \]
    We show that $\varphi$ is an isomorphism by constructing an inverse. Define
    \[
        \function{\psi}{\W(M)\otimes W}{\W(M\otimes W)}, \quad (m + \mathfrak{m}_\mu M)\otimes w_\lambda \mapsto m\otimes w_\lambda + \mathfrak{m}_{\mu + \lambda}(M\otimes W),
    \]
    for all $\mu \in \h^*$ and vector $w_\lambda \in W_\lambda$ of weight $\lambda\in \h^*$. The morphism $\psi$ is well defined and one can check  that $\varphi\psi = \Id_{\W(M)\otimes W}$ and  $\psi\varphi = \Id_{\W(M\otimes W)}$. Finally, for all weight vector $x \in \g_{\nu}$, we have:
    \begin{align*}
        \psi(x \cdot ((m + \mathfrak{m}_\mu M)\otimes w_\lambda))
            &= \psi( (xm + \mathfrak{m}_{\mu + \nu}M)\otimes w_\lambda + (m + \mathfrak{m}_\mu M)\otimes xw_\lambda)\\
            &= xm\otimes w_\lambda + \mathfrak{m}_{\mu + \nu + \lambda}(M \otimes W) + m\otimes xw_\lambda + \mathfrak{m}_{\mu +\lambda + \nu} (M \otimes W)\\
            &=x \cdot \psi((m + \mathfrak{m}_\mu M)\otimes w_\lambda),
    \end{align*}
    Hence $\psi$ is a morphism of $\g$-modules and then $\W(M\otimes W)$ and $\W(M)\otimes W$ are isomorphic as $\g$-modules.
    It is straightforward to check that $\psi$ and $\varphi$ defines natural transformations between the functors $\W(- \otimes W)$ and $\W(-)\otimes W$. Thus, the functors are naturally isomorphic.
        
    Now we prove the claim for $\W_k$ with $k > 0$. Let $P_*\xrightarrow{d_*} M\rightarrow 0$ be a free resolution of $M$ in $\g\mdfg$. I claim that $P_*\otimes W \xrightarrow{d_*\otimes \Id_W} M\otimes W\rightarrow 0$ is a free resolution of $M\otimes W$ in $\g\mdfg$.
    Indeed, from Lemma~\ref{lem:tensor-free}, we have that $P_*\otimes W$ is a free $\g$-module. Furthermore, the fact that $- \otimes W$ is an exact functor in the category of vector spaces implies that the resolution is free. Applying the functor $\W$ and using the fact that $\W(- \otimes W)$ and $\W(-)\otimes W$ are naturally isomorphic, we obtain the commutation of the following diagram,
    \[
        \xymatrix{
            {\W(P_*\otimes W)} \ar[rr]^-{\W(d_*\otimes \Id_W)} \ar[d]     && M\otimes W \ar[d] \ar[r]  & 0\\
            {\W(P_*)\otimes W} \ar[rr]_-{\W(d_*)\otimes \Id_W}            && M\otimes W \ar[r]         & 0
        }
    \]
    where the vertical maps are isomorphism. Therefore, the homology of the first line is naturally isomorphic to the homology of the second line. But the exactness of $-\otimes W$ implies that the homology of $(\W(P_*)\otimes W,\W(d_*)\otimes \Id_W)$ is $\W_*(M)\otimes W$. Therefore, $\W_*( - \otimes W)$ and $\W_*(-)\otimes W$ are naturally isomorphic.
    
    Let $\C$ be the trivial $\g$-module. Then, as $\U(\h)$-module, it is isomorphic to $\U(\h)/\mathfrak{m}_0$. Then $\Res_\h^\g \W_k(\C) \simeq \bigoplus_{\lambda \in \h^*} \Tor_k^{\U(\h)}(\U(\h)/\mathfrak{m}_0, \U(\h)/\mathfrak{m}_\lambda)$. Since the localization $(\U(\h)/\mathfrak{m}_\lambda)_{\mathfrak{m}_\mu}$ vanish if $\mu \ne \lambda$, we have that $\Tor_k^{\U(\h)}(\U(\h)/\mathfrak{m}_0, \U(\h)/\mathfrak{m}_\lambda) = 0$ if $\lambda \ne 0$. Therefore, $\W_k(\C) \simeq \Tor_k^{\U(\h)}(\C,\C) \simeq H_k(\h, \C)$ with the trivial action of $\U(\g)$. An easy computation with Chevalley-Eilemberg complex shows that $H_k(\h, \C) = \bigwedge^k\h$.
    Hence $\W_k(\C) = \C^{\binom{n}{k}}$, for $k > 0$. And finally, for any weight module $W$, we have that $\W_k(W) \simeq \W_k(\C \otimes W) \simeq \W_k(\C)\otimes W= W^{\oplus{\binom{n}{k}}}$, for all $k > 0$.
    
    \ref{lem:prop-Wk-annihilatator}
    Notice that $I = \Ann_{\U(\g)}(M)$ is a two-sided ideal of $\U(\g)$ and hence a submodule of the adjoint representation of $\g$ in $\U(\g)$. In particular, $I$ is a weight module. Furthermore, for any $N\in \g\mdfg$, the multiplication map $\function{\mu_N}{I\otimes N}{N}$, $\mu_N(u\otimes n) = un$, is a morphism of $\g$-modules.
    In addition, notice that the natural isomorphism between $\W(I\otimes N)$ and $I\otimes \W(N)$ interchanges $\W(\mu_N)$ and $\mu_{\W(N)}$.
    
    Now let $P_*\xrightarrow{d_*} M\rightarrow 0$ be a free resolution of $M$ in $\U(\g)\md$.
    Since $d_*$ commutes with left multiplication by elements of $\U(\g)$, we have $u d_* = d_* u$, for all $u \in \U(\g)$. Then the following is a commutative diagram in $\U(\g)\md$:
    \begin{equation}\label{diagram}
    \xymatrix{
        {I\otimes P_*} \ar[r]^{\Id_I\otimes d_*} \ar[d]_{\mu_{P^*}}     & I\otimes M \ar[d]^{\mu_M} \ar[r]  & 0\\
        P_* \ar[r]_{d_*}                                                & M \ar[r]                          & 0
    }
    \end{equation}
    We will construct inductively a family of morphism of $\U(\g)$-modules $\function{s_i}{I\otimes P_i}{P_{i+1}}$, for $i \ge 0$, such that
    \begin{equation}\label{eq:homotopia}
        d_{i+1} s_i - s_{i-1} \Id_I \otimes d_i = \mu_{P_i}, \quad i \ge 1.
    \end{equation}
    Notice that 
    \[
        d_0(\mu_{P_0}(u\otimes p)) = \mu_M(\Id_I\otimes d_0(u \otimes p)) = \mu_M(u\otimes d_0(p)) = ud_0(p) = 0,
    \]
    for all $u \in I$, $p \in P_0$. And since $I\otimes P_0$ is free, there exists a morphism of $\g$-modules $\function{s_0}{I\otimes P_0}{P_{1}}$ such that $d_1 s_0 = \mu_{P_0}$. Now assume that $s_j$ is defined, for all $0 \le j \le i-1$. Then, for every $u \in I$, $p \in P_i$, we have
    \[
        d_i\mu_{P_i}(u\otimes p) = ud_i(p) = \mu_{P_{i-1}}(u\otimes d_i(p)) = (d_i s_{i-1} - s_{i-2} \Id_I \otimes d_{i-1})(u \otimes d_i(p)) = d_i s_{i-1}(u \otimes d_i(p)),
    \]
    and then $\mu_{P_i}(u\otimes p) - s_{i-1}(u \otimes d_i(p)) \in \Ker d_i = \im d_{i+1}$. Therefore, there exists a morphism of $\g$-modules $\function{s_{i}}{I\otimes P_{i}}{P_{i+1}}$ such that
    \[
        d_{i+1}s_{i} = \mu_{P_i} - s_{i-1} \Id_I \otimes d_i,
    \]
    i.e., $s_{i}$ satisfies equation \eqref{eq:homotopia}.
    
    Applying the functor $\W$ in equations \eqref{diagram} and \eqref{eq:homotopia} and using the natural isomorphism $\W(I\otimes -) \simeq I\otimes \W(-)$ we obtain the commutative diagram
    \begin{equation}
    \xymatrix{
        {I\otimes \W(P_*)} \ar[d]_{\mu_{\W(P_*)}} \ar[r]   &I\otimes \W(M) \ar[d]^{\mu_{\W(M)}} \ar[r]     & 0 \\
	{\W(P_*)} \ar[r]                                    &\W(M) \ar[r]                                   & 0
    },
    \end{equation}
    and the class of morphism $\function{\widetilde{s_i}}{I\otimes \W(P_i)}{\W(P_{i+1})}$ satisfying
    \[
        \W(d_{i+1}) \widetilde{s_i} - \widetilde{s}_{i-1} \Id_I \otimes \W(d_i) = \mu_{\W(P_i)}, \quad i \le 1.
    \]
    In particular $\widetilde{s_i}$ defines a homotopy between $\mu_{\W(P_i)}$ and the zero morphism between complexes. Hence, the morphism $\function{\mu_{\W_*(M)}}{I\otimes \W_*(M)}{\W_*(M)}$ induced in the homology is the zero morphism. Therefore $I$ annihilates $\W_*(M)$.

    \ref{lem:prop-Wk-central-char}
    Suppose that $M$ has a generalized central character $\chi_M$. Since $M$ is finite generated as $\g$-module, we have that, given $z \in Z(\g)$, there is $N \in \N$ such that $(z - \chi_M(z))^N$ acts by zero in $M$. Hence so it does in $\W_*(M)$, by item~\ref{lem:prop-Wk-annihilatator}.

    \ref{lem:prop-Wk-projective-functor}
    The last item is a direct consequence of the additiveness of $\W_k$ and items \ref{lem:prop-Wk-commut-tensor-weight} and \ref{lem:prop-Wk-central-char}.
\end{proof}

\subsection{Weighting functor in $\frkA$}

Let $M\in \frkA$ be a simple $\U(\h)$-finite module of infinite dimension, and let $\chi_M$ be its central character.

\begin{prop}\label{prop:W(M)-is-admissible}
   The module $\W_*(M)[t]$ is admissible for all $t \in T^*$. 
\end{prop}
%proof with commutative algebra point of view
\rascunho{\begin{proof}
    For all $\lambda\in\h^*$ and $k \in \N$, Lemma~\ref{lem:LW-is-Tor} says that $\W_k(M)_\lambda = \Tor^{\U(\h)}_k(M,\U(\h)/\mathfrak{m}_\lambda)$. Then, for any projective resolution $(P_*,d_*)$ of $M$, we have
    \[
        \dim \W_k(M)_\lambda = \dim \Tor^{\U(\h)}_k(M,\U(\h)/\mathfrak{m}_\lambda) \le \dim P_k \otimes_{\U(\h)} \U(\h)/\mathfrak{m}_\lambda.
    \]
    Since $M$ is finitely generated as $\U(\h)$-module and $\U(\h)$ is Noetherian ring, we can assume that $P_k$ is free of finite rank, let say $r$. Then $\dim \W_k(M)_\lambda \le \dim \left(\U(\h)/\mathfrak{m}_\lambda\right)^{r}$. Since $\U(\h)/\mathfrak{m}_\lambda \simeq \C$ as vector space, we have that, for all $\lambda$ $\dim \W_k(M)_\lambda \le r$, and the proposition follows together with Remark~\ref{rem:adm-implies-g-fg}.
\end{proof}}%%% fim...
\begin{proof}
    Let $\lambda \in \h^*$ and $k \in \N$. Then Lemma~\ref{lem:LW-is-Tor} gives us an isomorphism of vector spaces $\W_k(M)_\lambda \simeq H_k(\h,M\otimes \C_{-\lambda})$. So we can use projective $\U(\h)$-resolution to compute the desired dimension. Since $M$ is $\U(\h)$-finitely generated, there exists a projective resolution $(P_*,d_*)$ of $M$ where $P_k$ is $\U(\h)$-free of finite rank. Clearly $(P_*\otimes \C_{-\lambda},d_*\otimes \C_{-\lambda})$ forms a projective $\U(\h)$ resolution of $M\otimes \C_{-\lambda}$. Hence
    \[
        \dim \W_k(M)_\lambda = \dim H_k(\h,M\otimes \C_{-\lambda}) \le \rank P_k \otimes \C_{-\lambda} = \rank P_k.
    \]
    And the proposition follows together with Remark~\ref{rem:adm-implies-g-fg}.
\end{proof}

The next step is to investigate how close to $\U(\h)$-free can $M$ be. The idea is to start by study $\W_1$, which is isomorphic to $\bigoplus_{\mathfrak{m}} \Tor_1^{\U(\h)}(-, \U(\h)/\mathfrak{m})$. Then, since flatness and freeness are equivalent properties in the category of modules over polynomial rings, the support of $\W_1(M)$ express the points where $M$ fails to be locally free.

\begin{lem}\label{lem:calM-not-dense-supp}
    Let $k \ge 1$. Then there exists a proper closed set $Z \subseteq \Specm \U(\h)$ such that $\Supp \W_k(M) \subseteq Z$.
\end{lem}
\begin{proof}
    First, I claim that there exists a non-empty open subset $U\subseteq \Specm \U(\h)$ such that $M_{\mathfrak{m}}$ is $\U(\h)_{\mathfrak{m}}$-free for all $\mathfrak{m} \in U$. Indeed, since $\U(\h)$ is a commutative domain, it follows from Generic freeness (see \cite[Lemma~2.6.3]{Dix96}) that there exists some $f \in \U(\h)\setminus \{0\}$ such that $M_f$ is $\U(\h)_f$-free. Then the set $U = \{\mathfrak{m} \in \Specm{\U(\h)}\mid f \notin \mathfrak{m}\}$ is non-empty and, for all $\mathfrak{m} \in U$, we have
    \[
        M_{\mathfrak{m}} \simeq (M_{\mathfrak{m}})_{\bar{f}} \simeq (M_f)_{\overline{\mathfrak{m}}},
    \]
    where $\overline{f}$ is the image of $f$ in $\U(\h)_{\mathfrak{m}}$ and $\overline{\mathfrak{m}}$ is the image of $\mathfrak{m} \in \U(\h)_f$. Since $M_f$ is free, so is $M_{\mathfrak{m}}$.

    Now, for all $\mathfrak{q} \in \Specm \U(\h)$ and $\mathfrak{m} \in U$, we have
    \[
        \left(\Tor_k^{\U(\h)}(M,\U(\h)/\mathfrak{m})\right)_{\mathfrak{q}} \simeq \Tor_k^{\U(\h)_{\mathfrak{q}}}(M_{\mathfrak{q}},(\U(\h)/\mathfrak{m})_{\mathfrak{q}}) = 0,
    \]
    since $(\U(\h)/\mathfrak{m})_{\mathfrak{q}} = 0$ if $\mathfrak{q} \ne \mathfrak{m}$, and $M_{\mathfrak{q}}$ is free otherwise.
    Hence $\Supp \W_k(M) \subseteq \Specm \U(\h)\setminus U$.
\end{proof}

\begin{prop}\label{prop:finite-dim-Wk}
    For $k \ge 1$ the module $\W_k(M)$ is finite dimensional.
    Furthermore, $\W_{k}(M)$ vanishes whenever $\chi_M$ is singular or non-integral central character.
\end{prop}
\begin{proof}
    Let $t \in T^*$. Since $\W_k(M)[t]$ is an admissible weight module, by Lemma~\ref{lem:first-prop-admW}\ref{lem:first-prop-admW:supp-ess-dense} it is either finite-dimensional or has Zariski-dense support. But by Lemma~\ref{lem:calM-not-dense-supp}, the second option can not happen. Thus $\W_k(M)[t]$ is finite-dimensional.
    
    In particular, $\W_k(M)[t]$ must have integral generalized central character. Then, together with the fact that $\W_k$ preserves generalized central character (Lemma~\ref{lem:prop-Wk}\ref{lem:prop-Wk-commut-tensor-weight}), it follows that $\W_k(M)[t] = 0$ if $\chi_M$ is singular or non-integral, for all $t \in T^*$. Otherwise, $\W_k(M)[t] \ne 0$ only for the case $t = \lambda + Q$, where $\lambda \in \h^*$ is a dominant integral weight such that $\chi_M = \chi_\lambda$. Hence $\W_k(M) = \W_k(M)[\lambda+ Q]$ is finite-dimensional
\end{proof}

Finally, we can prove that an infinite dimensional $\U(\h)$-finite module is locally $\U(\h)$-free in $\Specm \U(\h)$, outside a finite set of points.

\begin{theo}\label{theo:almost-free}
    The module $M$ is locally free over $\U(\h)$ in $\h^*\setminus \Supp \W_1(M)$. 
    Furthermore, if $\chi_M$ is either non-integral or integral and singular central character, then $M$ is $\U(\h)$-free.
\end{theo}
\begin{proof}
    Let $\mathcal{M} = \W_{1}(M)$. By Proposition \ref{prop:finite-dim-Wk}, we have that $\mathcal{M}$ is finite-dimensional. Therefore $\Supp(\mathcal{M})$ is finite. And for all $\mathfrak{m} \notin \Supp(\mathcal{M})$ we have that
    \[
        0 = \mathcal{M}_{\mathfrak{m}} = \Tor^{\U(\h)_{\mathfrak{m}}}_1(M_{\mathfrak{m}}, \U(\h)/{\mathfrak{m}}).
    \]
    From a minimal free resolution $(F_*, d_*)$ of $M_{\mathfrak{m}}$ we obtain that $\Tor^{\U(\h)_{\mathfrak{m}}}_1(M_{\mathfrak{m}}, \U(\h)/{\mathfrak{m}}) = \ker{d_1\otimes\U(\h)/{\mathfrak{m}}} = F_1/\mathfrak{m}F_1$. Then $M_{\mathfrak{m}}$ is free, since such Tor vanishes.
    Proving the first claim.
 
    In particular, if $M$ has non-integral or singular generalized central character, by Prop.~\ref{prop:finite-dim-Wk}, $\M$ must vanish. Then $M$ is locally free. Since $\U(\h)$ is Noetherian, $M$ is also finitely presented. Then $M$ is a projective $\U(\h)$-module, by commutative algebra theory. Thus, since $\U(\h)$ is isomorphic to a polynomial ring with variable in a basis of $\h$, Quillen-Suslin theorem (see \cite{Qui76}) implies that $M$ is $\U(\h)$-free.
\end{proof}

Next, we turn our attention to investigate the torsion of $M$ as $\U(\h)$-module. 

First, we need to prove some results about modules over commutative algebra. In the following, let $R$ be a polynomial $\C$-algebra in finite number of variables, and let $\Spec R$ denote the set of prime ideals of $R$. Consider a finitely generated $R$-module $N$ and its support $\Suppca N$, in the sense of commutative algebra theory, i.e, the set
\[
        \Suppca N = \{\mathfrak{p} \in \Spec R \mid N_\mathfrak{p} \ne 0\}.
\]

\begin{lem}\label{lem:finite-support-finite-dim}
    If $\Suppca N \cap \Specm R$ is finite, then $N$ is finite dimensional over $\C$.
\end{lem}
\begin{proof}
    Let $\mathfrak{p}\in \Spec R$ and let $V(\mathfrak{p})\subseteq \Spec R$ be the set of prime ideals containing $\mathfrak{p}$. We start by claiming that if $V(\mathfrak{p})\cap \Specm R$ is finite then $\dim_\C R/\mathfrak{p}$ is finite.

    Indeed, let $V(\mathfrak{p})\cap \Specm R = \{\mathfrak{m}_1,\dotsc, \mathfrak{m}_r\}$. By Nullstellensatz, we have $\mathfrak{p} = \sqrt{\mathfrak{p}} = \mathfrak{m}_1\cap \dotsb \cap \mathfrak{m}_r$. Then, by Chinese remainder theorem, it follows that $\frac{R}{\mathfrak{p}} = \frac{R}{\prod_{i}\mathfrak{m}_i} \simeq \bigoplus_i \frac{R}{\mathfrak{m}_i} \simeq \C^r$. Thus the claim.

    Since $N$ is finitely generated and $R$ is noetherian, there exists a finite filtration of $R$-submodules $0 = N_0 \subseteq N_1 \subseteq \dotsb \subseteq N_r = N$ and prime ideals $\mathfrak{p}_1,\dotsc, \mathfrak{p}_r \in \Spec R$, such that $N_i/N_{i-1} \simeq R/\mathfrak{p}_i$, for all $1 \le i \le r$. Then, for every $1 \le i \le r$ and $\mathfrak{m}\in V(\mathfrak{p}_i)\cap \Specm R$ we have
    \[
        \dfrac{N_\mathfrak{m}}{(N_{i-1})_\mathfrak{m}} \supseteq \left( \dfrac{N_i}{N_{i-1}}\right)_\mathfrak{m} \simeq \left( \dfrac{R}{\mathfrak{p_i}}\right)_\mathfrak{m} \simeq \dfrac{R_\mathfrak{m}}{\mathfrak{p_i}_\mathfrak{m}} \ne 0.
    \]
    Then $N_\mathfrak{m} \ne 0$, i.e, $\mathfrak{m} \in \Suppca N$. Thus $V(\mathfrak{p}_i)\cap \Specm R \subseteq \Suppca N \cap \Specm R$, for all $i$. Suppose $\Suppca N \cap \Specm R$ finite, then so is $V(\mathfrak{p}_i)\cap \Specm R$ and, by the claim, $N_i/N_{i-1} \simeq R/\mathfrak{p}_i$ is finite dimensional. Thus $\dim N = \sum_{i = 1}^r \dim_\C N_i/N_{i-1}$ is finite.
\end{proof}

\begin{lem}\label{lem:fin-dim-torsion}
    Let $N \in \frkA$ be simple and $\U(\h)$-torsion. Then $N$ must be finite-dimensional.
\end{lem}
\begin{proof}
    Let $\lambda \in \h^*\setminus \Supp \W_1(N)$. Suppose $N_{\mathfrak{m}_\lambda} \ne 0$, which is free by Theorem~\ref{theo:almost-free}. Then there is $v/g \in N_{\mathfrak{m}_\lambda}$ such that no element of $\U(\h)_{\mathfrak{m}_\lambda}$ annihilates it. However, since $N$ is $\U(\h)$-torsion and finitely generated, there exists an element $f \in \U(\h)\setminus \{0\}$ such that $fN = 0$, and in particular, $N_{\mathfrak{m}_\lambda} \ni fv/g = 0$, which is a contradiction. So $f/1 = 0$ in $\U(\h)_{\mathfrak{m}_\lambda}$, i.e, there is $f' \notin \mathfrak{m}_\lambda$ such that $f'f = 0$. This is another contradiction, since $\U(\h)$ is a domain. Thus $N_{\mathfrak{m}_\lambda} = 0$ for all $\lambda \in \h^*\setminus \Supp \W_1(N)$. Therefore
    \[
        \Suppca N \cap \Specm \U(\h) = \{\mathfrak{m} \in \Specm \U(\h) \mid N_{\mathfrak{m}} \ne 0\} \subseteq \Supp \W_1(N).
    \]
    But, by Prop~\ref{prop:finite-dim-Wk}, $\Supp \W_1(N)$ is finite, and thus so is $\Suppca N\cap \Specm \U(\h)$. Then the finite dimension of $N$ follows from Lemma~\ref{lem:finite-support-finite-dim}.
\end{proof}

\begin{cor}\label{cor:simple-tor-free}
    The module $M$ is $\U(\h)$-torsion free.
\end{cor}
\begin{proof}
    It is straightforward to show that the set of all $\U(\h)$-torsion elements of $M$ forms a $\g$-submodule of $M$. Then, since $M$ is simple, it is either $\U(\h)$-torsion free or $\U(\h)$-torsion. But the second option is contradictory to the infinite dimension of $M$ by Lemma~\ref{lem:fin-dim-torsion}, which follows the corollary.
\end{proof}

Infinite dimensional simple admissible weight modules only exist for $\g$ of type A and C. Since the weighting functor assigns a $\U(\h)$-finite module to admissible weight modules, it is natural to expect a similar condition on the type of $\g$ for category $\frkA$. The next results answer positively this question.

\begin{theo}\label{theo:U(h)-free-only-type-A-and-C}
    If there exists an infinite dimensional $\U(\h)$-finite $\g$-module that is $\U(\h)$-torsion free (and in particular, if exists a simple module in $\frkA$ of infinite dimension), then $\g$ must be of type A or C.
\end{theo}
\begin{proof}
    Let $N\in \frkA$ be infinite dimensional and $\U(\h)$-torsion free. Let any $t \in T^*$ and $V = \W(N)[t]$. I claim that $V$ is infinite dimensional. Indeed, if not, there exists $\lambda \in t$ such that $V_\lambda = 0$ and, equivalently, $\m_\lambda N = N$. Then, Nakayama's Lemma implies the existence of $r \in \m_\lambda + 1$, such that $rN = 0$, which contradicts the hypothesis over $N$. Thus $V$ is infinite dimensional. Then, since $V$ is admissible by Prop.~\ref{prop:W(M)-is-admissible}, we my apply Lemma~\ref{lem:first-prop-admW}\ref{lem:first-prop-admW:finite-lenght} to conclude that $V$ must have a simple subquotient of infinite dimension (which is also admissible). Thus $\g$ must be of type A or C, by Lemma~\ref{lem:first-prop-admW}\ref{lem:first-prop-admW:type-a-c}.
\end{proof}

For the next statement, assume $\g = \lie{sl}(n+1)$ and recall the class of parabolic subalgebras $\{\p_\bb,\tau(\p_\bb)\mid \bb \in (\C\setminus\{0\})^n\}$ given in Lemma~\ref{lem:dec-sln-h+p}. Moreover, recall the definition of the parabolic BGG category $\cO^\p$, also given in Section~\ref{sec:parabolic-induced} (see also \cite[Chap.~9]{Hum08}).
\begin{cor}\label{cor:parabolic-bgg-O^p-if-only-if}
    Let $\mathfrak{q}$ be a parabolic subalgebra of maximal dimension (i.e. $\dim \mathfrak{q} = \dim \lie{sl}(n+1) - n$). Then $\cO^\mathfrak{q}\subseteq \frkA$ if and only if $\mathfrak{q} \in \{\p_\bb,\tau(\p_\bb)\mid \bb \in (\C\setminus\{0\})^n\}$.
\end{cor}
\begin{proof}
    If $\mathfrak{q}$ is in such list, then Lemma~\ref{lem:dec-sln-h+p} implies $\cO^\mathfrak{q}\subseteq \frkA$.

    Now, suppose $\cO^\mathfrak{q}\subseteq \frkA$ and that $\h \cap \mathfrak{q} \ne 0$. There is an infinite dimensional simple module $N \in \cO^\mathfrak{q} \cap \frkA$. Since $N$ is simple in the parabolic BGG category $\cO^\mathfrak{q}$, there exists a finite dimensional $\mathfrak{q}$-module $V$ such that $N$ is a quotient of $\U(\g)\otimes_{\U(\mathfrak{q})}V$. Let $\widetilde V \subseteq N$ be the image of $V$ under such quotient. Then $\widetilde V$ is non-zero and finite dimensional. Furthermore, $\mathfrak{q}\widetilde V \subseteq \widetilde V$. Then let $h \in \h \cap \mathfrak{q}$ non-zero. Since $h\widetilde V \subseteq \widetilde V$ and $V$ is finite dimensional, there must exist $v \in \widetilde V \subseteq N$ and $p(h)\in \C[h]$ such that $p(h)v = 0$. Hence $N$ is not $\U(\h)$-torsion free, contradicting Corollary~\ref{cor:simple-tor-free}. Therefore $\mathfrak{q} \in \{\p_\bb,\tau(\p_\bb)\mid \bb \in (\C\setminus\{0\})^n\}$ by Lemma~\ref{lem:dec-sln-h+p}.
\end{proof}

%\newpage

%%%%%%%%%%%%%%%%%%%%%%%%%%%%%%%%%%%%%%%%%%%%%%%%%%%%%%%%%%%%%%%%%%%
% Section - Weighting Functor and almost coherent families
%%%%%%%%%%%%%%%%%%%%%%%%%%%%%%%%%%%%%%%%%%%%%%%%%%%%%%%%%%%%%%%%%%%
%%%%%%%%%%%%%%%%%%%%%%%%%%%%%%%%%%%%%%%%%%%%%%%%%%%%%%%%%%%%%%%%%%%
%%%%%%%%%%%%%%%%%%%%%%%%%%%%%%%%%%%%%%%%%%%%%%%%%%%%%%%%%%%%%%%%%%%
% U(h)-finite modules
%%%%%%%%%%%%%%%%%%%%%%%%%%%%%%%%%%%%%%%%%%%%%%%%%%%%%%%%%%%%%%%%%%%
%%%%%%%%%%%%%%%%%%%%%%%%%%%%%%%%%%%%%%%%%%%%%%%%%%%%%%%%%%%%%%%%%%%

\section{Weighting functor and almost-coherent family}\label{sec:weigting-functor-and-a.c.f}

Nilsson showed in \cite{N16} that $\W(M)$ is a coherent family for any $\U(\h)$-free $\g$-module. In this section, we investigate if the same happens for modules of $\frkA$. The difficulty here lies on the fact that such class of modules may not be $\U(\h)$-free if its central character is not integral-regular (see Theo.~\ref{theo:almost-free}). Fortunately, since a simple $\U(\h)$-finite module $M$ is not locally free in at most a finite set of points $\Supp \W_1(M) \subseteq \h^*$, we show that the weight module $\W(M)$ fails to be a coherent family exactly in the points $\Supp \W_1(M)$. This will lead us to the definition of \emph{almost-coherent families}.

\subsection{Geometry of $\U(\h)$-torsion free modules and almost-coherent family}
We denote by $\Frac \U(\h)$ the ring of fraction of $\U(\h)$ and, for any $\U(\h)$-module $M$, 
\[
    \rank M = \dim_{\Frac \U(\h)} M \otimes_{\U(\h)} \Frac \U(\h)
\]
denote its rank. For every element $f \in \U(\h)$, let $D(f) =\{\m \in \Specm \U(\h)\mid f \in \m\}$ be the Zariski open subset of $\Specm \U(\h)$ associated to $f$.

\begin{lem}\label{lem:suitable-basis-of-M_m}
    Let $M$ be a $\U(\h)$-torsion free module and $\mathfrak{m}_0 \in \Specm \U(\h)$ such that $M_{\mathfrak{m}_0}$ is $\U(\h)_{\mathfrak{m}_0}$-free. Let $r = \rank M$. There exist $f_0 \in \U(\h)$ and $\{v_1,\dotsc, v_r\} \in M$ satisfying:
    \begin{enumerate}[label=(\roman*),ref=(\roman*)]
        \item $\mathfrak{m}_0 \in D(f_0)$;
        \item $\{v_1,\dotsc, v_r\}$ is a $\Frac \U(\h)$-basis for $M \otimes_{\U(\h)} \Frac\U(\h)$;
        \item for all $\mathfrak{m} \in D(f_0)$, we have that $M_{\mathfrak{m}}$ is $\U(\h)_{\mathfrak{m}}$-free with basis $\{v_1,\dotsc, v_r\}$;
        \item $\{v_1,\dotsc, v_r\}$ is a $\U(\h)_{f_0}$-basis of $M_{f_0}$.
    \end{enumerate}
\end{lem}
\begin{proof}
    Let $\bar B = \{\bar v_1,\dotsc \bar v_r\}$ be a $\Frac \U(\h)$-basis of $M \otimes_\U(\h) \Frac \U(\h) = M_{(0)}$. Write $\bar v_i = v_i / g_i$, where $v_i \in M$ and $g_i \in \U(\h)\setminus{\mathfrak{m}_0}$.
    
    We claim that $B = \{v_1,\dotsc, v_r\}$ is a $\Frac \U(\h)$-basis of $M_{(0)}$. Indeed $\Frac \U(\h) B = M_{(0)}$, since $\bar B$ generates the $\Frac \U(\h)$-module  $M_{(0)}$. Further it is $\Frac \U(\h)$-linear independent: if $\sum \alpha_iv_i = 0$, then $0 = \sum \alpha_ig_iv_i/g_i = \sum \alpha_ig_i\bar v_i$ which implies that $\alpha_ig_i = 0$, for all $i$, and thus $\alpha_i = 0$ (since $\Frac \U(\h)$ is a field). So $B$ satisfies the second item of the Lemma.

    Further, we claim that $B$ is also a $\U(\h)_{\mathfrak{m}_0}$-basis for $M_{\mathfrak{m}_0}$. Clearly $B$ is $R_{\mathfrak{m}_0}$-linear independent, since it is over $\Frac \U(\h)$. Suppose now that exists $m' \in M_{\mathfrak{m}_0}\setminus \U(\h)_{\mathfrak{m}_0}B$. Then $B \cup \{m'\}$ is $\U(\h)_{\mathfrak{m}_0}$-linear independent. However, this is a contradiction since $\rank M_{\mathfrak{m}_0} = \rank M$ and the rank of a module is the maximum number of linear independent terms. Thus $M_{\mathfrak{m}_0}$ is $\U(\h)_{\mathfrak{m}_0}$-generated by $B$, and the claim is proved. 

    Now let $m_1,\dotsc, m_s \in M$ be a set of generators. Since $M$ is $\U(\h)$-torsion free, there is a natural monomorphism of $\U(\h)$-modules $M \hookrightarrow M_{\mathfrak{m}_0}$, and since $B$ is a $\U(\h)_{\mathfrak{m}_0}$-basis of $M_{\mathfrak{m}_0}$, there exists $f_0 \in \U(\h)\setminus{\mathfrak{m}_0}$ (which, in particular, satisfies the first item of the Lemma) such that
    \[
        m_i = \dfrac{1}{f_0} \sum_{j=0}^{r}\beta_{i,j} v_j,
    \]
    where $\beta_{i,j} \in R$, for all $1 \le i \le s$ and $1 \le j \le r$.
    In particular $M \subseteq \dfrac{1}{f_0}\U(\h) B$, which implies that
    \[
        M_{\mathfrak{m}} \subseteq \dfrac{1}{f_0}\U(\h)_{\mathfrak m} B = \U(\h)_{\mathfrak m} B,
    \]
    for all $\mathfrak{m} \in D(f_0)$. But $\U(\h)_{\mathfrak m} B \subseteq M_{\mathfrak{m}}$ and $B$ is $\Frac \U(\h)$-linear independent. Thus $B$ is a $\U(\h)_{\mathfrak{m}}$-basis of $M_{\mathfrak{m}}$, for all $\mathfrak{m} \in D(f_0)$, i.e, $B$ satisfies the third item of the Lemma.

    Finally, $M_{f_0} = f_0 M_{f_0} \subseteq \U(\h)_{f_0} B \subseteq M_{f_0}$. Thus $M_{f_0} = \U(\h)_{f_0} B$, and since $B$ is $\Frac \U(\h)$-linear independent, it follows that $f_0$ and $B$ satisfies the last item, and the Lemma is proved.
\end{proof}

The next goal is to study the action of $\A = \U(\g)_0$ (the commutant of $\U(\h)$ in $\U(\g)$) in the localization of a module $M\in \frkA$. Before doing it, we may extend the $\g$-action of a module $M\in \frkA$ to the module $M \otimes_{\U(\h)} \Frac \U(\h)$.

Fix $\{h_1,\dotsc,h_n\}$ a basis of $\h$ and, for every $\lambda \in \h^*$, define the automorphism of rings
\begin{equation}\label{def:sigma_lambda}
    \function{\sigma_\lambda}{\U(\h)}{\U(\h)},\quad h_i \mapsto h_i - \lambda(h_i).
\end{equation}
Then, given a $\g$-module $M$, we have that
\[
    x_\alpha (hm) = (h x_\alpha - [h,x_\alpha])m = (h - \alpha(h))x_\alpha m = \sigma_\alpha(h)x_\alpha m,
\]
for all $h \in \h^*$, $\alpha \in \Delta$, $m \in M$, and weight vector $x_\alpha \in \g_\alpha$. And by induction at the degree of polynomials in $\U(\h) \simeq \C[h_1,\dotsc,h_n]$ it can be proved that
\begin{equation}\label{eq:dif-op}
    x_\alpha p m = \sigma_\alpha(p)x_\alpha m, \quad p \in \U(\h), \alpha \in \Delta, m \in M, x_\alpha \in \g_\alpha.
\end{equation}
Furthermore, a straightforward computation shows that
\begin{equation}\label{eq:dif-op-comm}
    \sigma_\alpha \sigma_\beta = \sigma_{\alpha + \beta} = \sigma_\beta \sigma_\alpha, \quad \forall \alpha,\beta \in \h^*.
\end{equation}
With these equations, we may prove:

\begin{lem}\label{lem:g-mod-in-O-localization}
    Let $M \in \frkA$. Then the $\U(\h)$-module $M_{(0)} = M \otimes_{\U(\h)} \Frac \U(\h)$ becomes a $\g$-module by setting
        \[
            x_\alpha \cdot (m \otimes f/g) = x_\alpha m \otimes \sigma_\alpha(f)/\sigma_\alpha(g),
        \]
        for all $\alpha \in \Delta$, $m \in M$, $f/g \in \Frac \U(\h)$ and weight vector $x_\alpha \in \g_\alpha$.
\end{lem}

Given a Zariski-open set $U \subseteq \C^n$, writes by $\mathcal{O}_{\C^n}(U)$ the set of regular functions in $U$ \cite{Harr92}. Then, by identifying $\U(\h) = \C[h_1,\dotsc,h_n] = \mathcal{O}_{\C^n}(\C^n)$, we have that $\mathcal{O}_{\C^n}(D(f)) = \C[h_1,\dotsc,h_n]_{f}$, for all $f \in \U(\h)\setminus\{0\}$, where $D(f)$ denotes the corresponding basic Zariski-open set.

\begin{lem}\label{lem:regular-trace}
    Let $M \in \frkA$ which is $\U(\h)$-torsion free. Let $\mathfrak{m}_0 \in \Specm \U(\h)$ such that $M_{\mathfrak{m}_0}$ is $\U(\h)_{\mathfrak{m}_0}$-free and let $f_0$ satisfying Lemma~\ref{lem:suitable-basis-of-M_m}. Then, for all $u \in \A$, the map
    \[
        \psi_{\mathfrak{m}_0} \colon D(f_0) \ni \lambda \mapsto \tr \restr{u}{M/\mathfrak{m}_\lambda M},
    \]
    is a regular function in $D(f_0)$,
\end{lem}
\begin{proof}
    Let $u \in \A$. From Lemma~\ref{lem:g-mod-in-O-localization} we have that the action of $u$ in $M$ extends to $M_{(0)}$. Furthermore, since $u$ lies in $\A$, the definition of the action of $\g$ in $M_{(0)}$ implies that $uM_{f_0} \subseteq M_{f_0}$.
   
    Fix a suitable $\U(\h)_{f_0}$-basis $B = \{v_1,\dotsc,v_r\}\subseteq M$ of $M_{f_0}$ obtained from Lemma~\ref{lem:suitable-basis-of-M_m}. Since $u$ commutes with $\U(\h)_{f_0}$ and $M_{f_0}$ is $\U(\h)_{f_0}$-free, we have that $u \in \End_{\U(\h)_{f_0}}(M_{f_0})$ and that exists $\{\alpha_{i,j}\}_{1 \le i,j \le r}$ in $\U(\h)_{f_0}$ such that $u v_i = \sum_{i,j} \alpha_{i,j}v_j$. Since $B$ is a $\U(\h)_{\mathfrak{m}}$-basis of $M_{\mathfrak{m}}$, for all $\mathfrak{m} \in D(f_0)$, we obtain 
    \[
        M/\mathfrak{m}M \simeq M_{\mathfrak{m}}/\mathfrak{m} M_{\mathfrak{m}} \simeq \bigoplus_{i=1}^{r} \U(\h)/\mathfrak{m},
    \]
    and $u$ acts in $M/\mathfrak{m}M$ via the matrix $(\alpha_{i,j} + \mathfrak{m}M)_{i,j}$.

    Then, thought the identification $\C^n \ni \lambda \mapsto (h_1 - \lambda_1,\dotsc, h_n-\lambda_n) \in \Specm \U(\h)$ we have that $u$ acts in $M/\mathfrak{m}_\lambda M$ via the matrix $(\alpha_{i,j}(\lambda))_{i,j}$, for all $\lambda$ such that $f_0(\lambda) \ne 0$.
    Thus
    \[
        \tr{\restr{u}{M/\mathfrak{m}_\lambda M}} = \alpha_{1,1}(\lambda) + \dotsb + \alpha_{r,r}(\lambda).
    \]
    Therefore $\psi_{\mathfrak{m}_0} \in \C[h_1,\dotsc,h_n]_{f_0}$, and the lemma is proved.
\end{proof}

\begin{theo}\label{theo:W(M)-quasi-coh-fam}
    Let $M \in \frkA$ be simple of infinite dimension. Let $r = \rank M$, $\mathcal{M} = \W(M)$ and $F = \Supp \W_1(M)$. Then
    \begin{enumerate}[label=(\roman*),ref=(\roman*)]
        \item $\dim \mathcal{M}_\lambda = r$, for all $\lambda \in \h^* \setminus F$; and \label{theo:W(M)-quasi-coh-fam:dim}
        \item For all $u \in \A$, the function $\psi_u \colon \lambda \mapsto \tr{\restr{u}{\mathcal{M}_\lambda}}$ is regular in $\lambda \in \h^* \setminus F$.\label{theo:W(M)-quasi-coh-fam:trace-regular}
    \end{enumerate}
    In particular, if $n \ge 2$ we may replace item \ref{theo:W(M)-quasi-coh-fam:trace-regular} by
    \begin{enumerate}[resume, label=(\roman*),ref=(\roman*)]
        \item For all $u \in \A$, the function $\psi_u \colon \lambda \mapsto \tr{\restr{u}{\mathcal{M}_\lambda}}$ is \emph{polynomial} in $\lambda \in \h^* \setminus F$.\label{theo:W(M)-quasi-coh-fam:trace-polynomial}
    \end{enumerate}
\end{theo}
\begin{proof}
    By Theorem~\ref{theo:almost-free}, we have that $M_{\mathfrak{m}_\lambda}$ is $\U(\h)_{\mathfrak{m}_\lambda}$-free if, and only if, $\lambda \in \h^* \setminus F$. Notice that, as $\U(\h)$-modules,
    \[
        \mathcal{M}_\lambda= M/\mathfrak{m}_\lambda M \simeq M_{\mathfrak{m}_\lambda}/\mathfrak{m}_\lambda M_{\mathfrak{m}_\lambda} \simeq \bigoplus_{i=1}^{r} \U(\h)/\mathfrak{m} \simeq \C^r,
    \]
    whenever $M_{\mathfrak{m}_\lambda}$ is $\U(\h)_{\mathfrak{m}_\lambda}$-free. Furthermore, by Corollary~\ref{cor:simple-tor-free}, $M$ is $\U(\h)$-torsion free.
    Then, by Lemmas~\ref{lem:suitable-basis-of-M_m} and \ref{lem:regular-trace}, it follows that, for all $\lambda \in \h^* \setminus F$, there is an element $f_{\lambda} \in \U(\h)$ such that
    \begin{enumerate}[label=(\roman*),ref=(\roman*)]
        \item $\dim \mathcal{M}_\mu = r$, for all $\lambda \in D(f_\lambda)$; and
        \item For all $u \in \A$, the function $\psi_u \colon D(f_\lambda) \ni \lambda \mapsto \tr{\restr{u}{\mathcal{M}_\lambda}}$ is an element of $\U(\h)_{f_0}$.
    \end{enumerate}

    One can see that $\h^*\setminus F \subseteq \bigcup_{\lambda \in \h^*\setminus F} D(f_\lambda)$. Suppose exists an element $\mu \in F\cap D(f_\lambda)$, for some $\lambda \in \h^*\setminus F$. Then $M_{\mathfrak{m}_\mu}$ is $\U(\h)_{\mathfrak{m}_\mu}$-free, by Lemma~\ref{lem:suitable-basis-of-M_m}. But this is a contradiction, by Theorem~\ref{theo:almost-free}. Thus $\h^*\setminus F = \bigcup_{\lambda \in \h^*\setminus F} D(f_\lambda)$ and items \ref{theo:W(M)-quasi-coh-fam:dim} and \ref{theo:W(M)-quasi-coh-fam:trace-regular} follow.

    Notice that, by Prop.~\ref{prop:finite-dim-Wk}, $\Supp \W_1(M) \subseteq \h^*$ is a cofinite subset. Then the last item follows from the fact that a regular functions defined in a cofinite subset of an affine variety of dimension at least two are in fact polynomial. 
\end{proof}

The previous theorem gives us the following idea for a more general definition of coherent families:

\begin{defin}\label{def:U-coherent-family}
    Let $U \subseteq \h^*$ be a subset with finite complement. A weight module $\M$ is called a \emph{$U$-coherent family} of degree $d = \deg \M$ if satisfies
    \begin{enumerate}[label=(\roman*),ref=(\roman*)]
        \item $\dim \mathcal{M}_\lambda = d$, for all $\lambda \in U$,
        \item For all $u \in \A$, the function $\bigtr_u \colon \lambda \mapsto \tr{\restr{u}{\mathcal{M}_\lambda}}$ is polynomial in $U$.
    \end{enumerate}
    We say that $\M$ is an \emph{almost-coherent family} if it is a $U$-coherent family, for some cofinite subset $U \subseteq \h*$.
\end{defin}

\begin{rem}
    Notice that the definition of $\h^*$-coherent family recovers the usual coherent family defined by O. Mathieu in \cite{Mat00}.
\end{rem}

Then, Theorem~\ref{theo:W(M)-quasi-coh-fam} implies:

\begin{cor}\label{Cor:W(M)-U-coherent-family}
    Let $M \in \frkA$ that is $\U(\h)$-torsion free and assume $\rank \g \ge 2$. Then $\W(M)$ is an $\h^* \setminus{\Supp \W_1(M)}$-coherent family of degree equal to $\rank M$. 
\end{cor}

Recall that $\frkA_\fl$ is the subcategory of $\frkA$ consisting of $\U(\h)$-finite modules of finite length.

\begin{theo}\label{Theo:W(M)-U-coherent-family}
    Assume $\rank \g \ge 2$. The weighting functor $\W$ assign any infinite dimensional object of $\frkA_\fl$ to an almost-coherent family.
\end{theo}
\begin{proof}
    Let $N,M \in \frkA_\fl$ of infinite dimension such that $N\subseteq M$, $M/N$ is simple and $\W(N)$ is $U$-coherent family, for cofinite subset $U\subseteq \h^*$. Next, we prove that $\W(M)$ is an almost coherent family.
    
    Consider the exact sequence $0 \rightarrow N \rightarrow M \rightarrow M/N \rightarrow 0$ in $\frkA_\fl$. Then the weighting functors induce the exact sequence
    \[
        \W_1(M/N) \rightarrow \W(N) \rightarrow \W(M) \rightarrow \W(M/N) \rightarrow 0.
    \]
    By Prop.~\ref{prop:finite-dim-Wk}, the module $\W_1(M/N)$ is finite dimensional. Thus, we have an exact sequence
    \begin{equation}\label{eq:Theo:W(M)-U-coherent-family}
        0 \rightarrow \W(N)/W \rightarrow \W(M) \rightarrow \W(M/N) \rightarrow 0,
    \end{equation}
    where $W$ is the image of $\W_1(M/N) \rightarrow \W(N)$.
    Moreover, since $W$ is finite dimensional, the module $\W(N)/W$ is a $U\setminus\Supp W$-coherent family.
    
    Now, if $M/N$ is finite dimensional, it is a weight module and then $\W(M/N) = M/N$, by Lemma~\ref{lem:prop-Wk}. Then $\dim (\W(N)/W)_\lambda = \dim \W(M)_\lambda$ for all $\lambda \notin \Supp M/N$, and $\W(M)$ is a $U\setminus(\Supp W\cup \Supp M/N)$-coherent family. 
    
    If $M/N$ is infinite dimensional, it is $\U(\h)$-torsion free by Corollary~\ref{cor:simple-tor-free}. In particular, by Corollary~\ref{Cor:W(M)-U-coherent-family}, there exists a cofinite subset $U'\subseteq\h^*$ such that $\W(M/N)$ is a $U'$-coherent family.
    Then the exact sequence \eqref{eq:Theo:W(M)-U-coherent-family} implies, for all $\lambda \in U\cap U'$:
    \begin{itemize}
        \item $\dim \W(M)_\lambda = \dim\W(N)_\lambda + \dim\W(M/N)_\lambda$;
        \item $\tr \restr{a}{\W(M)_\lambda} = \tr \restr{a}{(\W(N)\oplus\W(M/N))_\lambda} = \tr \restr{a}{\W(N)_\lambda} + \tr \restr{a}{\W(M/N)_\lambda}$, for all $a \in \A$.
    \end{itemize}
    And the fact that, for all $a \in \A$, the maps $\lambda \mapsto \tr \restr{a}{\W(N)_\lambda}$ and $\lambda \mapsto \tr \restr{a}{\W(M/N)_\lambda}$ are polynomial in $U$ and $U'$, respectively, implies that $\W(M)$ is a $U\cap U'$-coherent family.

    Finally, the general case of the theorem follows by the claim above together with induction on the length of objects of $\frkA_\fl$.
\end{proof}

\subsection{Almost-equivalence and semi-simple almost-coherent family}\label{sec:almost-equiv}

For a weight module $M$ define the \emph{trace map associated} to $M$ to be the function
\[
    \bigtr^M \colon \h^* \times \A \rightarrow \C,\quad \bigtr^M(\lambda,u) = \tr \restr{u}{M_\lambda}.
\]
\begin{defin}
    Two weight modules $M$ and $N$ are \emph{almost-equivalent}, and write $M\sim N$, if there exists a cofinite subset $U\subseteq \h^*$ such that $\restr{\bigtr^M}{U\times \A} = \restr{\bigtr^N}{U\times \A}$. When it is of interest to emphasize such a cofinite subset $U$, we write $M\sim_U N$.
\end{defin}

It is straightforward to check that $\sim$ defines an equivalence relation. 
Moreover, almost equivalence preserves generalized central character. Precisely:

\begin{lem}\label{lem:almost-equiv-central-char}
    Let $M$ and $N$ be two almost-equivalent weight modules with generalize central characters $\chi_M$ and $\chi_N$, respectively. Then $\chi_M = \chi_N$.
\end{lem}
\begin{proof}
    Let $\lambda \in \h^*$ be such that $\restr{\bigtr^M}{\{\lambda\}\times \A} =\restr{\bigtr^N}{\{\lambda\}\times \A} \ne 0$. Then, for every $z \in Z(\g)$ we have
\[
    \chi_M(z)\dim M_\lambda = \tr\restr{z}{M_\lambda} =\bigtr^M(\lambda,z) = \bigtr^N(\lambda,z) = \tr\restr{z}{N_\lambda} = \chi_N(z)\dim N_\lambda.
\]
Since $\chi_M(1) = 1 =\chi_N(1)$, it follows that $\dim M_\lambda = \dim N_\lambda$ and thus $\chi_M = \chi_N$.
\end{proof}

Furthermore, the next lemma shows that we can focus our study for semi-simple weight modules.
\begin{lem}\label{lem:module-almost-equiv-to-semi-simply}
    Let $0 \rightarrow L \rightarrow M \rightarrow N \rightarrow 0$ be an exact sequence of weight modules. Then $M$ is almost equivalent to the direct sum $L \oplus N$. 
    
    Furthermore, suppose $M$ has finite length and let $M^\semisim$ be its semi-simplification, i.e, be the semi-simple weight modules with same composition series of $M$. Then $M\sim M^\semisim$.
\end{lem}
\begin{proof}
    For the first claim, let $\lambda \in \h^*$ and $a \in \A$. The vector space $M_\lambda$ is isomorphic to $L_\lambda \oplus N_\lambda$. Moreover, $L_\lambda$ and $N_\lambda$ are stable under the action of $a$. Then we have
    \[
        \bigtr^M(\lambda,a) = \tr \restr{a}{M_\lambda} = \tr \restr{a}{L_\lambda} + \tr \restr{a}{N_\lambda} = \tr\restr{a}{(L \oplus N)_\lambda} = \bigtr^{L\oplus N}(\lambda,a),
    \]
    and the claim follows. The second claim follows by induction on the length of the module.
\end{proof}

In \cite[Lemma~2.3]{Mat00}, it is proved that two \emph{semi-simple} weight modules have same trace maps if, and only if, they are \emph{isomorphic}. Next we show that the almost-equivalence of two semi-simple weight modules translates as ``isomorphism up to finite dimensional modules''. More precisely:

\begin{lem}\label{lem:U-approx}
    Two semi-simple weight modules $M$ and $N$ are almost-equivalent if, and only if, exist finite dimensional $\g$-modules submodules $V_M\subseteq M$ and $V_N\subseteq N$ such that $M/V_M$ and $N/V_N$ are isomorphic. 
    
    Emphasizing the cofinite subset in the almost-equivalence, we may state that $M\sim_U N$, for cofinite $U\subseteq \h^*$ if, and only if, $M/V_M \simeq N/V_N$ for some finite dimensional $\g$-modules $V_M\subseteq M$ and $V_N\subseteq N$ with supports contained in the complement of $U$, i.e., $\Supp V_M,\,\Supp V_N \subseteq \h^*\setminus U$.
\end{lem}
\begin{proof}
    Let $\mathcal{S}$ be a set of representatives of isomorphism class of simple weight $\g$-modules. Consider two isomorphism of $\g$-modules
    \[
        \function{\psi}{M}{\bigoplus_{L\in S} L^{\oplus r_L}}\quad \text{and}\quad \function{\varphi}{\bigoplus_{L\in S} L^{\oplus s_L}}{N},
    \]
    where $r_L,s_L \in \N$ for all $L \in S$, and define the submodules $M' \subseteq M$ and $N' \subseteq N$ such that
    \[
        \psi(M') = \bigoplus_{L \in \mathcal{S}} L^{\oplus \min\{r_L,s_L\}} = \varphi^{-1}(N').
    \]
    Clearly $\function{\varphi^{-1}\circ\psi}{M'}{N'}$ define an isomorphism. And, by setting $V_M = M/M'$ and $V_N = N/N'$, the semi-simplicity of $M$ and $N$ implies that $M' \simeq M/V_M$ and $N' \simeq N/V_N$. Therefore it is enough to prove that $M\sim N$ if and only if $V_M$ and $V_N$ are finite dimensional.

    Notice that $M_\lambda = M'_\lambda$ and $N_\lambda = N'_\lambda$ for every $\lambda$ outside $\Supp V_M\cup\Supp V_N$. Then, since $M'$ and $N'$ are isomorphic as $\g$-modules they have the same trace. This implies that $\restr{\bigtr^M}{U\times \A} = \restr{\bigtr^N}{U\times \A}$, where $U = \h^* \setminus (\Supp V_M\cup\Supp V_N)$. Then, if $V_M$ and $V_N$ are finite dimensional, the subset $U$ is cofinite and thus $M$ is almost-equivalent to $N$ (precisely, $M\sim_U N$).
    
    Reciprocally, assume $M \sim_U N$ for a cofinite subset $U\subseteq \h^*$.
    For any $L \in \frkW$ simple, \cite[Lemma~2.2]{Mat00} implies that the multiplicities of $L$ in $M$ and $N$ are the same as the multiplicities of $L_\lambda$ in $M_\lambda$ and $N_\lambda$ as $\A$-module, for all $\lambda \in \Supp L$. Theses last numbers are determined by $\restr{\bigtr^M}{\{\lambda\}\times \A}$ and $\restr{\bigtr^N}{\{\lambda\}\times \A}$, respectively. Therefore, for every simple weight module $L$ with support intersecting $U$, we have that $L^{\oplus r} \hookrightarrow M$ if and only if $L^{\oplus r} \hookrightarrow N$.   
    Thus $\dim M_\lambda = \dim M'_\lambda = \dim N'_\lambda = \dim N_\lambda$, for all $\lambda \in U$. This implies that $\Supp V_M$ and $\Supp V_N$ are contained in $\h^*\setminus U$ and therefore have finite dimension. 
\end{proof}

Let $M,\,N \in\frkW^\chi$ be two semi-simple weight modules with \emph{non-integral} or \emph{singular} generalized central character. Suppose that $M$ is almost equivalent to $N$ and consider finite dimensional submodule modules $V_M\subseteq M$ and $V_N\subseteq N$ such as Lemma~\ref{lem:U-approx}. Then $V_M$ and $V_N$ have central character $\chi$. But non-trivial finite dimensional modules must be direct sum of modules with integral regular central character. Hence $V_M = V_N = 0$. Therefore $M \simeq M/V_M \simeq N/V_N \simeq N$. Moreover, it is clear that isomorphism of modules implies almost-equivalency. Hence, the following holds:
\begin{cor}\label{cor:u-aprox-non-intregral-singular}
    Let $\chi \in \X$ be a non-integral or a singular central character. Then, two semi-simple weight modules with central character $\chi$ are isomorphic if and only if they are almost equivalent.
\end{cor}

Another interesting property of $\sim$ is that this equivalence is preserved by the translation functor when it is an equivalence of category.

\begin{lem}\label{lem:u-aprox-translation}
    Let $(\lambda, \mu) \in \Upsilon$ and let $M, N$ be two semi-simple weight modules in $\g\mdfg^{\chi_\mu}$. If $M \sim N$, then $T_\mu^\lambda M \sim T_\mu^\lambda N$.
\end{lem}
\begin{proof}
    By Lemma~\ref{lem:U-approx}, we can consider finite dimensional submodules $V_M\subseteq M$ and $V_N\subseteq N$ such that $M/V_M \simeq N/V_N$. Write $T = T_\mu^\lambda$. Since $T$ is an equivalence of category, we have
    \[
        TM \simeq T(M/V_M) \oplus T(V_M),\quad TN \simeq T(N/V_N) \oplus T(V_N) \quad \text{and}\quad T(M/V_M)\simeq T(N/V_N).
    \]
    And the fact that translation functors preserve the category of finite dimensional $\g$-modules implies that $T(V_M) \subseteq TM$ and $T(V_N) \subseteq TN$ are finite dimensional. Finally, applying Lemma~\ref{lem:U-approx}, we obtain the almost equivalence of $TM$ and $TN$.
\end{proof}

The next goal is to show that a semi-simple almost-coherent family is almost equivalent to a standard semi-simple coherent family. For the argumentation, we make use of \emph{cuspidal modules}, which by definition are admissible weight $\g$-modules on which every weight vector of $\g$ act bijectively. In particular, every cuspidal module $M$ has support equal to exactly one $Q$-coset $t \in T^*$ and $\dim M_\lambda$ is the same for all $\lambda \in t$. This later property allows us, by comparing dimensions, to conclude almost-equivalence from inclusions. More precisely, if $L$ is a cuspidal module and $\M$ is an almost-coherent family containing $L$ and with the same degree, then $\M[\Supp L] \sim L$, since their trace maps $\bigtr^{\M[\Supp L]}(-,\A)$ and $\bigtr^L(-,\A)$ differ only by finitely many points of $\h^*$.

\begin{lem}\label{lem:cuspidal-almost-equiv}
    Let $\M,\mathcal{N}$ be two semi-simple almost-coherent family of degree $d$. If there exists a cuspidal module $L$ of degree $d$ that is contained in both families $\M$ and $\mathcal{N}$, then these families are almost-equivalent.
\end{lem}
\begin{proof}
    Let $U \subseteq \h^*$ be a cofinite such that both $\M$ and $\mathcal{N}$ are $U$-coherent families, and let $t = \Supp L$. Then, since $L\subseteq\M,\mathcal{N}$ and they have the same degree, we have the equality $\M[t] \sim L \sim \mathcal{N}[t]$. Then, there exists a cofinite subset $U'\subseteq t$ such that $\restr{\bigtr^\M}{U'\times \A} = \restr{\bigtr^{\mathcal{N}}}{U'\times \A}$. Since $\restr{\bigtr^\M}{U \times \{a\}}$ and $\restr{\bigtr^{\mathcal{N}}}{U\times \{a\}}$ are polynomial in $U$ for all $a \in \A$, and $U'$ is dense in $U$ (because $t$ is dense and $U'$ is cofinite subset of $t$), it follows that $\restr{\bigtr^\M}{U\times \A} = \restr{\bigtr^{\mathcal{N}}}{U\times \A}$. Thus $\M \sim_U \mathcal{N}$.
\end{proof}

The previous lemma hints that we need to study how to find cuspidal submodules in an almost-coherent family. The approach is the same of \cite[Sec.~5]{Mat00}.
So, for an almost-coherent family $\M$ we define 
\[
    \Sing \M = \{t \in T^* \mid \restr{f_\alpha}{\M[t]} \text { is not injective for some }\alpha \in \Delta\}.
\]
Equivalently, the set $\Sing\M$ is precisely the subset of $T^*$ consisting of cosets $t$ such that $\M[t]$ is not cuspidal. 

As in \cite[Sec.~5]{Mat00}, we define an exotic topology in $T^* = \h^*/\Z\Delta$ in the following way. For any Zariski open $\Omega \subseteq \h^*$, we set $T(\Omega) = \cap_{\mu \in \Z\Delta} (\mu + \Omega)$. It is satisfied that $T(\Omega) \cap T(\Omega') = T(\Omega \cap \Omega')$ for any $\Omega, \Omega' \subseteq \h*$ Zariski open. Then 
\[
    \{T(\Omega)/\Z\Delta \mid \Omega \subseteq \h^* \text{ is Zariski open}\}
\]
form a basis of a topology in $T^*$. We call such (exotic) topology as the \emph{torus topology} in $T^*$.
\begin{lem}[{\citealp[Lemma~5.2]{Mat00}}]\label{lem:T(omega)}
    For any Zariski open subset $\Omega \subseteq \h^*$, the subset $T(\Omega)\subseteq T^*$ is non-empty. And every non-empty open subset of $T^*$ is dense.
\end{lem}

\begin{lem}\label{lem:sing-proper-close}
    The set $\Sing \M$ is a proper closed subset of $T^*$
\end{lem}
\begin{proof}
    Let $U \subseteq \h^*$ be the biggest (by inclusion) cofinite subset such that $\M$ is a $U$-coherent family. Start by noticing that 
    \begin{equation}\label{eq:complement-singM}
        F = \{ \lambda + \Z\Delta \in T^* \mid \lambda + \Z\Delta \subseteq \h^*\setminus U\} \subseteq \Sing \M.
    \end{equation}
    Indeed, if exists $t \in F \cap (T^*\setminus \Sing \M)$, then $\M[t]$ is cuspidal and $\M$ is a $U \cup t$-coherent family, contradicting the hypothesis over $U$
    
    For every $\alpha \in \Delta$, consider the map $\function{p_\alpha}{\h^*}{\C}$ such that $p_\alpha(\lambda) = \restr{e_\alpha f_\alpha}{\M_\lambda}$. I claim that $p_\alpha$ is a non-zero polynomial in $U$. Indeed, the determinant of a given operator $A \in \End_\C(V)$ in a $d$-dimensional space $V$ may be expressed as a polynomial in the variables of the trace of the powers of $A$ (see the Leverrier-Faddeev Characteristic Polynomial Algorithm in \cite{Hou98}). Further, we have that $\restr{\bigtr^{\M}}{U\times \{(e_\alpha f_\alpha)^l\}}$ is polynomial for all $l \ge 1$. Then the claim is proved.
    
    Then 
    \[
        \Omega_\alpha = \{ \lambda \in U \mid \restr{e_\alpha f_\alpha}{\M_\lambda} \text{ is bijective}\} = p_\alpha^{-1}(\C\setminus \{0\}) \cap U
    \]
    is non-empty open subset of $\h^*$, for all $\alpha \in \Delta$.
    Then Lemma~\ref{lem:T(omega)} implies that $\widetilde{\Omega}_\alpha = T(\Omega_\alpha)/\Z\Delta$ is non-empty open and dense subset of $T^*$, for all $\alpha \in \Delta$. Furthermore, for all $t \in \widetilde \Omega \coloneqq \cap_{\alpha \in \Delta} \widetilde{\Omega}_\alpha$ and $\alpha \in \Delta$, it follows that $\restr{e_\alpha f_\alpha}{\M[t]}$ is bijective. Thus 
    \begin{equation}\label{eq:omega-in-set-of-cuspidal}
        \widetilde\Omega \subseteq \{ t \in T^* \mid \restr{f_\alpha,f_{-\alpha}}{\M[t]} \text{ are injective, }\forall\alpha\in\Delta\}. 
    \end{equation}
    Therefore, $T^*\setminus \widetilde\Omega$ contains $\Sing \M$. Suppose that $t \notin \widetilde\Omega$, then there is a weight $\lambda \in t$ such that $\lambda \notin U$ or such that $\restr{e_\alpha f_\alpha}{\M_\lambda}$ is not bijective. Both cases imply that $t\in \Sing\M$, by equations \eqref{eq:complement-singM} and \eqref{eq:omega-in-set-of-cuspidal}, respectively. Thus $T^*\setminus \widetilde\Omega = \Sing \M$. Moreover, $\widetilde{\Omega}$ is a non-empty open dense subset of $T^*$, since it is the intersection of finitely many open dense subsets. Thus the Lemma is proved.
\end{proof}

Finally, we can precisely state how working with semi-simple almost-coherent family is essentially the same as working with semi-simple ($\h^*$-)coherent family.

\begin{theo}\label{theo:U-aprox-ss-U-c.f.}
    Let $U \subseteq \h^*$ be a cofinite subset and let $\M$ be a semi-simple $U$-coherent family. Then there exists a semi-simple coherent family $\mathcal{N}$ that is almost-equivalent to $\M$ (precisely, $\M \sim_U \mathcal{N}$). Furthermore, $\M$ and $\mathcal{N}$ have the same degree.
\end{theo}
\begin{proof}
    By Lemmas~\ref{lem:T(omega)} and \ref{lem:sing-proper-close}, there exists $t \in T^*$ such that $t \subseteq U$ and $\restr{f_\alpha}{\M[t]}$ is injective for all $\alpha \in \Delta$. In particular $\M[t]$ is cuspidal. Then there exist simple cuspidal weight modules $L_1,\dotsc,L_r$ such that $\M[t] \simeq L_1 \oplus \dotsb \oplus L_r$. Let $\mathcal{N} = \cEXT(L_1) \oplus \dotsb \oplus \cEXT(L_r)$. Clearly $\M[t] = \mathcal{N}[t]$, and the theorem follows from Lemma~\ref{lem:cuspidal-almost-equiv}.
\end{proof}

\begin{defin}
    A $U$-coherent family $\M$ is called \emph{irreducible} if $\M_\lambda$ is a simple $\A$-module, for some $\lambda \in \h^*\setminus(U+ \Z\Delta)$.  We simply say that an almost-coherent family is irreducible if it is an irreducible $U$-coherent family, for some cofinite subset $U \subseteq \h^*$.
\end{defin}

\begin{cor}\label{cor:irr-ss-U-c.f.}
    Every irreducible almost-coherent family is almost-equivalent to a unique irreducible semi-simple coherent family.
    %Let $\M$ be an irreducible almost-coherent family. Then there exists an irreducible semi-simple coherent family $\mathcal{N}$ that is almost-equivalent to $\M$.
\end{cor}
\begin{proof}
    Since $\M[t]$ is admissible weight module for all $t \in T^*$, and admissible weight modules have finite length by Lemma~\ref{lem:first-prop-admW}\ref{lem:first-prop-admW:finite-lenght}, we may define the ``semi-simplification'' of $\M$ as the semi-simple weight module $\M^\semisim$ such that $\M[t]$ and $\M^\semisim[t]$ have the same composition series, for all $t \in T^*$. And Lemma~\ref{lem:module-almost-equiv-to-semi-simply} implies $\M\sim \M^\semisim$. Then we may assume that $\M$ is semi-simple.
    
    Now let $U \subseteq \h^*$ be a cofinite subset such that $\M$ is a semi-simple $U$-coherent family and let $\lambda \in \h^*\setminus(U+ \Z\Delta)$ such that $\M_\lambda$ is $\A$-simple. By Theorem~\ref{theo:U-aprox-ss-U-c.f.} and the definition of almost-equivalence, $\M$ is almost equivalent to $\mathcal{N} = \cEXT(L_1) \oplus \dotsb \oplus \cEXT(L_r)$, for some simple cuspidal modules $L_1,\dotsc, L_r$, and $\M_\lambda \simeq \mathcal{N}_\lambda$ as $\A$-module. Then the simplicity of $\M_\lambda$ implies that $r = 1$, i.e, $\mathcal{N} = \cEXT(L_1)$, hence the Corollary. 
\end{proof}

We finish this section by showing how irreducible almost-coherent family behaves with translation functor. 

\begin{lem}\label{lem:T_mu^lambda(EXT(L)}
    Let $(\lambda,\mu)\in \Upsilon$ and let $L \in \frkW$ be a simple admissible of infinite dimension and central character $\chi_\mu$. Let $\M$ be an irreducible almost-coherent family such that $\M \sim \cEXT(L)$. Then $T_\mu^\lambda \M \sim T_\mu^\lambda \cEXT(L)\simeq \cEXT(T_\mu^\lambda L)$.
\end{lem}
\begin{proof}
    By Lemma~\ref{lem:u-aprox-translation}, it is enough to show the isomorphism $T_\mu^\lambda \cEXT(L)\simeq \cEXT(T_\mu^\lambda L)$.
    
    Let $V$ be a finite dimensional $\g$-module and writes $\cE = \cEXT(L)$. We claim that $\cE \otimes V$ is a coherent family. Indeed, for every $\nu \in \h^*$, it is straightforward to show that $\dim(\cE\otimes V)_\nu = \deg L \dim V$. Let $u \in \A$. We have that (using Sweedler's notation)
    \[
        \Delta(u) = \sum u_{(1)} \otimes u_{(2)} + \sum u'_{(1)} \otimes u'_{(2)},
    \]
    where $u_{(1)},u_{(2)} \in \A$ and $u'_{(1)}\otimes u'_{(2)} \in \U(\g)_{\nu'}\otimes\U(\g)_{-\nu'}$, for some $\nu' \in \h^*\setminus \{0\}$.
    Since the matrix of $\restr{u'_{(1)}\otimes u'_{(2)}}{(\cE\otimes V)_\nu}$ has null diagonal, for a suitable basis, it clearly has zero trace. Furthermore,
    \[
        \tr \restr{u_{(1)}\otimes u_{(2)}}{(\cE\otimes V)_\nu} = \sum_{\nu' \in \Supp V} \tr \restr{u_{(1)}}{\cE_{\nu - \nu'}} \tr \restr{u_{(2)}}{V_{\nu'}},
    \]
    which is polynomial in $\nu$, since $\cE$ is a coherent family. Therefore $\tr \restr{u}{(\cE \otimes V)_\nu}$ is polynomial in $\nu \in \h^*$.

    Thus, for $V = L(\widehat{\lambda - \mu})$, there exist simple admissible weight modules $L_1,\dots,L_k$ of infinite dimension, such that
    \[
        \left( \cE \otimes V\right)^\semisim \simeq \cEXT(L_1)\oplus \dotsb \oplus \cEXT(L_k),
    \]
    where $(-)^\semisim$ denotes the ``semi-simplification'' we defined in the proof of Corollary~\ref{cor:irr-ss-U-c.f.}.
    Since $(-)^{\chi_\lambda}$ and $T_\mu^\lambda$ commute with semi-simplification, $\cE$ is semi-simple and $\cEXT(L_i)$ has the same central character of $L_i$, we obtain:
    \[
        T_\mu^\lambda \cE \simeq \left(\left( \cE \otimes V\right)^\semisim\right)^{\chi_\lambda} \simeq \bigoplus_{\chi_{L_i} = \chi_\lambda} \cEXT(L_i).
    \]
    
    Writes $\cL_i = \cEXT(L_i)$. By Lemma~\ref{lem:sing-proper-close} (and the fact that non-empty open subsets of $T^*$ are dense), we can choose $t \in T^*$ such that $\cE[t]$ is cuspidal. Moreover, cuspidal modules are preserved by translation functor, by \cite[Lemma~2.7]{GS10}. Thus $T_\mu^\lambda(\cE[t]) \subseteq T_\mu^\lambda(\cE)$ is simple and cuspidal module and thus isomorphisc to $(T_\mu^\lambda\cE)[t + (\lambda - \mu)]$. But $(T_\mu^\lambda\cE)[t + (\lambda - \mu)] = \bigoplus_{\chi_{L_i} = \chi_\lambda} \cL_i[t + (\lambda - \mu)]$. So this sum must have only one term, let say $\cL_1[t + (\lambda - \mu)]$. Hence $T_\mu^\lambda \cE = \cL_1$.

    Finally, notice that $T_\mu^\lambda L \subseteq T_\mu^\lambda \cE$ and that $T_\mu^\lambda L$ is an infinite dimensional simple admissible weight module with central character $\chi_\mu$, as $T_\mu^\lambda$ is equivalence of category. Then from Prop.~\ref{prop:EXT(L)}~\ref{prop:EXT(L):uniq-ss-irr-c-f} we conclude $T_\mu^\lambda \cE \simeq \cEXT(T_\mu^\lambda L)$.
\end{proof}

%\newpage

%%%%%%%%%%%%%%%%%%%%%%%%%%%%%%%%%%%%%%%%%%%%%%%%%%%%%%%%%%%%%%%%%%%
% Section - Classification
%%%%%%%%%%%%%%%%%%%%%%%%%%%%%%%%%%%%%%%%%%%%%%%%%%%%%%%%%%%%%%%%%%%
%%%%%%%%%%%%%%%%%%%%%%%%%%%%%%%%%%%%%%%%%%%%%%%%%%%%%%%%%%%%%%%%%%%
% Section - Classification
%%%%%%%%%%%%%%%%%%%%%%%%%%%%%%%%%%%%%%%%%%%%%%%%%%%%%%%%%%%%%%%%%%%

\section{Category $\frkA^{\irr}$ and its simple objects}\label{sec:almost-comp-class}

Theorem~\ref{Theo:W(M)-U-coherent-family} shows that, after applying the weighting functor on an $\U(\h)$-finite module of infinite dimension and finite length, we obtain an almost-coherent family. So in this section we study which of this $\U(\h)$-finite modules correspond to irreducible families. In other words, we will study the simple objects of the full subcategory $\frkA^{\irr}$ of $\frkA$ defined by
\[
    \frkA^{\irr} = \{M \in \frkA\mid \W(M)\text{ is an irreducible almost-coherent family}\}.
\]

For the case $\g = \lie{sl}(n+1)$ and $n \ge 2$, the complete classification of these modules is left open. However, combining results from \cite{GN22}, \cite{N15} and \cite{BL01} we give the classification when considering modules with non-integral or integral singular central character. For $\g = \lie{sp}(2n)$, we conclude the classification with no restriction on the central character. 

Finally, we finish the chapter with a theorem that can be interpreted as the surjectivety of $\W(-)$, for the case $\g = \lie{sl}(n+1)$. More specifically, for a given irreducible semi-simple coherent family $\M$ we construct a family of simple $\U(\h)$-finite modules such that $\W(M) \sim \M$, for all module $M$ from this family.

\subsection{The case of type C}

In this section assume $\g = \lie{sp}(2n)$. Let $\pi = \{\alpha_1,\dotsc, \alpha_n\}$ be a basis of $\Delta$ indexed such that two consecutive roots are connected, $\alpha_1,\dotsc, \alpha_{n-1}$ are the short roots and $\alpha_n$ is the longest root. Let $h_1,\dotsc,h_n$ be the corresponding coroots.

Recall, from Theorem \ref{theo:class-c.f.-sp2n}, the description of $\X^\infty(C_n)\subseteq \X$ in terms of the weight $\wt(\chi)$, for $\chi \in \X^\infty(C_n)$, and the classification of irreducible semi-simple coherent family.
Notice that the weight $\wt(\chi)$ is dominant, for all $\lambda \in \X^\infty(C_n)$. Moreover, the difference $\wt(\chi_1)-\wt(\chi_2)$ is integral, for any $\chi_1,\chi_2 \in \X^\infty(C_n)$. Therefore, for all $\lambda, \mu \in \wt(\X^{\infty}(C_n))$ the pair $(\lambda, \mu)$ lies in $\Upsilon$, and the corresponding translation functor $\function{T_{\lambda}^{\mu}}{\g\mdfg^{\chi_{\lambda}}}{\g\mdfg^{\chi_{\mu}}}$ is an equivalence of categories. Since $\frkA$ is abelian and preserved by translation functor, then
\[
    \function{T_{\lambda}^{\mu}}{\frkA\cap\g\mdfg^{\chi_{\lambda}}}{\frkA\cap\g\mdfg^{\chi_{\mu}}}
\]
is also an equivalency of category.

Now turn the attention to the simple module $M_0 \in \frkA$ from section~\ref{sec:example-typeC}. Let $\chi_0$ be its central character. Then, $\W(M_0)$ has central character $\chi_0$ and contains an infinite dimensional admissible simple module. Thus, $\chi_0$ must be in $\X^\infty(C_n)$. Let $\lambda_0 = \wt(\chi_0)$ and define the list
\begin{equation}
    M(\lambda) \coloneqq T_{\lambda_0}^{\lambda} M_0 \in \frkA \cap\g\mdfg^{\chi_\lambda}, \quad \lambda \in \wt(\X^\infty(C_n)).
\end{equation}
This list consists of simple modules since $M_0$ is simple (by Prop.~\ref{prop:M_0-simple}) and equivalence of category preserves simplicity.

Recall that, for any $\g$-module $M$ and any automorphism of Lie algebras $\function{\varphi}{\g}{\g}$, we denote by $M^\varphi$ the $\g$-module obtained from $M$ by twisting the action of $\g$ by $\varphi$.
\begin{theo}\label{theo:class-frkA-typeC}
    The set
    \[
        M(\lambda)^\varphi \quad \colon \quad \varphi \in \Aut(\g,\h),\, \lambda \in \wt(\X^\infty(C_n))
    \]
    forms a complete list, up to isomorphism, of simple objects of $\frkA^{\irr}$. Furthermore, for $\varphi \in \Aut(\g,\h)$ and $\lambda \in \wt(\X^\infty(C_n))$, we have the almost-equivalency $\W(M(\lambda)^\varphi) \sim \cEXT(L(\lambda))$
\end{theo}
\begin{proof}
    Let $M \in \frkA$ be simple such that $\W(M)$ is irreducible almost-coherent family. Let $\chi = \chi_\lambda$ be the central character of $M$. By Cor.~\ref{cor:irr-ss-U-c.f.} and Theo.~\ref{theo:class-c.f.-sp2n}, there exists $\lambda' \in \wt(\X^\infty(C_n))$ such that $\W(M) \sim \cEXT(L(\lambda'))$. And since weighting functor and almost-equivalence preserve central characters, we must have $\chi_\lambda = \chi_{\lambda'}$. So let $\lambda = \lambda'$.

    Let $\mathcal{V}$ be the semi-simple coherent family $\cEXT(L(\omega^+))$, where $\omega^+ \in \h^*$ is such that $\omega^+(h_n) = -1/2$ and $\omega^+(h_i) = 0$, for $1 \le i \le n-1$. Clearly $\omega^+ \in\wt(\X^\infty(C_n))$. The coherent family $\mathcal{V}$ is the semi-simplification of the Shale-Weil coherent family defined by O. Mathieu in \cite[Sec.~12]{Mat00} and it is the unique irreducible semi-simple coherent family of \emph{degree one}.
    
    Then the facts that $T_{\lambda}^{\omega^+}$ is equivalence of category and the functor $\W$ commutes with translation functors, together with Lemmas~\ref{lem:u-aprox-translation} and \ref{lem:T_mu^lambda(EXT(L)}, it follows:
    \[
        \W(T_\lambda^{\omega^+} M) \sim T_\lambda^{\omega^+} \cEXT(L(\lambda)) \simeq \cEXT(T_\lambda^{\omega^+}L(\lambda)) \simeq \mathcal{V}.
    \]
    
    In particular, $T_\lambda^{\omega^+} M$ has degree one. And by Theorem~\ref{theo:almost-free}, $T_\lambda^{\omega^+} M$ is actually $\h$-free. Then, by \cite[Theorem~22]{N16}, there exists $\varphi \in \Aut(\g,\h)$ such that 
    \[
        T_\lambda^{\omega^+} M \simeq M_0^\varphi.
    \]
    But, for type C, every automorphism of $\g$ is elementary by \cite[Ch.VIII, \S 5.2, Cor.~ 1 and \S 5.3, Cor.~ 2]{Bor78}. Therefore $(-)^\varphi$ preserves (generalized) central characters by Lemma~\ref{lem:centra-char-with-automorphisms}. In particular, this implies that $\chi_0$, the central character of $M_0$, is actually $\chi_{\omega^+}$. And since $\wt(\chi_{\omega^+}) = \omega^+$, it follows that $\lambda_0 = \omega^+$. Moreover, the preservation of (generalized) central character of $(-)^\varphi$ implies that $V^\varphi \simeq V$ for any finite dimensional $\g$-module $V$ and in particular $(-)^\varphi$ commutes with translation functors. Then we obtain:
    \[
        M \simeq   T_{\omega^+}^\lambda \left(T_\lambda^{\omega^+} M\right) \simeq T_{\omega^+}^\lambda \left(M_0^\varphi\right) \simeq \left(T_{\omega^+}^\lambda M_0\right)^\varphi = M(\lambda)^\varphi,
    \]
    proving the theorem.
\end{proof}

\subsection{The case of type A}

In the following, we fix the Cartan subalgebra $\h \subseteq \lie{sl}(n+1)$ that is the span of $\{h_k = E_{k,k} - E_{k+1,k+1} \mid 1 \le k \le n\}$, where $E_{i,j}$ stands as the $(i,j)$-elementary matrix of $\lie{gl}(n+1)$. The root system is $\Delta= \{ \varepsilon_i - \varepsilon_j \mid 1 \le i,j \le n+1,\, i \ne j\}$, where $\varepsilon_i E_{l,l} = \delta_{i,l}$, and $\{\alpha_i = \varepsilon_i - \varepsilon_{i+1} \mid 1 \le i \le n\}$ forms a basis of $\Delta$ such that $h_k$ is the co-root of $\alpha_k$.
Also, let $\omega_k = \sum_{i=1}^k \varepsilon_i$ be the $k$-th fundamental weight of $\lie{sl}(n+1)$.

\subsubsection{Exponential tensor modules}

Next, we define a subclass of tensor modules that will appear in our classification.

Recall the definition of \emph{$\mathfrak{t}$-finite tensor modules} from Example~\ref{eg:tensor-modules}. For each $A_n$-module $P$ that is finitely generated as $\U(\mathfrak{t})$-module (where $\mathfrak{t} \subseteq A_n$ is subspace spanned by $\{x_k\del_k\mid 1 \le k \le n\}$), finite dimensional $\lie{gl}(n)$-module $V$, and $S\subseteq \{1,\dotsc,n\}=\setn{n}$, the tensor module $T(P,V,S)$ is a $\lie{sl}(n+1)$-module defined in $P\otimes V$. See Section \ref{sec:tensor-modules} for the definitions of the actions.

Since the Weyl algebra $A_n$ is an algebra of differential operators in $\cO = \C[x_1,\dotsc,x_n]$ (where $\del_k$ is identified with $\frac{\del}{\del x_k}$), then $\cO$ is clearly an $A_n$-module. Now, for any $g \in \cO$, define the $A_n$-module $\cO e^g$ by
\[
    x_i \cdot pe^g = x_ipe^g,\quad \del_i \cdot pe^g = \left(\frac{\del}{\del x_i}p+ p\frac{\del}{\del x_i}g\right)e^g,
\]
for all $pe^g \in \cO e^g$ and $1 \le i \le n$. For any $\bb = (b_1,\dotsc,b_n) \in \C^n$, let $\bb\bx = \sum_{i = 1}^n b_ix_i \in \cO$.

\begin{lem}\label{lem:Oe^bx-cyclicy-as-t-mod}
    For every $\bb \in (\C\setminus\{0\})^n$, as a $\U(\mathfrak{t})$-module, the module $\cO e^{\bb\bx}$ is generated by $e^{\bb\bx}$.
\end{lem}
\begin{proof}
    Consider the standard filtration $\{\U_r(\mathfrak{t})\}_{r \ge 0}$ of $\U(\mathfrak{t})$, i.e., the filtration given by $\U_0(\mathfrak{t})_0 = \C$, $\U_1(\mathfrak{t}) = \C + \mathfrak{t}$ and $\U_r(\mathfrak{t}) = \U_1(\mathfrak{t})^r$, for $r \ge 1$. Additionally, consider the natural filtration of $\cO$ defined by $\cO_0 = \C$, $\cO_1 = \Vspan_\C\{1,x_1,\dotsc,x_n\}$ and $\cO_r = (\cO_1)^r$, for all $r \ge 1$.
    To prove the lemma, it is enough to show that $\U_r(\mathfrak{t}) e^{\bb\bx} \supseteq \cO_r e^{\bb\bx}$. We prove it by induction.
    
    Notice that,
    \begin{equation}\label{eq:Oe^b-finit-gen}
        x_k\del_k p e^{\bb\bx} = x_i\left(\frac{\del}{\del x_i}p+ b_ip\right)e^g,
    \end{equation}
    for all $1 \le k \le n$. In particular, $x_k\del_k e^{\bb\bx} = b_ix_ie^{\bb\bx}$. So 
    \[
        \U_i(\mathfrak{t}) e^{\bb\bx} = \Vspan_\C\{e^{\bb\bx}, x_1\del_1 e^{\bb\bx},\dotsc,x_n\del_n e^{\bb\bx}\} = \Vspan_\C\{e^{\bb\bx}, b_1x_1 e^{\bb\bx},\dotsc,b_nx_n e^{\bb\bx}\}.
    \]
    And since $b_i \ne 0$, for all $1 \le i \le n$, we conclude that $\U_1(\mathfrak{t}) e^{\bb\bx} \supseteq \cO_1 e^{\bb\bx}$, proving the step of induction.

    Assume $\U_r(\mathfrak{t}) e^{\bb\bx} \supseteq \cO_r e^{\bb\bx}$, for some $r \ge 1$. Then
    \begin{align*}
        \U_{r+1}(\mathfrak{t})e^{\bb\bx} 
            &= \U_1(\mathfrak{t})\U_r(\mathfrak{t})e^{\bb\bx} \supseteq \U_1(\mathfrak{t})\cO_{r}e^{\bb\bx}\\
            &\overset{\eqref{eq:Oe^b-finit-gen}}{\supseteq} \Vspan_\C\left\{p e^{\bb\bx}, \left(\left(x_i\frac{\del}{\del x_i} + b_ix_i\right)p \right)e^{\bb\bx} \mid 1 \le i \le n, p \in \cO_r\right\}.
    \end{align*}
    But since $x_i\frac{\del}{\del x_i} p \in \cO_r$ and $b_i \ne 0$, for all $p \in \cO_r$ and $1 \le i \le n$, it follows $\U_{r+1}(\mathfrak{t})e^{\bb\bx} \supseteq
    \Vspan_\C\{p e^{\bb\bx}, b_ix_ipe^{\bb\bx} \mid 1 \le i \le n, p \in \cO_r\} = \cO_{r+1} e^{\bb\bx}$. And the lemma is proved.
\end{proof}

\begin{defin}
    For any $\bb \in (\C\setminus\{0\})^n$, finite dimensional $\lie{gl}(n)$-module $V$ and $S \subseteq \setn{n}$, we define
    \[
        E(\bb,V,S) \coloneqq T(\cO e^{\bb\bx}, V,S).
    \]
    We call such modules as \emph{exponential tensor modules}.
\end{defin}

By the previous lemma, exponential tensor modules are $\mathfrak{t}$-finite tensor modules and, therefore, objects of $\frkA$.

\begin{rem}[Ambiguity of exponential tensor modules in \cite{GN22}]
    In \cite{GN22}, the authors defined the exponential tensor modules by $T(\cO^{\psi_S} e^{\bb\bx}, V,\emptyset)$, where $\function{\psi_S}{A_n}{A_n}$ is the algebra automorphism given by
    \[
        \psi_S(x_i) = 
        \begin{cases}
            x_i,& i \notin S;\\
            \del_i,& i \in S;
        \end{cases}
        \quad \text{and} \quad
        \psi_S(\del_i) = 
        \begin{cases}
            \del_i,& i \notin S;\\
            -x_i,& i \in S;
        \end{cases}
    \]
    for all $1\le i \le n$. And $E(\bb,V,S)$ is actually $T((\cO e^{\bb\bx})^{\psi_S}, V,\emptyset)$ which is not isomorphic to $T(\cO^{\psi_S} e^{\bb\bx}, V,\emptyset)$, if $S \ne \emptyset$. But in some results of \cite{GN22}, such as Corollary~5.7 (which will be very important for us later), the authors use our definition of exponential tensor module. Therefore, we must be extra careful when applying propositions from \cite{GN22}.
\end{rem}

Next we are going to investigate the simplicity of a subclass of exponential tensor modules. We will start by looking at modules of the form $E(\bb,V_a,S)$ where $V_a$ denotes the one dimensional $\lie{gl}(n)$-module $\C v_a$ where $x\cdot v_a = a\tr(x)v_a$, for all $x \in\lie{gl}(n)$.

The following lemma is a direct consequence of results of \cite{GN22} and \cite{N15}.
\begin{prop}\label{prop:simplicity-E(b,V_a,S)-non-int-reg}
    Let $\bb \in (\C\setminus\{0\})^n$, $a \in\C$ and $S \subseteq \setn{n}$. 
    \begin{enumerate}[label=(\roman*)]
        \item If $S \ne \emptyset$, then the module $E(\bb,V_a,S)$ is simple.
        \item If $(n+1)(a-1)\notin \Z$ or $(n+1)(a-1) \in \Z_{<0}$, then the module $E(\bb,V_a,\emptyset)$ is simple.
        %\item If $(n+1)(1-a) \in\Z$ and $(n+1)(1-a) < 1$, then $E(\bb,V_a,\emptyset)$ has a unique simple submodule, say $E'(\bb,V_a,\emptyset)$, and the corresponding quotient $E''(\bb,V_a,\emptyset) = E(\bb,V_a,\emptyset)/E'(\bb,V_a,\emptyset)$ has dimension $\binom{(n+1)a - 1}{n}$.
    \end{enumerate}
\end{prop}
\begin{proof}
    Nilsson defined in \cite{N15} a class of $\U(\h)$-finite modules $M_c^I$ (in fact $\U(\h)$-free of rank one), for $c \in \C$ and $I \subseteq \setn{n}$, and gave their simplicity criteria.
    For $\ba \in (\C\setminus\{0\})^{n+1}$, let $\phi_\ba\in \Aut(\lie{sl}(n+1))$ be the automorphism such that $\phi_\ba(E_{i,j}) = \frac{a_i}{a_j}E_{i,j}$. By \cite[Cor.~5.7]{GN22}, we have that
    \[
        E(\bb,V_a,S) \simeq \left(M_{a-1}^{\setn{n}\setminus S}\right)^{\phi_{\tilde\bb}},
    \]
    where $\tilde\bb = (\tilde b_1,\dotsc,\tilde b_n, 1) \in (\C\setminus\{0\})^{n+1}$ satisfies $\tilde b_i = -\frac{1}{b_i}$, if $i \in S$, and $\tilde b_i = b_i$, otherwise.
    Then the theorem follows from \cite[Theorems.~3.2 and 3.4]{N15}.
\end{proof}

Next, we identify $\lie{gl}(n)$ with the subalgebra of $\lie{sl}(n+1)$ generated by $\{E_{i,j},\tilde h_k, \mid 1 \le i,j,k \le n,\, i\ne j\}$, so that $\h$ is a Cartan subalgebra of $\lie{gl}(n)$.
The fixed basis of $\Delta$ determines the Borel subalgebra $\mathfrak{b}\subseteq \lie{sl}(n+1)$ given by the span of $E_{i,j}$ with $1 \le i < j \le n$. Then $\mathfrak{b}' = \mathfrak{b}\cap\lie{gl}(n)$ defines a Borel subalgebra so that the simple weight $\lie{gl}(n)$-module $L_{\lie{gl}(n)}(\lambda)$ of highest weight $\lambda$ is well-defined.

Let $P_{\lie{gl}(n)}^+$ be the set of all $\lambda \in \h^*$ such that $L_{\lie{gl}(n)}(\lambda)$ is finite dimensional. Then 
\[
    P_{\lie{gl}(n)}^+ = \{a_1\omega_1+\dotsb + a_n\omega_n \in \h^* \mid a_i \in \Z_{\ge 0},\, \text{for }1 \le i \le n-1\}.
\]
Recall the following results from \cite{GN22}:
\begin{lem}[{\citealp[Cor.~4.2]{GN22}}]\label{lem:central-char-T(P,V)}
    For every $A_n$-module $P$ and $\lambda \in P_{\lie{gl}(n)}^+$, the modules $T(P,L_{\lie{gl}(n)}(\lambda))$ and $L(\lambda - (n+1)\omega_n)$ have the same central character.
\end{lem}
\begin{lem}[{\citealp[Prop.~4.5]{GN22}}]\label{lem:trans-tensor-mod}
    Let $\lambda,\mu \in P_{\lie{gl}(n)}^+$ such that $(\lambda - (n+1)\omega_n, \mu - (n+1)\omega_n) \in \Upsilon$. For every $A_n$-module $P$ we have the isomorphism
    \[
        T_{\lambda - (n+1)\omega_n}^{\mu - (n+1)\omega_n}T(P,L_{\lie{gl}(n)}(\lambda)) \simeq T(P,L_{\lie{gl}(n)}(\mu)).
    \]
\end{lem}

Fix the notation
\[
    E(\bb,\lambda,S) = E(\bb,L_{\lie{gl}(n)}(\lambda + (n+1)\omega_n), S),\quad \bb \in (\C\setminus\{0\})^n,\,\lambda \in P_{\lie{gl}(n)}^+,\,S\subseteq\setn{n}.
\]
Notice that, for $a \in \C$, the $\lie{gl}(n)$-module $V_a$ is isomorphic to $L_{\lie{gl}(n)}(a(n+1)\omega_n)$. Then, Lemma~\ref{lem:trans-tensor-mod} implies that
\[
    P_{\lie{gl}(n)}^{+,1} = \{ \lambda \in P_{\lie{gl}(n)}^+ \mid (\lambda,a(n+1)\omega_n) \in \Upsilon\, \text{for some }a \in \C\}
\]
is precisely the subset of $P_{\lie{gl}(n)}^{+}$ consisting of weights $\lambda$ such that $E(\bb,\lambda,S)$ can be translated to a module isomorphic to $E(\bb,V_{a_{\lambda}},S)$, for some ${a_\lambda} \in \C$. See \cite[Prop.~4.5]{GN22} for the existence of the pairs $(\lambda , a_\lambda)$. Moreover, the pair $(\lambda , a_\lambda)$ satisfies
\begin{itemize}
    \item If $\lambda$ is non-integral, then $(n+1)(a_\lambda-1) \notin \Z$; and
    \item If $\lambda$ is integral and singular, then $(n+1)(a_\lambda-1) \in \Z_{<0}$.
\end{itemize}
Then, these properties and Proposition~\ref{prop:simplicity-E(b,V_a,S)-non-int-reg} imply:

\begin{prop}\label{prop:simplicity-E(b,lambda,S)-non-integ}
    For every $\bb \in (\C\setminus\{0\})^n$, $S \subseteq \setn{n}$ and $\lambda \in P_{\lie{gl}(n)}^{+,1}$, the module $E(\bb,\lambda,S)$ is simple whenever $S \ne \emptyset$ or $\lambda$ is either non-integral or integral and singular.
\end{prop}

Let $\function{\tau}{\lie{sl}(n+1)}{\lie{sl}(n+1)}$ be the negative transpose of $\lie{gl}(n+1)$ restricted to $\lie{sl}(n+1)$. Since $\tau$ preserves $\h$, the functor $-^\tau$ is also an auto-equivalence of the categories $\frkA$.

\begin{prop}[{\cite{GN22}, Theo~5.6, Cor.~5.7} and {\cite{N15}, Theo.~30}]\label{prop:h-free-mod-rank-one}
    For every $\bb \in (\C\setminus \{0\})^n$, finite dimensional $\lie{gl}(n)$-module $V$ and $S \subseteq \setn{n}$, the $\lie{sl}(n+1)$-module $E(\bb,V,S)$ is $\h$-free of rank equal to $\dim V$. Furthermore, 
    \[
        E(\bb,a(n+1)\omega_n,S),\,E(\bb,a(n+1)\omega_n, S)^\tau \quad \colon \quad \bb \in (\C\setminus \{0\})^n, a \in \C, S\subseteq\setn{n},
    \]
     forms a complete list (up to isomorphism) of mutually non-isomorphic $\lie{sl}(n+1)$-module that are $\h$-free of rank one.
\end{prop}

The above proposition gives another interpretation of $P_{\lie{gl}(n)}^{+,1}$, that justify the upper index $1$ in the notation: $P_{\lie{gl}(n)}^{+,1}$ is precisely the subset of $P_{\lie{gl}(n)}^{+}$ consisting of weights $\lambda$ such that $E(\bb,\lambda,S)$ can be translated to a module that is $\U(\h)$-free of rank $1$.

\subsubsection{The (almost complete) classification}
Throughout this section, we assume that $n = \rank \g > 1$ and we classify simple objects of $\frkA^{\irr}$ with \textbf{non-integral} or \textbf{integral and singular} central character.

First, recall from Proposition~\ref{prop:central-char-adm-sln+1}, the description of $\X^\infty(A_n)\subseteq \X$ in terms of the weight  $\wt(\chi)$ (and $i_\chi$, in the integral singular case), for $\chi \in \X^\infty(A_n)$. Also recall, from equation~\eqref{eq:def-w_k}, the definition of the elements $w_0,w_1, \dotsc,w_n$ of the Weyl group $W$.

Notice that, for $\chi \in \X^\infty(A_n)$, we have
\begin{align*}
    w_n \cdot \wt(\chi) \in P_{\lie{gl}(n)}^+,\qquad
        &\text{if $\chi$ is non-integral,}\\
    w_{n-i_\chi} \cdot \wt(\chi) \in P_{\lie{gl}(n)}^+,\qquad
        &\text{if $\chi$ is integral singular.}
\end{align*}
Let $\wlongest \in W$ be the longest element. One has that the action of $-\wlongest$ in $\h^*$ sends $\omega_i$ to $\omega_{n-i+1}$, for all $1 \le i \le n$. Then, considering $c_i = s_1s_2\dotsb s_{i-1}$ (setting $c_1$ as the identity), a straightforward computation shows that
\begin{align*}
    -\wlongest\wt(\chi) \in P_{\lie{gl}(n)}^+,\qquad
        &\text{if $\chi$ is non-integral,}\\
    -\wlongest(c_{i_\chi} \cdot \wt(\chi)) \in P_{\lie{gl}(n)}^+,\qquad
        &\text{if $\chi$ is integral singular.}
\end{align*}

As a consequence of the preservation of central character by $\W$, together with Lemma~\ref{lem:almost-equiv-central-char}, we have that the set of central character of simple objects of $\frkA^{\irr}$ is contained in $\X^\infty(A_n)$. 
Furthermore, recall that every irreducible semi-simple coherent family of non-integral or integral and singular central character $\chi$ is isomorphic to $\cEXT(L(\wt(\chi)))$, by Theorem~\ref{theo:class-c.f.-sl(n+1)}. 
Then may we state:

\begin{theo}\label{theo:class-simple-multipli-one-non-integral-reg}
    Let $\chi \in \X^\infty(A_n)$ be either non-integral or integral and singular. Let $\lambda = \wt(\chi)$ (and $i_\chi$ in the integral singular case). For each of the next cases, the following form a \textbf{complete} list (up to isomorphism) of mutually non-isomorphic simple objects of $\frkA^{\irr}$ and central character $\chi$:
    \begin{itemize}
        \item if $\chi$ is non-integral:
        \[
            E(\bb, w_n\cdot \lambda,S),\, E(\bb,-\wlongest\lambda,S)^\tau, \qquad \bb \in (\C\setminus\{0\})^n, S\subseteq\setn{n};
        \]
        \item if $\chi$ is singular and integral:
        \[
            E(\bb, w_{n-i}\cdot \lambda,S),\,E(\bb, - \wlongest(c_i\cdot \lambda),S)^\tau, \qquad \bb \in (\C\setminus\{0\})^n, S\subseteq\setn{n}.
        \]
    \end{itemize}
    Moreover, we have that $\W(M)$ is almost-equivalent to $\cEXT(L(\lambda))$ for every module $M$ in the list.
\end{theo}

Before proving the theorem, we recall the following result given by Britten and Lemire.
\begin{prop}[{\citealp[Prop.~1.6]{BL01}}]\label{prop:c-f-deg-1}
    A coherent family has degree one if and only if it has a semi-simplification isomorphic to 
    \[
        \cEXT(L(a\omega_1)) \quad \text{or} \quad \cEXT(L(-(N+2)\omega_1 + (N+1)\omega_2)),
    \]
    for some $a \in \C\setminus\Z_{\ge 0}$ or $N \in \Z_{\ge 0}$.
\end{prop}

\begin{proof}[Proof of Theorem~\ref{theo:class-simple-multipli-one-non-integral-reg}]
    Let an infinite dimensional simple module $M \in \frkA^{\irr}$ of central character $\chi$ (and let $\lambda$ and $i_\chi$ as the theorem's statement). Then $\W(M) \sim \cEXT(L(\lambda))$ by Theorem~\ref{theo:class-c.f.-sl(n+1)} and Corollary~\ref{cor:irr-ss-U-c.f.}.
    
    To unify the argument for both cases, fix the notations
    \begin{itemize}
        \item{if $\lambda$ is non-integral:
        \[
            \tilde \lambda \coloneqq w_n\cdot \lambda,\quad
            \hat \lambda \coloneqq \lambda,\quad
            \mu \coloneqq \lambda(h_1)\omega_1,\quad
            \tilde \mu \coloneqq w_n \cdot \mu,\quad
            \text{and}\quad\hat{\mu} \coloneqq \mu.
        \]}
        \item{if $\lambda$ is singular:
        \[
            \tilde \lambda \coloneqq w_{n-i_\chi}\cdot \lambda,\quad
            \hat \lambda \coloneqq c_{i_\chi}\cdot \lambda,\quad
            \mu \coloneqq -\omega_{i_\chi},\quad
            \tilde \mu \coloneqq w_{n-i_{\chi}} \cdot \mu,\quad
            \text{and}\quad\hat{\mu} \coloneqq c_{i_\chi} \cdot \mu = -i_{\chi}\omega_1.
        \]
        }
    \end{itemize}
    
    Notice that $(\lambda, \mu) \in \Upsilon$. Moreover, by Theorem~\ref{theo:class-c.f.-sl(n+1)} and Proposition~\ref{prop:c-f-deg-1}, it follows that $\cEXT(L(\mu)) \simeq \cEXT(L(\hat{\mu}))$ which has degree one. We have:
    \begin{align*}
        \W(T_\lambda^\mu M)
            &\simeq T_\lambda^\mu\W(M)\quad
                &\text{(Lemma~\ref{lem:prop-Wk}\ref{lem:prop-Wk-projective-functor})}\\
            &\sim T_\lambda^\mu \cEXT(L(\lambda))\quad
                &\text{(Lemma~\ref{lem:u-aprox-translation} and $\W(M) \sim \cEXT(L(\lambda))$)}\\
            &\sim \cEXT(L(\mu)) \simeq \cEXT(L(\hat{\mu}))\quad
                &\text{(Lemma~\ref{lem:T_mu^lambda(EXT(L)}).}
    \end{align*}
    Since, for any simple module $N \in \frkA$ of infinite dimension, $\W(N)$ is an almost-coherent family of degree equal to $\rank N$ (by Cor.~\ref{Cor:W(M)-U-coherent-family}), we conclude that $T_\lambda^\mu M\in \frkA$ is simple of rank one. Furthermore, Theorem~\ref{theo:almost-free} and the fact that $\mu$ is non-integral or singular imply that $T_\lambda^\mu M$ is $\U(\h)$-free. So, by Prop.~\ref{prop:h-free-mod-rank-one} it follows that
    \[
        i)\, T_\lambda^\mu M \simeq E(\bb,a(n+1)\omega_n,S); \quad \text{or} \quad  ii)\, T_\lambda^\mu M \simeq E(\bb,a(n+1)\omega_n,S)^{\tau},
    \]
    for some $\bb \in (\C\setminus\{0\})^n$, $a \in \C$ and $S \subseteq \setn{n}$.

    Consider the case $i)$. By Harish-Chandra theorem, $\chi_\lambda = \chi_{\tilde \lambda}$ and $\chi_\mu = \chi_{\tilde \mu}$. Moreover, $(\tilde\lambda, \tilde\mu) \in \Upsilon$ and $\tilde \mu \in P_{\lie{gl}(n)}^+$. Then $E(\bb,a(n+1)\omega_n,S) \simeq T_\lambda^\mu M  = T_{\tilde\lambda}^{\tilde\mu} M$ and Lemma~\ref{lem:central-char-T(P,V)} implies that $a(n+1)\omega_n = \tilde\mu + (n+1)\omega_n$. Therefore, 
    \[
        M \simeq T^{\tilde\lambda}_{\tilde\mu}E(\bb,\tilde\mu,S) \simeq E(\bb,\tilde\lambda,S),
    \]
    where the last isomorphism follows from Lemma~\ref{lem:trans-tensor-mod} together with the fact that $\tilde\lambda \in P_{\lie{gl}(n)}^+$.

    Suppose now we are in the case $ii)$. Then Lemma~\ref{lem:centra-char-with-automorphisms} and the fact that $(-\wlongest\hat{\lambda},-\wlongest\hat{\mu})$ lies in $\Upsilon$ imply that $(T_\lambda^\mu - )^\tau \simeq T_{-\wlongest\hat{\lambda}}^{-\wlongest\hat{\mu}}(-^\tau)$. Then, since $-\wlongest\hat{\mu} \in P_{\lie{gl}(n)}^+$ and the fact that $(-^\tau)^\tau$ is the identity functor, we have that $E(\bb,a(n+1)\omega_n,S) \simeq (T_\lambda^\mu M)^\tau \simeq T_{-\wlongest\hat{\lambda}}^{-\wlongest\hat{\mu}}(M^\tau)$. Then Lemma~\ref{lem:central-char-T(P,V)} implies that $a(n+1)\omega_n = -\wlongest\hat{\mu} + (n+1)\omega_n$.
    Therefore,
    \[
        M \simeq (T^{-\wlongest\hat{\lambda}}_{-\wlongest\hat{\mu}}E(\bb,-\wlongest\hat{\mu},S))^\tau \simeq E(\bb,-\wlongest\hat{\lambda},S)^\tau,
    \]
    where the last isomorphism follows from Lemma~\ref{lem:trans-tensor-mod} together with the fact that $-\wlongest\hat{\lambda} \in P_{\lie{gl}(n)}^+$.

    The non-redundancy of the given lists in the Theorem's statement follows by Prop.~\ref{prop:h-free-mod-rank-one}.
\end{proof}

\begin{rem}
Although the almost equivalence of $\W(M)$ for the modules listed in Theorem~\ref{theo:class-simple-multipli-one-non-integral-reg} follows from Proposition 5.10 and Remark 5.5 of \cite{GN22}, the \textbf{completeness} of the stated list does not follow directly from the results in \cite{GN22}.
\end{rem}

\begin{rem}
    The argument used for proving Theorem~\ref{theo:class-simple-multipli-one-non-integral-reg} has it core in the possibility to translate a coherent family to another of degree one. Such property is true whenever the central character is non-integral or integral and singular. When the central character of a coherent family is integral regular, we generally can not make such translation. This is the reason our classification of simple objects of $\frkA^\irr$ does not consider all possible central characters.
\end{rem}

\subsection{Surjectivity of $\W$ for type A via parabolic induced modules}
Let again $\g$ be $\lie{sl}(n+1)$ and we consider the same definitions and notations of elements and subalgebras of $\g$ established in the previous section. Moreover, by fixing the canonical basis $\{e_1,\dotsc,e_n,e_{n+1}\}$ of $\C^{n+1}$ and assigning $E_{i,j}$ to the linear map $e_k \mapsto \delta_{k,j}e_i$, we also identify $\lie{gl}(n+1)$ (and hence $\g$) with $\lie{gl}(\C^{n+1})$ (respectively, $\lie{sl}(\C^{n+1})$). So that we also use the notations and definitions introduced in Section~\ref{sec:parabolic-induced}. Recall that $\langle ,\rangle$ is the complex inner product of $\C^{n+1}$ such that $\{e_1,\dotsc,e_{n+1}\}$ is orthonormal. 

For every $\bb \in (\C\setminus\{0\})^n$, define $v_\bb$ to be the normalization of $\sum_{i = 1}^{n}b_ie_i + e_{n+1} \in \C^{n+1}$ and fix an orthonormal basis $B_{\bb} = \{\tilde e_1, \dotsc, \tilde e_n,v_\bb\}$ of $\C^{n+1}$. Let $S'_\bb\in GL(\C^{n+1})$ be the linear operator that defined by $e_i \mapsto \tilde e_j$, for $1 \le j \le n$ and $e_{n+1} \mapsto v_\bb$, and let $S_\bb$ be the unitary operator $\frac{1}{\det S'_\bb}S'_\bb \in SU(\C^{n+1})$. Then, by \cite[Chap VII,\S 3, Remmark 1.2]{Bor78}, the automorphism
\[
    \function{\varphi_\bb}{\g}{\g}\quad x\mapsto\left(\map{\varphi_\bb(x)}{v}{S_\bb x S_\bb^{-1}(v)}\right)
\]
is elementary, i.e, $\varphi_\bb \in\Aut_e(\g)$. One can check that:
\begin{lem}\label{lem:propers-varphi_bb}
    For every $\bb\in(\C\setminus\{0\})^n$, the automorphism $\varphi_\bb$ commutes with $\tau$ and send $\p_{\C e_n}$ to $\p_\bb$:
\end{lem}

For every $\bb \in (\C\setminus\{0\})^n$, define $\h_\bb = \varphi_\bb(\h)$. Then the pullback of $\varphi_\bb$ defines an isomorphism $\function{\varphi_\bb^*}{\h}{\h_\bb^*}$. Then the Weyl group $W$ acts on $\h_\bb$ by pulling back the action on $\h$.
Consider the parabolic subalgebra $\p_0 = \tau(\p_{\C e_{n+1}})$ and its nilradical $\mathfrak{u}^+$. Notice that $\p_0 = \p$ from Section~\ref{sec:class-irr-ss-cf-A} and that the defined Levi factor of $\p$ is $\lie{l} = \p_0 \cap \tau(\p_0)$. 
Write $\mathfrak{u}_\bb^+ = \varphi_\bb(\mathfrak{u}^+)$ and $\mathfrak{l}_\bb = \varphi_\bb(\mathfrak{l})$. 
For $w_0,\dotsc,w_n \in W$ defined in equation~\eqref{eq:def-deg_w_k}, it is straightforward to check that $L_{\mathfrak{l}_\bb}(w_k \cdot \varphi_\bb^*(\lambda))$, the simple $\mathfrak{l}_\bb$-module of highest $\h_\bb$-weight $w_k \cdot \varphi_\bb^*(\lambda)$, is finite dimensional.
For every $\lambda \in P^+$ and $1 \le k \le n$, we define
\[
    L_\bb(k,\lambda) \coloneqq L_{\p_\bb}(L_{\mathfrak{l}_\bb}(w_k \cdot \varphi_\bb^*(\lambda))) \in \cO^{\p_\bb}
    \quad
    \text{and}
    \quad
    L_\bb^\tau(k,\lambda) \coloneqq L_\bb(n-k+1,-\wlongest\lambda)^\tau \in \cO^{\tau(\p_\bb)},
\]
where $L_{\mathfrak{l}_\bb}(w_k \cdot \varphi_\bb^*(\lambda))$ is viewed as $\p_\bb$-module by defining $\mathfrak{u}_\bb^+$ to act trivially on it. These modules are simple by definition and lie in $\frkA$ by Corollary~\ref{cor:parabolic-bgg-O^p}.

\begin{lem}\label{lem:iso-and-central-char-parabolic-induced}
    For $1 \le k \le n$ and $\lambda \in P^+$, we have the isomorphisms of $\g$-modules:
    \[
        L_\bb(k,\lambda) \simeq L_{\p}(L_{\mathfrak{l}}(w_k\cdot \lambda))^{\varphi_\bb^{-1}}
        \quad
        \text{and}
        \quad
        L_\bb^\tau(k,\lambda) \simeq L_{\p}(L_{\mathfrak{l}}(w_{n-k+1}\cdot (-\wlongest\lambda))^{\varphi_\bb^{-1}\tau}.
    \]
    Moreover, both modules have central character $\chi_\lambda$.
\end{lem}
\begin{proof}
    The first claim follows from the fact that
    \begin{align*}
        M_{\p}(V)^\varphi = \U(\g)\otimes_{\p}V&\longrightarrow \U(\g)\otimes_{\varphi^{-1}(\p)}V^\varphi = M_{\varphi^{-1}(\p)}(V^\varphi)\\
        u\otimes v &\longmapsto \varphi^{-1}(u)\otimes v,
    \end{align*}
    defines an isomorphism of modules, for all $\p$-modules $V$ and automorphism $\varphi \in \Aut(\g)$.

    Next, we show that there exists a Cartan subalgebra $\tilde\h_\bb \subseteq \g$ stabilized by $\varphi_\bb$. Since $S_\bb$ is unitary and invertible, there exist a unitary invertible operator $\tilde S \in \Aut(\C^{n+1})$ and an operator $\map{D}{e_i}{a_ie_i}$ (where $a_1,\dotsc,a_{n+1} \in \C$ have norm one) such that $S_\bb = \tilde S D {\tilde S}^{-1}$. Then, one can show that the following is a Borel subalgebra:
    \[
        \tilde{\mathfrak{b}}_\bb = \{f \in \g \mid f(\tilde{S}(V_j)) \subseteq \tilde{S}(V_j),\, 1 \le j \le n+1\},
    \]
    where $V_j = \Vspan_\C\{e_1,\dotsc,e_j\}$. Moreover, a straightforward computation shows that $\tilde \h_\bb = \tilde{\mathfrak{b}}_\bb \cap \tau(\tilde{\mathfrak{b}}_\bb)$ forms a Cartan subalgebra and that $\varphi_\bb(\tilde{\mathfrak{b}}_\bb) = \tilde{\mathfrak{b}}_\bb$. Then the commutation of $\varphi_\bb$ and $\tau$ (by Lemma~\ref{lem:propers-varphi_bb}) implies the equality $\varphi_\bb(\tilde \h_\bb) = \tilde \h_\bb$.

    Then $\varphi_\bb \in \Aut_e(\g,\tilde \h_\bb)$ and Lemma~\ref{lem:centra-char-with-automorphisms} implies the statement about central character.
\end{proof}

By Theorem~\ref{theo:class-simple-multipli-one-non-integral-reg}, there exists a module $M \in \frkA$ such that $\W(M)$ is almost equivalent to a given irreducible coherent family $\M$, whenever their central character are not integral and regular. Thus, the following statement, together with the classification of irreducible semi-simple coherent family in Theorem~\ref{theo:class-c.f.-sl(n+1)}, concludes the ``subjectivity'' of $\W$:

\begin{theo}\label{theo:surjectivety-of-W}
    Let $\lambda \in P^+$ and $1 \le k \le n$. Then, every module $M$ in the list
    \[
        \mathcal{P}_{k,\lambda} = \{L_\bb(k,\lambda),L_\bb^\tau(k,\lambda) \mid \bb \in \C^n\}
    \]
    satisfies $\W(M)\sim \cEXT(L(w_k\cdot\lambda))$.
\end{theo}

Before proving, we need to consider the degree of coherent families involved. Recall Theorem~\ref{theo:deg-c.f.} and the description of the set $P^+_{\mathfrak{l}}$. By Weyl's dimension formula, we have that $\dim L_{\mathfrak{l}}(\nu)$ is polynomial in $\nu \in P^+_{\mathfrak{l}}$. Since $P^+_{\mathfrak{l}}$ is a dense subset of $\h^*$, the map $P^+_{\mathfrak{l}} \ni \nu\mapsto \dim L_{\mathfrak{l}}(\nu) \in \C$ can be extended to a polynomial function in $\h^*$. Therefore, there exists a unique polynomial $\function{\deg_{k}}{\h^*}{\C}$ such that
\begin{equation}\label{eq:def-deg_w_k}
    \deg_{k}(\lambda) = \deg \cEXT(L(w_k\cdot \lambda)) = \deg L(w_k\cdot \lambda), \quad \lambda \in P^+.
\end{equation}

\begin{lem}\label{lem:deg_w_i's-linear-independent}
    The polynomials $\deg_{1},\dotsc,\deg_{n}$ are linear independent over $\C$.
\end{lem}
\begin{proof}
    For $1 \le k \le n$, define $\function{d_k}{\h^*}{\C}$ to be the polynomial such that $d_k(\lambda) \coloneqq \dim L_{\mathfrak{l}}(w_k \cdot \lambda)$, for all $\lambda \in P^+$. We have $\deg_{w_k}(\lambda) = \sum_{i=0}^{n-k}(-1)^{i} d_{k+i}(\lambda)$, for all $\lambda \in \h^*$.
    The lemma follows by open the formulas and finding elements $\lambda^{(1)},\dotsc, \lambda^{(n)} \in \h^*$ such that
    \[
        d_k(\lambda^{(k)}) \ne 0 \qquad \text{and}\qquad d_{\ell}(\lambda^{(k)}) = 0, \quad \forall k < \ell \le n.
    \]
    For more details in the computation, see the author's thesis \cite{Men}.
\end{proof}

We also need the fallowing result about the behavior of the functor $(-)^\tau$ with irreducible coherent families.
\begin{lem}\label{lem:EXT(L(lambda))^tau}
    Let $\lambda \in P^+$ and $1 \le k \le n$. Then we have the isomorphism of modules
    \[
        \cEXT(L(w_k\cdot \lambda))^\tau \simeq \cEXT(L(w_{n-k+1}\cdot(-\wlongest\lambda))).
    \]
\end{lem}
\begin{proof}
    For every simple admissible weight module $L$ of infinite dimension, one can check that $\cEXT(L)^\tau$ is as semi-simple coherent family that contains $L^\tau$, and thus $\cEXT(L)^\tau \simeq \cEXT(L^\tau)$ by Prop.~\ref{prop:EXT(L)}. Thus $\cEXT(L(w_k\cdot \lambda))^\tau \simeq \cEXT(L(w_k\cdot \lambda)^\tau)$.

    Any basis $\pi'$ for the set of roots $\Delta$ define the set $\Delta^+_{\pi'} = \Z_{\ge 0}\pi' \cap \Delta$ of positive roots in relation to $\pi'$ and thus a Borel subalgebra $\frak{b}_{\pi'} = \h\bigoplus_{\beta \in \Delta^+_{\pi'}} \g_{\beta}$.
    Then, using the notation from \cite{Mat00}, we let $L_{\pi'}(\lambda)$ be the simple weight module of $\pi'$-highest weight $\lambda\in\h^*$, i.e., the simple quotient of $\U(\g)\otimes_{\U(\frak{b}_{\pi'})}V_\lambda$, where $V_\lambda$ is the unidimensional $\frak{b}_{\pi'}$-module with $\h$-weight $\lambda$. Then, since $L(w_k\cdot \lambda)^\tau$ is a simple weight module of lowest weight $-w_k\cdot \lambda$ we obtain that it is isomorphic to $L_{\wlongest\pi}(-w_k\cdot\lambda))$, where $\pi = \{\alpha_1,\dotsc,\alpha_n\}$. By \cite[Lemma~6.1]{Mat00}, we have that $\cEXT(L_\pi(\lambda)) \simeq \cEXT(L_{w\pi}(w\lambda))$, for every $w \in W$. In particular, we have
    \[
        \cEXT(L(w_k\cdot \lambda)^\tau) \simeq \cEXT(L_{\wlongest\pi}(-w_k\cdot\lambda)) \simeq \cEXT(L(\wlongest(-w_k\cdot\lambda)))
    \]

    Finally, check that $\wlongest(-w_k\cdot\lambda) = s_1s_2\dotsc s_k\cdot (-\wlongest\lambda)$. Then by Theorem~\ref{theo:class-c.f.-sl(n+1)} we have
    \[
        \cEXT(L(w_k\cdot \lambda))^\tau \simeq \cEXT(L(s_1s_2\dotsc s_k\cdot (-\wlongest\lambda))) \simeq \cEXT(L(w_{n-k+1}\cdot (-\wlongest\lambda))).
    \]
\end{proof}

\begin{proof}[Proof of Theorem~\ref{theo:surjectivety-of-W}]
    
    Since $(-)^\varphi$ is an exact functor, for all $\varphi \in\Aut(\g)$, it follows from Lemma~\ref{lem:iso-and-central-char-parabolic-induced} and Theorem~\ref{theo:deg-c.f.} that we have exact sequences
    \[
        0 \rightarrow L_\bb(k+1,\lambda) \rightarrow M_{\p_\bb}(L_{\mathfrak{l}_\bb}(w_k \cdot \varphi_\bb^*(\lambda))) \rightarrow L_\bb(k,\lambda) \rightarrow 0
    \]
    for all $1 \le k \le n$ and $\lambda \in P^+$, where $L_\bb(n+1,\lambda) = 0$. Write $\rank_k(\lambda) = \rank L_\bb(k,\lambda)$, where $\rank_{n+1} = 0$. Then
    \[
        \rank_{k+1}(\lambda)+ \rank_k(\lambda) = \rank M_{\p_\bb}(L_{\mathfrak{l}_\bb}(w_k \cdot \varphi_\bb^*(\lambda))).
    \]
    But $M_{\p_\bb}(L_{\mathfrak{l}_\bb}(w_k \cdot \varphi_\bb^*(\lambda)))$ is $\U(\h)$-free of rank equal to $\dim L_{\mathfrak{l}_\bb}(w_k \cdot \varphi_\bb^*(\lambda))$, by Example~\ref{eg:parabolic-induced}. Moreover, $L_{\mathfrak{l}_\bb}(w_k \cdot \varphi_\bb^*(\lambda)) \simeq  L_{\mathfrak{l}}(w_k \cdot \lambda)^{\varphi_\bb^{-1}}$ and thus its dimension is equal to $\dim L_{\mathfrak{l}}(w_k \cdot \lambda)$. Thus, 
    \[
        \rank_{k+1}(\lambda)+ \rank_k(\lambda) = \dim L_{\mathfrak{l}}(w_k \cdot \lambda)
    \]
    for all $1 \le k \le n$ and $\lambda \in P^+$. Then, with the formula of the degree of coherent families given in Theorem~\ref{theo:deg-c.f.}, we may conclude that
    \begin{equation}\label{eq:rank-parabolic-induced}
        \rank_k(\lambda) = \deg_k(\lambda), \quad \lambda \in P^+.
    \end{equation}

    Since $L_\bb(k,\lambda)$ is simple of central character $\chi_\lambda$, the Theorems \ref{theo:W(M)-quasi-coh-fam} and \ref{theo:U-aprox-ss-U-c.f.} and the classification of irreducible semi-simple coherent family of integral regular central character implies that
    \[
        \W(L_\bb(k,\lambda)) \sim \bigoplus_{i=1}^n \cEXT(L(w_i\cdot \lambda))^{m_{i,k}(\lambda)}
        \quad \text{and}\quad
        \rank L_\bb(k,\lambda) = \sum_{i = 1}^n m_{i,k}(\lambda)\deg \cEXT(L(w_i\cdot \lambda))
    \]
    for $m_{1,k}(\lambda),\dotsc,m_{n,k}(\lambda) \in \N$. Moreover, the commutation of the translation functors with $\W$ and $\cEXT$ (Lemmas~\ref{lem:prop-Wk}\ref{lem:prop-Wk-projective-functor} and \ref{lem:T_mu^lambda(EXT(L)}) implies that $m_{i,k}(\lambda)= m_{i,k}$ does not depend on $\lambda$. Thus:
    \[
        \rank_k(\lambda) = \sum_{i = 1}^n m_{i,k}\deg_i(\lambda).
    \]
    Finally, equation \eqref{eq:rank-parabolic-induced} and the linear independence of $\deg_1,\dotsc,\deg_n$ (Lemma~\ref{lem:deg_w_i's-linear-independent}) implies that $m_{i,k} = \delta_{i,k}$. Therefore, it follows that $\W(L_\bb(k,\lambda)) \sim \cEXT(L(w_k\cdot \lambda))$.

    It is straightforward to show that $(-)^\tau$ preserves almost-equivalence and commutes with $\W$. Then, for every $1 \le k \le n$ and $\lambda \in P^+$:
    \begin{align*}
        \W(L_\bb^\tau(k,\lambda)) 
            &=\W(L_\bb(n-k+1,-\wlongest\lambda)^\tau)\\
            &\simeq\W(L_\bb(n-k+1,-\wlongest\lambda))^\tau\\
            &\sim \cEXT(L(w_{n-k+1}\cdot (-\wlongest\lambda)))^\tau\\
            &\simeq \cEXT(L(w_k\cdot \lambda))),
    \end{align*}
    where the last isomorphism follows from Lemma~\ref{lem:EXT(L(lambda))^tau}, which concludes the proof
\end{proof}

\begin{rem}[Formula for $\deg_1$]\label{rem:formula-deg_1}
    A direct consequence of the above proof is that $\deg_1(\lambda) = \rank_1(\lambda)$, for every $\lambda \in P^+$. By Theorem~\ref{theo:deg-c.f.}, there exists an exact sequence
    \[
        0 \rightarrow L_\bb(1,\lambda) \rightarrow M_{\p_\bb}(L_{\mathfrak{l}_\bb}(\varphi_\bb^*(\lambda))) \rightarrow L(\lambda) \rightarrow 0,
    \]
    and we have that $L(\lambda)$ is finite dimensional. Thus $\rank_1(\lambda)$ is equal to the rank of $M_{\p_\bb}(L_{\mathfrak{l}_\bb}(\varphi_\bb^*(\lambda)))$, which is equal to $\dim L_{\mathfrak{l}_\bb}(\varphi_\bb^*(\lambda))$. Therefore $\deg_1(\lambda) = \dim L_{\mathfrak{l}}(\lambda)$.
\end{rem}

%\newpage

%%%%%%%%%%%%%%%%%%%%%%%%%%%%%%%%%%%%%%%%%%%%%%%%%%%%%%%%%%%%%%%%%%%
% Section - Conjectures
%%%%%%%%%%%%%%%%%%%%%%%%%%%%%%%%%%%%%%%%%%%%%%%%%%%%%%%%%%%%%%%%%%%
%%%%%%%%%%%%%%%%%%%%%%%%%%%%%%%%%%%%%%%%%%%%%%%%%%%%%%%%%%%%%%%%%%%
% Section - Conjectures
%%%%%%%%%%%%%%%%%%%%%%%%%%%%%%%%%%%%%%%%%%%%%%%%%%%%%%%%%%%%%%%%%%%
\section{Conjectures for the classification of $\frkA^{\irr}$ for type A}\label{sec:conjectures}

Let $\g = \lie{sl}(n+1)$ and consider the same definitions and notations of its elements and subalgebras that were established in the previous section. 

For any abelian category $\mathcal{C}$ with object of finite length, let $K_0(\mathcal{C})$ be the Grothendieck group of $\mathcal{C}$, i.e, the free abelian group on the isomorphism classes of objects of $\mathcal{C}$, modulo the short exact sequences. Define $K_\infty(\frkA_\fl)$ to be the quotient group 
\[
    K_\infty(\frkA_\fl) = K_0(\frkA_\fl)/K_0(\g\fmd),
\]
where $\g\fmd$ denotes the category of finite dimensional $\g$-modules, and for every $M\in \frkA_\fl$ denote by $[M]$ the corresponding element in $K_\infty(\frkA_\fl)$.

It remains still to complete the classification of $\frkA^{\irr}$ when considering integral regular central character. However, Theorem~\ref{theo:class-simple-multipli-one-non-integral-reg} suggest that every simple object of $\frkA^{\irr}$ must be submodule of some exponential tensor module. Then we state the conjectures above:
\begin{conj}\label{conj:exponential-mod}
    Let $\bb \in (\C\setminus\{0\})^n$, $\lambda \in P^+$ and $S \subseteq \setn{n}$. For every $0 \le k \le n$, $E(\bb,w_k\cdot \lambda,S)$ has finite length and contains an infinite dimensional simple subquotients $E'(\bb,w_k\cdot \lambda,S)$ satisfying the following equality in $K_\infty(\frkA_\fl)$:
    \begin{align*}
        [E(\bb,\lambda,S)] &= [E'(\bb,w_1\cdot \lambda,S)],\\
        [E(\bb,w_k\cdot \lambda,S)] &= [E'(\bb,w_k\cdot \lambda,S)] + [E'(\bb,w_{k+1}\cdot \lambda,S)], \quad 1 \le k \le n-1;\\
        [E(\bb,w_n\cdot \lambda,S)] &= [E'(\bb,w_n\cdot \lambda,S)].
    \end{align*}
\end{conj}

Suppose that Conjecture \ref{conj:exponential-mod} is true. Then we may state:
\begin{conj}\label{conj:full-classification-of-frkA1}
    For every $\lambda \in P^+$ and $1 \le k \le n$, the set
    \[
        \cE_{k,\lambda} = \left\{E'(\bb,w_k\cdot \lambda,S), E'(\bb,w_{n-k+1}\cdot (-\wlongest\lambda),S)^\tau\mid \bb\in (\C\setminus\{0\})^n,\, S\subseteq \setn{n}\right\}
    \]
    forms a \textbf{complete} list (up to isomorphism) of simple modules of $M \in \frkA^{\irr}$ satisfying
    \[
        \W(M)\sim \cEXT(L(w_k\cdot \lambda)).
    \]
\end{conj}

In this section, we prove some results that support both conjectures. Precisely, we show that
\begin{itemize}
    \item Conjecture \ref{conj:exponential-mod} is true when $S = \emptyset$;
    \item If Conjecture \ref{conj:exponential-mod} is true, then Conjecture \ref{conj:full-classification-of-frkA1} is also true, except for the \emph{completeness} of the lists; and
    \item Exists the subquotients $E'(\bb,w_1\cdot \lambda,S)$ and $E'(\bb,w_n\cdot \lambda,S)$ from Conjecture \ref{conj:exponential-mod} and Conjecture \ref{conj:full-classification-of-frkA1} is true for $k \in \{1,n\}$.
\end{itemize}

As a direct consequence of the last item above, we obtain the complete classification of simple modules of $\frkA^\irr$ for $\g = \lie{sl}(3)$ (Corollary~\ref{cor:full-class-frkA-irr-sl(3)}).

\subsection{Parabolic induced modules and exponential tensor modules}

We show in this section that every exponential tensor module $E(\bb,\lambda,\emptyset)$ is isomorphic to a parabolic Verma module.

Let $\bb=(b_1,\dotsc,b_n) \in (\C\setminus\{0\})^n$ and define 
\[
    \theta_\bb = \exp(-x_\bb) \in \Aut_e(\g),\quad \text{where }x_\bb = \sum_{j=1}^nb_jE_{n+1,j}.
\]
Consider the subalgebras of $\g$
\begin{gather*}
    \mathfrak{l} = \h \oplus \Vspan_\C\{E_{i,j}\mid 1 \le i,j \le n,\, i \ne j\} \simeq \lie{gl}(n),\\
    \mathfrak{u}^+ = \Vspan_\C\{E_{1,n+1},\dotsc,E_{n,n+1}\},\quad \mathfrak{u}^- = \Vspan_\C\{E_{n+1,1},\dotsc,E_{n+1,1}\},\quad \p = \mathfrak{l}\oplus\mathfrak{u}^+\\
    \mathfrak{q}_\bb = \theta_\bb(\p),\quad \h_\bb = \theta_\bb(\h),\quad \mathfrak{l}_\bb = \theta_\bb(\mathfrak{l})\quad \text{and} \quad\mathfrak{u}^{\pm}_\bb = \theta_\bb(\mathfrak{u}^{\pm}).
\end{gather*}
The pullback of $\theta_\bb$ defines a linear isomorphism $\function{\theta_\bb^*}{\h}{\h_\bb^*}$ and we write $\lambda_\bb = \theta_\bb^*(\lambda)$, for every $\lambda \in \h^*$.
For every $m \in \N^n$, define $\nu_m \in \h^*$ to be the weight such that $\nu_m(\tilde h_k) = -m_k$, for every $1 \le k \le n$.

One can check that $\theta_\bb(\tilde h_k) = \tilde h_k - b_kE_{n+1,k}$, for all $1 \le k \le n$. Then, a straightforward computation using the formulas of Section~\ref{sec:tensor-modules} shows that
\begin{lem}[{\citealp[Corollary~5.11]{GN22}}]\label{lem:exp-mod-is-weight-mod}
    For every $\lambda \in P_{\lie{gl}(n)}^+$, the module $E(\bb,\lambda,\emptyset)$ is an $\h_\bb$-weight module.
    Moreover, for every $m \in \N^n$, $\mu \in \Supp L_{\lie{gl}(n)}(\lambda + (n+1)\omega_n)$ and $v_\mu \in L_{\lie{gl}(n)}(\lambda + (n+1)\omega_n)_\mu$, we have that
    \[
        \wtvec_\bb(x^m e^{\bb\bx}\otimes v_\mu) = \nu_\bb + \mu_\bb - (n+1)(\omega_n)_\bb,
    \]
    where $\wtvec_\bb(v) \in \h_\bb^*$ denotes the weight of a weight vector $v$.
\end{lem}

Notice that
\begin{align}
    \label{eq:img-theta_bb-u^-}\theta_\bb(E_{n+1,j}) &= E_{n+1,j},\\
    \label{eq:img-theta_bb-u^+}\theta_\bb(E_{j,n+1}) &= E_{j,n+1} - b_j(\tilde h_1 + \tilde h_2 + \dotsb + \tilde h_n + \tilde h_j) + \sum_{k=1,\,k\ne j}^n b_k E_{j,k} - \sum_{k=1}^n b_k E_{n+1,k},\\
    \label{eq:img-theta_bb-l^+}\theta_\bb(E_{i,i+1}) &= E_{i,i+1} - b_iE_{n+1,i+1}\text{and}\\
    \label{eq:img-theta_bb-l^-}\theta_\bb(E_{i,i+1}) &= E_{i+1,i} - b_{i+1}E_{n+1,i}, 
\end{align}
for all $1 \le j \le n$ and $1 \le i \le n-1$. 

\begin{prop}\label{prop:exp-mod-is-parabolic-verma}
    For every $\lambda \in P_{\lie{gl}(n)}^+$ we have the isomorphism of $\g$-modules 
    \[
        E(\bb,\lambda,\emptyset) \simeq M_{\mathfrak{q}_\bb}(L_{\mathfrak{l}_\bb}(\lambda_\bb)).
    \]
\end{prop}
\begin{proof}
    From the equations~\eqref{eq:img-theta_bb-l^+} and \eqref{eq:img-theta_bb-l^-}, we can compute the following images through $\omega_{\emptyset}$
    \begin{align*}
        \theta_\bb(E_{i,i+1}) &= 1 \otimes E_{i,i+1} - (x_{i+1}(\partial_i-b_i))\otimes 1\quad \text{and}\\
        \theta_\bb(E_{i,i+1}) &=1 \otimes E_{i+1,i} - (x_{i}(\partial_{i+1}-b_{i+1}))\otimes 1, 
    \end{align*}
    for all $1 \le i \le n-1$. In particular, since $\partial_i$ acts as multiplication by $b_i$ on $e^{\bb\bx}$, for all $1 \le i \le n$, the above equations and Lemma~\ref{lem:exp-mod-is-weight-mod} implies that
    \[
        V_0 = \{e^{\bb\bx}\otimes v \mid v \in L_{\lie{gl}(n)}(\lambda + (n+1)\omega_n)\}
    \]
    is an $\mathfrak{l}_\bb$-module isomorphic to $L_{\mathfrak{l}_\bb}(\lambda_\bb)$. Moreover, from direct computation using equation~\eqref{eq:img-theta_bb-u^+} and the definition of the $\g$-module structure on $E(\bb,\lambda,\emptyset)$, we can show that $\theta_\bb(E_{j,n+1})$ acts on $V_0$ by zero, for all $1 \le j \le n$. Hence $\mathfrak{u}^+_\bb V_0 = 0$ and thus $V_0 \simeq L_{\mathfrak{l}_\bb}(\lambda_\bb)$ as $\mathfrak{q}_\bb$-module.
    Then we obtain an isomorphism of $\g$-modules
    \[
         M_{\mathfrak{q}_\bb}(L_{\mathfrak{l}_\bb}(\lambda_\bb)) \simeq \U(\g) \otimes_{\U(\mathfrak{q}_\bb)} V_0.
    \]
    By PBW theorem, we have that $\U(\g) \otimes_{\U(\mathfrak{q}_\bb)} V_0 = \U(\mathfrak{u}^-_\bb)V_0$. And $\U(\mathfrak{u}^-_\bb)V_0 = E(\bb,\lambda,\emptyset)$, since $\theta_\bb(E_{n+1,j})$ acts on it by $-x_j\otimes 1$, for all $1 \le j \le n$. Thus the theorem follows.
\end{proof}

From the previous theorem and Theorem~\ref{theo:deg-c.f.}, it follows that, for any $\lambda \in P_{\lie{gl}(n)}^+$ and $1 \le k\le n$, the modules
\[
    E'(\bb,w_k \cdot \lambda,\emptyset) \coloneqq L_{\mathfrak{q}_\bb}(L_{\mathfrak{l}_\bb}(w_k \cdot \lambda_\bb))
\]
satisfy the conditions of Conjecture~\ref{conj:exponential-mod}. Therefore:
\begin{cor}
    For $S = \emptyset$ the Conjecture~\ref{conj:exponential-mod} is true. 
\end{cor}

\subsection{Supporting results for Conjecture \ref{conj:full-classification-of-frkA1}}

\subsubsection{The list $\cE_{k,\lambda}$}
Let $\admW$ be the full subcategory of $\frkW$ consisting of admissible weight modules. It is easy to see, by definition of admissible weight modules, that $\admW$ is abelian. Moreover, by Lemma~\ref{lem:first-prop-admW}, the objects of $\admW$ have finite length. Thus the abelian group
\[
    K_\infty(\admW) = K_0(\admW)/K_0(\g\fmd)
\]
is well-defined. And given $M\in \admW$, we denote by $[M]$ the corresponding element in $K_\infty(\admW)$. 

Now, notice that Lemmas~\ref{lem:module-almost-equiv-to-semi-simply} and \ref{lem:U-approx} imply that equality in $K_\infty(\admW)$ correspond to almost equivalence in $\admW$, i.e: 
\[
    M \sim N \Leftrightarrow [M] = [N], \quad M,N \in \admW.
\]
In particular, for any two almost-coherent families $\M$ and $\mathcal{N}$, we have:
\begin{equation}\label{eq:sim-equal-to-=}
    \M \sim \mathcal{N} \Leftrightarrow [\M[t]] = [\mathcal{N}[t]], \quad \text{for all }t \in T^*
\end{equation}
\begin{lem}\label{lem:W-induce-a-homo-of-groups}
    For every $t \in T^*$, the weighting functor $\W$ induces a morphism of groups
    \[
        \function{[\W[t]]}{K_\infty(\frkA_\fl)}{K_\infty(\admW)},\quad [M] \mapsto [\W(M)[t]].
    \]
\end{lem}
\begin{proof}
    It follows by induction on the length of objects of $\frkA_\fl$ together with \eqref{eq:sim-equal-to-=} and the fact that, for every exact sequence $0 \rightarrow M \rightarrow N \rightarrow L \rightarrow 0$ in $\frkA$ with $L$ simple, the weighting functors induce an exact sequence
    \[
        \W_1(L) \rightarrow \W(M) \rightarrow \W(N) \rightarrow \W(L) \rightarrow 0,
    \]
    such that $\dim \W_1(L) < \infty$, by Proposition~\ref{prop:finite-dim-Wk}.
\end{proof}

\begin{prop}\label{prop:W(cE_{k,lambda})=EXT(L(w_klambda))}
    Let $\lambda \in P^+$ and $1 \le k \le n$. If Conjecture~\ref{conj:exponential-mod} verifies, then $\W(M) \sim \cEXT(L(w_k\cdot \lambda))$ for every module $M$ from $\cE_{k,\lambda}$.
\end{prop}
\begin{proof}
    Let $\bb \in (\C\setminus\{0\})^n$ and $S\subseteq \setn{n}$.
    Then, for every $\lambda \in P^+$ and $1 \le k \le n$ we have
    \begin{align*}\label{eq:rank-E(bb,wk-dot-lambda,S)-deg}
        \rank E(\bb,w_k\cdot \lambda,S) 
            &= \dim L_{\lie{gl}(n)}(w_k\cdot \lambda + (n+1)\omega_n)\quad 
                &\text{(by Proposition \ref{prop:h-free-mod-rank-one})}\\
            &= \dim L_{\lie{gl}(n)}(w_k\cdot \lambda)\quad 
                &\text{(by Weyl's dimention formula for $\lie{gl}(n)$)}\\
            &= \deg_k(\lambda) + \deg_{k+1}(\lambda)\quad 
                &\text{(by Theorem~\ref{theo:deg-c.f.})},
    \end{align*}
    where $\deg_{n+1} = 0$. 
    Now, with the same type of argument used in the proof of Theorem~\ref{theo:surjectivety-of-W} (with $E(\bb,w_k\cdot \lambda,S)$ in the place of the simple quotients of parabolic Verma modules) we conclude that
    \begin{align*}
        \W(E(\bb,w_k\cdot \lambda,S)) &\sim \cEXT(L(w_k\cdot \lambda)) + \cEXT(L(w_{k+1}\cdot \lambda)),\\
        \W(E(\bb,w_n\cdot \lambda,S)) &\sim \cEXT(L(w_n\cdot \lambda)),
    \end{align*}
    for all $1 \le k \le n-1$ and $\lambda \in P^+$.
        
    Assume that conjecture~\ref{conj:exponential-mod} verifies. It follows from the equation above and Lemma~\ref{lem:W-induce-a-homo-of-groups} that
    \[
        [\W(E'(\bb,w_k\cdot \lambda,S))[t]] =[\cEXT(L(w_k\cdot \lambda))[t]],\quad 1 \le k \le n,
    \]
    for all $t \in T^*$, and the fact \eqref{eq:sim-equal-to-=} implies 
    \[
        \W(E'(\bb,w_k\cdot \lambda,S)) \sim \cEXT(L(w_k\cdot \lambda)).
    \]
    Finally, the almost equivalence 
    \[
        \W(E'(\bb,w_{n-k+1}\cdot (-\wlongest\lambda),S)^\tau) \sim \cEXT(L(w_k\cdot \lambda))
    \]
    follows from the commutation of $(-)^\tau$ with $\W$ and Lemma~\ref{lem:EXT(L(lambda))^tau}.
\end{proof}

We now claim that, for $\lambda \in P^+$, the modules from $\cE_{1,\lambda}$ and $\cE_{n,\lambda}$ exist, independently of Conjecture~\ref{conj:exponential-mod}. Indeed, the modules with $S = \emptyset$ from these lists exist by the previous section. So we need to prove the existence for $S \ne \emptyset$. But notice that, for every $\lambda \in P^+$, we have that $(\lambda, 0) \in \Upsilon$, and since $w_n\cdot 0 = -(n+1)\omega_n$, we also have $(w_n\cdot \lambda,-(n+1)\omega_n) \in \Upsilon$. In particular, $\lambda$ and $w_n \cdot \lambda$ belong to $P_{\lie{gl}(n)}^{+,1}$. Then Prop.~\ref{prop:simplicity-E(b,V_a,S)-non-int-reg} implies the existence of $E'(\bb,w_1\cdot \lambda,S)$ and $E'(\bb,w_n\cdot \lambda,S)$, which are equal to $E(\bb,\lambda,S)$ and $E(\bb,w_n\cdot \lambda,S)$, respectively.

\begin{prop}\label{prop:if-M-covered-by-free-then-classification}
    Let $k \in \{1,n\}$, let $\lambda \in P^+$ and $M \in \frkA^{\irr}$ simple such that $\W(M) \sim \cEXT(L(w_k\cdot \lambda))$. Suppose that there exists a $\U(\h)$-free $\g$-module $\widetilde M \in \frkA$ satisfying
    \begin{itemize}
        \item $M \subseteq \widetilde M$ and $\widetilde M/M$ is finite dimensional;
        \item $\widetilde M$ admits central character and $\chi_{\widetilde M} = \chi_M = \chi_\lambda$.
    \end{itemize}
    Then, up to isomorphism, $M$ belongs to $\cE_{k,\lambda}$
\end{prop}
\begin{proof}
    First, notice that $[M] = [\widetilde M]$ in $K_\infty(\frkA_\fl)$. Then, by Lemma~\ref{lem:W-induce-a-homo-of-groups} we conclude that
    \[
        [\W(\widetilde M)[t]] = [\W(M)[t]] = [\cEXT(L(w_k\cdot \lambda))[t]],
    \]
    for every $t \in T^*$.
    This, together with equation~\eqref{eq:sim-equal-to-=}, implies that $\W(\widetilde M) \sim \cEXT(L(w_k\cdot \lambda))$.
    
    Now, let any $\mu \in P^+$ and consider the translation functor $T_\lambda^\mu$, which is an equivalence of category (since $(\lambda,\mu)\in \Upsilon$). Then, since $T_\lambda^\mu$ commutes with $\W$ and $\cEXT$ and preserves $\sim$, we have
    \[
        \W(T_\lambda^\mu\widetilde M) \sim \cEXT(L(w_k\cdot \mu)).
    \]
    It easy to check that $T_\lambda^\mu$ preserves the category of $\U(\h)$-free modules. So $T_\lambda^\mu\widetilde M$ is $\U(\h)$-free and 
    \begin{equation}\label{eq:deg-rank-tildeM}
        \rank T_\lambda^\mu\widetilde M = \deg_k(\mu),
    \end{equation}
    for all $\mu \in P^+$.
    
    By Theorem~\ref{theo:class-c.f.-sl(n+1)} it follows that $\cEXT(L(w_n\cdot 0)) \sim \cEXT(L(s_1\cdot 0))$. Since $w_1\cdot 0 = s_n \cdot 0 = -\omega_1$ and $s_1\cdot 0 = -2\omega_1 + \omega_2$, we can conclude from Prop.~\ref{prop:c-f-deg-1} that $\cEXT(L(w_n\cdot 0))$ and $\cEXT(L(w_1\cdot 0))$ have degree one. Therefore $T_\lambda^0\widetilde M$ is $\U(\h)$-free of rank one. Then the classification of $\U(\h)$-free modules of rank one (Prop.~\ref{prop:h-free-mod-rank-one}) and Lemma~\ref{lem:central-char-T(P,V)} implies that 
    \[
        T_\lambda^0\widetilde M \in \left\{E(\bb,w_n\cdot 0,S), E(\bb,w_1\cdot 0,S)^\tau, E(\bb,0,S),E(\bb,0,S)^\tau \mid \bb \in (\C\setminus\{0\})^n,\,S \subseteq \setn{n}\right\}.
    \]
    Applying $T_\lambda^\mu T_0^\lambda$ and using Lemma~\ref{lem:trans-tensor-mod}, we obtain that $T_\lambda^\mu\widetilde M \in \widetilde{\cE}_{1,\mu} \cup \widetilde{\cE}_{n,\mu}$, where
    \begin{align*}
        \widetilde{\cE}_{1,\mu} &= \left\{E(\bb,\mu,S),E(\bb,w_n\cdot(-\wlongest\mu),S)^\tau \mid \bb \in (\C\setminus\{0\})^n,\,S \subseteq \setn{n}\right\}\\
        \widetilde{\cE}_{n,\mu} &= \left\{E(\bb,w_n\cdot \mu,S), E(\bb,w_1\cdot (-\wlongest\mu),S)^\tau \mid \bb \in (\C\setminus\{0\})^n,\,S \subseteq \setn{n}\right\}.
    \end{align*}
    From Proposition~\ref{prop:h-free-mod-rank-one}, the fact that $\tau(\h) = \h$ and equation~\eqref{eq:def-deg_w_k}, we have that the element of $\widetilde{\cE}_{n,\mu}$ and $\widetilde{\cE}_{1,\mu}$ have rank equal to $\deg_n$ and $\dim L_{\mathfrak{l}}(\mu)$, respectively. And $\dim L_{\mathfrak{l}}(\mu)=\deg_1(\mu)$, by Remark~\ref{rem:formula-deg_1}. Thus $\widetilde M \in \widetilde{\cE}_{k,\lambda}$, by linear independence of $\deg_1,\dotsc,\deg_n$ (by Lemma~\ref{lem:deg_w_i's-linear-independent}). Finally, since $\cE_{k,\lambda}$ is exactly the set of those unique infinite dimensional simple subquotient of modules of $\widetilde{\cE}_{k,\lambda}$, we conclude that $M$ must be isomorphic to one module of $\cE_{k,\lambda}$.
\end{proof}

\subsubsection{Duality in $\frkA$}

This subsection is motivated by Prop.~\ref{prop:if-M-covered-by-free-then-classification}. More precisely, the goal is to create a tool to find $\U(\h)$-free hulls of simple modules of $\frkA$. The idea is to use $\U(\h)$-reflexive hulls.

For $M \in \frkA$, define the $\U(\h)$-module
\[
    M^\vee = \Hom_{\U(\h)}(M,\U(\h)).
\]
Since $M$ is $\U(\h)$-finitely generated, so is $M^\vee$. Consider the anti-automorphism $\function{\hat \tau = -\tau}{\U(\g)}{\U(\g)}$ and the ring automorphisms $\function{\sigma_\lambda}{\U(\h)}{\U(\h)}$, for all $\lambda \in \h^*$, defined in \ref{def:sigma_lambda}.

\begin{prop}\label{prop:defin-duality}
    The module $M^\vee$ becomes a $\U(\g)$-module by define the action
    \[
        \function{x_\alpha \cdot F}{M}{\U(\h)},\quad m \mapsto \sigma_\alpha(F(\hat\tau(x_\alpha)m)),
    \]
    for all $F \in M^\vee$, $\alpha \in \Delta$ and root-vector $x_\alpha \in \g_\alpha$. Moreover, the assignments
    \[
        M \mapsto M^\vee \quad \text{and} \quad \Hom_\mathfrak{A}(M,N)\ni f \mapsto f^\vee = \Hom_{\U(\h)}(f,\U(\h)) \in \Hom_{\frkA}(N^\vee,M^\vee)
    \]
    defines an additive contra-variant functor $\function{(-)^\vee}{\frkA}{\frkA}$.
\end{prop}
\begin{proof}
    The proof is a standard verification of the definitions and relations, using equation \eqref{eq:dif-op} and the fact that $\hat\tau(\g_\alpha) = \g_{-\alpha}$, for all $\alpha \in \Delta$.
\end{proof}

\begin{lem}\label{lem:propriedade-duality}
    The functor $(-)^\vee$ satisfies the following properties:
    \begin{enumerate}[label=(\roman*)]
        \item Preserves central character;
        \item Commutes with the functor $-\otimes V$, for every finite dimensional $\g$-module $V$; and
        \item Commutes with translation functor.
    \end{enumerate}
\end{lem}
\begin{proof}
    The first item follows from the fact that $\hat\tau$ fix $Z(\g)$ pointwise (see \cite[Exercise 1.10]{Hum08}).

    For the second item, consider the duality $\function{(-)^*}{\cO}{\cO}$ in the BGG category $\cO$ (see \cite[Section~3.2]{Hum08}). It is well-known that $L^* \simeq L$, for any simple object of $\cO$. So, for $V \in \g\md$, fix an isomorphism of $\g$-modules $\function{\varphi}{V^*}{V}$. Given $M \in \frkA$, define $\function{\Phi_M}{M^\vee\otimes V}{(M\otimes V)^\vee}$ by
    \[
        \Phi_M(F\otimes v_\mu)\colon M\otimes V \ni m\otimes v\mapsto \varphi^{-1}(v_\mu)\sigma_\mu(F(\mu)),
    \]
    for every $F \in M^\vee$ and weight vector $v_\mu \in V_\mu$. One can check that $\Phi_M(F\otimes v_\mu)$ is a homomorphism of $\U(\h)$-module and that $\Phi_M$ is a homomorphism of $\g$-modules. Now, fix a basis $B \subseteq V$ consisting of weight vectors such that $\varphi^{-1}(v)(w) = \delta_{v,w}$, for all $v,w \in V$. For each $\widetilde F \in (M\otimes V)^\vee$ and $v_\mu \in B$ of weight $\mu$, define $\function{\widetilde F_{v_\mu}}{M}{\U(\h)}$ by setting $\widetilde F_{v_\mu}(m) = \sigma_\mu^{-1}(\widetilde F(m\otimes v_\mu))$, for all $m \in M$, and define 
    \[
        \function{\Psi_M}{(M\otimes V)^\vee}{M^\vee\otimes V},\quad 
        \widetilde F \mapsto \sum_{v_\mu\in B}\widetilde F_{v_\mu} \otimes v_\mu.
    \]
    Again, one can check that $\widetilde F_{v_\mu}$ is indeed a homomorphism of $\U(\h)$-module.
    Now, notice that
    \[
        (\Phi_M(F\otimes v_\mu))_{w_\lambda}(m) = \sigma_\mu^{-1}\varphi^{-1}(\Phi_M(F\otimes v_\mu)(m\otimes w_\lambda)) = \sigma_\mu^{-1}(\varphi^{-1}(v_\mu)(w_\lambda)\sigma_\mu(F(m))) = \delta_{v_\mu,w_\lambda}F(m),
    \]
    for all $F\in M^\vee$, $v_\mu,w_\lambda \in B$ and $m \in M$. Then 
    \[
    \Psi_M(\Phi_M(F\otimes v_\mu)) = \sum_{w_\lambda\in B}\widetilde (\Phi_M(F\otimes v_\mu))_{w_\lambda} \otimes w_\lambda = F\otimes v_\mu
    \]
    and $\Psi_M \circ \Phi_M = \Id_{M^\vee \otimes V}$. Analogously, it shows that $\Phi_M\circ \Psi_M = \Id_{(M\otimes V)^\vee}$. It also can be shown that these maps are natural in $M$, and hence the second item follows.

    Finally, the last item follows from the two previous one, together with additiveness of $(-)^\vee$.
\end{proof}

If $M \in \frkA$ is simple of infinite dimension, it is $\U(\h)$-torsion free by Corollary~\ref{cor:simple-tor-free}. In particular, the map
\[
    \function{\ev}{M}{M^{\vee\vee}},\quad m\mapsto \left(\ev_m\colon M^\vee \ni F\mapsto F(m)\right),
\]
is a monomorphism of $\U(\h)$-module. We may prove that it is in fact a homomorphism of $\U(\g)$-module. Indeed, for every weight vector $x_\alpha \in \g_\alpha$, $F \in M^\vee$ and $m \in M$ we have
\[
    (x_\alpha \cdot \ev_m)(F) = \sigma_\alpha(\ev_m(\hat\tau(x_\alpha)\cdot F)) = \sigma_\alpha((\hat\tau(x_\alpha)\cdot F)(m)) =\sigma_\alpha(\sigma_\alpha^{-1}(F(x_\alpha \cdot m))) =\ev_{x_\alpha \cdot m}(F).
\]

\begin{lem}\label{lem:ev-finite-dim-coker}
    For every simple $\U(\h)$-finite module $M$ of infinite dimension, the injective homomorphism $\function{\ev}{M}{M^{\vee\vee}}$ has cokernel of finite dimension.
\end{lem}
\begin{proof}
    Identify $M$ with $\ev(M)$. Since $M^{\vee\vee}_{\m} = \Hom_{\U(\h)_\m}(\Hom_{\U(\h)_\m}(M_{\m},\U(\h)_\m),\U(\h)_\m)$, for all $\m \in \Specm \U(\h)$, it follows that $M^{\vee\vee}_{\m} =M_\m$, whenever $M_\m$ is $\U(\h)_\m$-free. So, by Theorem \ref{theo:almost-free}, we conclude that $(M^{\vee\vee}/M)_{\m_\lambda} = 0$ for all $\lambda \in \h^* \setminus \Supp(\W_1(M))$. But $\Supp(\W_1(M))$ is finite by Prop.~\ref{prop:finite-dim-Wk}. Thus $\Suppca M^{\vee\vee}/M\cap \Specm \U(\h)$ is finite and the proposition follows from Lemma~ \ref{lem:finite-support-finite-dim}.
\end{proof}

Beck to our motivation, we state the following:
\begin{lem}\label{lem:rank-1-implies-free}
    Let $M \in \frkA$ simple of infinite dimension. If $M$ can be translated to a module of rank 1, then $M^{\vee\vee}$ is $\U(\h)$-free.
\end{lem}
\begin{proof}
    Let $(\lambda,\mu) \in \Upsilon$ such that the central character of $M$ is $\chi_\lambda$ and $T_\lambda^\mu M$ has rank one. By \cite[Prop. 1.9]{Har80}, $(T_\lambda^\mu M)^{\vee\vee}$ is $\U(\h)$-free. Since $(-)^\vee$ commutes with translation functors, by Lemma~\ref{lem:propriedade-duality}, we have that $T_\lambda^\mu (M^{\vee\vee})$ is $\U(\h)$-free. Therefore, since translation functors preserve the category of $\U(\h)$-free modules, the bidual $M^{\vee\vee}$ is $\U(\h)$-free.
\end{proof}

\begin{theo}\label{theo:class-frkA^(1)-integ-reg-k=1,n}
    Conjecture~\ref{conj:full-classification-of-frkA1} holds for $k \in \{1,n\}$. 
\end{theo}
\begin{proof}
    Let $\lambda \in P^+$ and let $M \in \frkA^{\irr}$ simple of infinite dimension such that $\W(M)\sim \cEXT(L(w_k\cdot \lambda))$, where $k = 1$ or $n$. Then, by applying the translation functor $T_\lambda^0$ and noticing that $\deg_n(0) = \deg_1(0) = 1$, we conclude that $\rank(T_\lambda^0 M) = 1$. In particular, $M$ satisfies the conditions of Lemma~\ref{lem:rank-1-implies-free}. Thus $M^{\vee\vee}$ is $\U(\h)$-free. Moreover
    \begin{itemize}
        \item $M^{\vee\vee}/M$ is finite dimensional, by Lemma~\ref{lem:ev-finite-dim-coker}; and
        \item $M^{\vee\vee}$ has the same central character of $M$, since $(-)^{\vee\vee}$ preserves central character by Lemma~\ref{lem:propriedade-duality}.
    \end{itemize}
    Thus we may apply Proposition~\ref{prop:if-M-covered-by-free-then-classification}, obtaining the classification.    
\end{proof}

Since, for $n=2$, the almost-coherent families obtained after applying $\W$ to the modules of
\[
    \cE_{1,\lambda} \cup \cE_{n,\lambda}, \quad \lambda \in P^+
\]
correspond (under almost equivalence) with all possibles irreducible semi-simple coherent families of integral regular central character, we see that Theorems~\ref{theo:class-frkA^(1)-integ-reg-k=1,n} and \ref{theo:class-simple-multipli-one-non-integral-reg} complete the classification for $\lie{sl}(3)$:

\begin{cor}\label{cor:full-class-frkA-irr-sl(3)}
    For $\g= \lie{sl}(3)$, every simple module of $\frkA^\irr$ is of the form $E(\bb,\mu,S)$, $E'(\bb,w_1\cdot\lambda,S)$, $E'(\bb,w_2\cdot\lambda,S)$, or a $\tau$-twist of them (for $\bb \in (\C\setminus \{0\})^n$, $S\subseteq \{1,2\}$, $\mu \in \h^*$ non-integral such that $\chi_\mu \in \X^\infty(A_2)$ and $\lambda \in P^+$ is integral dominant).
\end{cor}
%\newpage

%%%%%%%%%%%%%%%%%%%%%%%%%%%%%%%%%%%%%%%%%%%%%%%%%%%%%%%%%%%%%%%%%%%
% Further directions
%%%%%%%%%%%%%%%%%%%%%%%%%%%%%%%%%%%%%%%%%%%%%%%%%%%%%%%%%%%%%%%%%%%
%\section*{Further directions}
%\import{sections_paper/}{introduction}
%\newpage

%%%%%%%%%%%%%%%%%%%%%%%%%%%%%%%%%%%%%%%%%%%%%%%%%%%%%%%%%%%%%%%%%%%
% References
%%%%%%%%%%%%%%%%%%%%%%%%%%%%%%%%%%%%%%%%%%%%%%%%%%%%%%%%%%%%%%%%%%%
%\bibliographystyle{plain}
%\bibliography{ref.bib}

\def\cprime{$'$} \newcommand{\doi}[1]{\href{http://dx.doi.org/#1}{\tt \nolinkurl{{http://dx.doi.org/#1}}}}

%\nocite{*}

\Addresses

\end{document}